\documentclass{article}
\usepackage[english]{babel}
\usepackage[utf8]{inputenc}

\usepackage{bm}
\usepackage{amsmath}
\usepackage{graphicx}
\usepackage{color}
\usepackage{subcaption}
\usepackage{amsfonts}
\newtheorem{remark}{Remark}
\graphicspath{ {./figs/} }

\begin{document}
\title{Multiscale dimension reduction for flow and transport problems in thin domain with reactive boundaries}

\author{
Maria Vasilyeva \thanks{Multiscale model reduction Laboratory, North-Eastern Federal University, Yakutsk, Russia, Institute for Scientific Computation, Texas A\&M University, College Station, TX,  USA. 
Email: {\tt vasilyevadotmdotv@gmail.com}.} 
\and
Valentin Alekseev \thanks{Yakutsk branch of the Regional Scientific and Educational Mathematical Center "Far Eastern Center of Mathematical Research", North-Eastern Federal University, Yakutsk, Russia, Multiscale model reduction Laboratory, North-Eastern Federal University, Yakutsk, Russia} 
\and
Eric T. Chung 
\thanks{Department of Mathematics,
The Chinese University of Hong Kong (CUHK), Hong Kong SAR. 
Email: {\tt tschung@math.cuhk.edu.hk}.}
\and
Yalchin Efendiev \thanks{Department of Mathematics \& Institute for Scientific Computation (ISC), Texas A\&M University, College Station, Texas, USA. Email: {\tt efendiev@math.tamu.edu}.} 
}

\maketitle

\begin{abstract}
In this paper, we consider flow and transport problems in thin domains. 
Modeling problems in thin domains occur in many applications, where
thin domains lead to some type of reduced models. A typical example is
one dimensional reduced-order model for flows in pipe-like geometries 
(e.g., blood vessels). In many reduced-order models, the model equations
are described apriori by some analytical approaches. In this paper,
we propose the use of multiscale methods, which are alternative to 
reduced-order models and can represent reduced-dimension modeling by using 
fewer basis functions (e.g., the use of one basis function corresponds
to one dimensional approximation).

The mathematical model considered in the paper
 is described by a system of equations for velocity, 
pressure, and concentration, where the flow is described by the 
Stokes equations, and the transport is described by an unsteady 
convection-diffusion equation with non-homogeneous boundary 
conditions on walls (reactive boundaries).
We start with the finite element approximation of the problem on unstructured grids and use it as a reference solution for two and three-dimensional model problems.
Fine grid approximation resolves complex geometries on the grid level and leads to a large discrete system of 
equations that is computationally expensive to solve.
To reduce the size of the discrete systems, we develop a multiscale model reduction technique, where we construct local multiscale basis functions to generate a  lower-dimensional model on a coarse grid.
The proposed multiscale model reduction is based on the Discontinuous Galerkin Generalized Multiscale Finite Element Method (DG-GMsGEM).
In DG-GMsFEM for flow problems, we start with constructing the snapshot space 
for each interface between coarse grid cells to capture possible flows. 
For the reduction of the snapshot space size, we perform a dimension reduction via a solution of the local spectral problem and use eigenvectors corresponding to the smallest eigenvalues as multiscale basis functions for the approximation on the coarse grid.
For the 
 transport problem, we construct multiscale basis functions for each interface between coarse grid cells and present additional basis functions to capture non-homogeneous boundary conditions on walls.
Finally, we will present some numerical simulations for three test geometries for two and three-dimensional problems to demonstrate the method's performance. 
\end{abstract}

\section{Introduction}

Mathematical models in thin domains occur in many real-world applications, 
scientific and engineering problems. 
Fluid flow and transport in thin tube structures are widely used in biological applications, for example, to simulate blood flow in vessels  \cite{quarteroni2004mathematical, nachit2013asymptotic, oshima2001finite}.
In engineering problems, flow simulation is used to study fluid flow in complex pipe structures, for example in pipewise industrial installations, wells in oil and gas industry, heat exchangers, etc. 
In reservoir simulations, thin domains are related to the fractures that usually have complex geometries with very small thickness compared to typical reservoir sizes \cite{martin2005modeling, formaggia2014reduced, Quarteroni2008coupling}.
Such problems are often 
considered with complex interaction processes with surrounding media or 
walls. In many applications, these problems
are transformed to reduced (e.g., one) dimensional problems via some type of apriori postulated
models. Our goal is to present an alternative approach to analytical model reduction by using multiscale
basis functions.

For the applicability of the convenient numerical methods for simulations of such problems, a very fine grid should be constructed to resolve the geometry's complex structure  on the grid level. 
Moreover, a very small domain thickness provides an additional complexity in the grid construction for thin and long domains. 
For a fast and accurate solution of the presented problem, a homogenization (upscaling) technique or multiscale models are used, which are based on 
constructing the lower dimensional coarse grid models. 
The asymptotic behavior of the solutions in thin domains is intensively studied.
In \cite{panasenko2012asymptotic}, the authors consider the problem of complete asymptotic expansion for the time-dependent Poiseuille flow in a thin tube. In \cite{panasenko2017method}, the method of asymptotic partial decomposition of a domain (MAPDD) was presented to reduce the computational complexity of the numerical solution of such problems. This method combines the three-dimensional description in some neighborhoods of bifurcations and the one-dimensional description out of these small subdomains, and it prescribes some special junction conditions at the interface between these 3D and 1D submodels. Numerical results were presented in \cite{nachit2013asymptotic} for the method of asymptotic partial domain decomposition for thin tube structures with Newtonian and non-Newtonian flows in large systems of vessels. Our goal is to provide an alternative systematic approach to model reduction
for thin domain problems that can add complexity via additional multiscale basis functions.

Model reduction techniques usually depend on a coarse grid approximation, which can be obtained by 
discretizing the problem on a coarse grid and choosing a suitable coarse-grid formulation of the problem.
In the literature, several approaches have been developed to obtain the coarse-grid formulation for the problems in heterogeneous domains, including multiscale finite element method \cite{hou1997multiscale, efendiev2009multiscale,chung2010reduced},
mixed multiscale finite element method \cite{chen2003mixed, aarnes2004use},
generalized multiscale finite element method \cite{efendiev2013generalized,chung2016adaptive,chung2013sub,cheung2020constraint,chung2020convergence},
multiscale mortar mixed finite element method \cite{arbogast2007multiscale},
multiscale finite volume method \cite{jenny2005adaptive, hajibeygi2008iterative},
variational multiscale methods \cite{hughes1998variational,kornhuber2018analysis,maalqvist2014localization},
and heterogeneous multiscale methods \cite{abdulle2012heterogeneous} 
and etc.
The non - conforming multiscale method is considered for the solution of the Stokes flow problem in a heterogeneous domain in \cite{muljadi2015nonconforming, jankowiak2018non}.
In \cite{chung2016generalized, chung2018multiscale}, we presented the Generalized Multiscale Finite Element Method (GMsFEM) for the solution of the flow problems in perforated domains with continuous multiscale basis functions. GMsFEM shows a good accuracy for solving the nonlinear (non - Newtonian) fluid flow problems \cite{chung2016generalizednn}.
In \cite{chung2017conservative}, we presented the Discontinuous Galerkin Generalized Multiscale Finite Element Method (DG-GMsFEM)  for the solving the 
two - dimensional problems in perforated domains.

In this paper, we consider flow and transport processes in thin structures with reactive boundaries. The mathematical model is described using Stokes equations and the unsteady 
convection-diffusion equation with non - homogeneous boundary conditions. 
Non-homogeneous boundary conditions occur in many applications.  For example, the pore-scale modeling and simulation of reactive flow in porous media have many applications in many branches of science such as biology, physics, chemistry, geomechanics, and geology \cite{iliev2017pore, hornung1991diffusion, allaire2010homogenization}. 
In order to handle the complex geometry of the walls and non-homogeneous boundary conditions on them, we present additional spectral multiscale basis functions. 
In this work, we continue developing the multiscale model reduction techniques for problems with multiscale features and developing the generalization of the techniques for DG-GMsFEM. 
In our previous work \cite{vasilyeva2019upscaling, spiridonov2019generalized}, we considered elliptic problems in perforated domains and constructed additional basis functions to capture non-homogeneous boundary conditions on perforations. In \cite{vasilyeva2019upscaling}, we proposed a non-local multi-continua (NLMC) method for Laplace, elasticity, and parabolic equations with non-homogeneous boundary conditions on perforations. 
In \cite{spiridonov2019generalized}, we considered the Continuous Galerkin Generalized Multiscale Finite Element Method (CG-GMsFEM) for problems in perforated domains with non-homogeneous boundary conditions, where we constructed one additional basis function for local domains with perforations. Recently, we extended this technique for solving unsaturated flow problems in heterogeneous domains with rough boundary \cite{spiridonov2020multiscale}.  
In this work, we consider DG-GMsFEM and present additional spectral basis functions for rough non-homogeneous boundaries in transport problems. The concept is based on the separation of the snapshots for each feature and shares a lot of similarities with multiscale methods for fractured
media presented in our previous works \cite{chung2018non, vasilyeva2019nonlocal}. Our work is also motivated by a recently developed Edge GMsFEM, where multiscale basis functions are constructed for each coarse grid interface \cite{fu2019wavelet, fu2019edge}.

In this work, we use the DG-GMsFEM for constructing multiscale basis functions for problems in thin domains. 
The goals of this paper are in constructing the general approach for problems in complex thin geometries with an accurate approximation of the velocity space and transport processes. 
We construct local multiscale basis functions to generate a  lower-dimensional model on a coarse grid. 
In DG-GMsFEM for flow problem, we start with constructing the snapshot space that captures possible flows between coarse cell interfaces. 
After constructing the snapshot space, we perform a dimension reduction by a solution of the local spectral problems.  
For the pressure approximation, we use piecewise constant functions on the coarse grid. 
For the transport problem, we construct multiscale basis functions for each interface between coarse grid cells and present additional basis functions to capture non-homogeneous boundary conditions on walls. 
The presented snapshot spaces can accurately capture processes on the rough wall boundaries with non-homogeneous boundary conditions on them. 
We present the results of the numerical simulations for three test geometries for two and three-dimensional problems.

As we mentioned earlier, our goal is to investigate the use of multiscale and generalized multiscale
methods for dimension reduction for problems in this domains. Many existing approaches propose
analytical or semi-analytical reduced-dimensional problems, where the dimensions are determined
apriori. Our idea is to use generalized multiscale method and identify the dimension. The proposed approach
combined with aposteriori error estimate can further identify the dimension across thin layers and allow
obtaining  more accurate solution.

The paper is organized as follows. In Section 2, we describe the problem formulation and the fine-scale approximation. In Section 3, we present the multiscale method for flow and transport processes in thin structures with rough reactive boundaries.
In Section 4, we present numerical results. The paper ends with a conclusion.

\section{Problem formulation} \label{s-mmfg}

\begin{figure}[h!]
\centering
\includegraphics[width=0.9\textwidth]{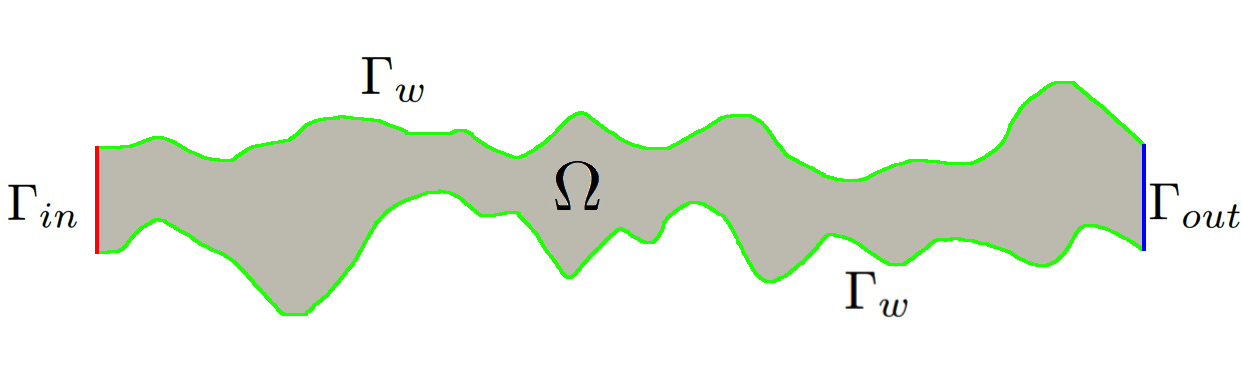}
\caption{Illustration of a thin domain $\Omega$ with multiscale features. }
\label{domain}
\end{figure}

Let $\Omega$ be a thin domain with multiscale features, 
where the thickness is small compared to the domain length (See Figure \ref{domain} for an illustration).  
We will consider the following flow and transport equations in the thin domain $\Omega$
\begin{equation}
\label{eq:mm1}
\begin{split}
\rho \frac{\partial u}{\partial t} - \mu \Delta u + \nabla p = 0 , \quad 
x \in \Omega,  \\
\nabla \cdot  u = 0, \quad 
x \in \Omega,\\
\frac{\partial c}{\partial t} + u \nabla c - \nabla (D \nabla c)  = 0, \quad 
x \in \Omega,
\end{split}
\end{equation}
where
$\mu$ is the fluid viscosity,
$\rho$ is the fluid density,
$D$ is the diffusion coefficient,
$c$ is the concentration,
$u$ and $p$ are the velocity and pressure.

The system  \eqref{eq:mm1} is equipped with following initial conditions
\[
c = c_0, \quad 
u = u_0, \quad x \in \Omega, \quad t = 0.
\]
Moreover, 
We consider the following boundary conditions for the flow problem
\[
\begin{split}
u = g, \quad &x \in \Gamma_{in}, \\ 
(\nabla u - p \mathcal{I}) \cdot n = 0, \quad &x \in \Gamma_{out},\\
u = 0, \quad &x \in \Gamma_w, 
\end{split}
\]
and for the transport problem
\[
\begin{split}
c = c_{in}, \quad &x \in \Gamma_{in}, \\ 
- D \nabla c  \cdot n = 0, \quad &x \in \Gamma_{out},\\
- D \nabla c  \cdot n = \alpha (c - c_w), \quad &x \in \Gamma_w, 
\end{split}
\]
where $n$ is the unit outward normal vector on $\partial \Omega$, 
$\mathcal{I}$ is the  identity matrix, 
$\Gamma_{in}$ is the inflow boundary, 
$\Gamma_{out}$ is the outflow boundary, 
$\Gamma_w$ is the reactive boundary of the thin domain, 
$\Gamma_w \cup \Gamma_{in} \cup \Gamma_{out} = \partial \Omega$ (see Figure \ref{domain}).

\begin{figure}[h!]
\centering
\includegraphics[width=0.9\textwidth]{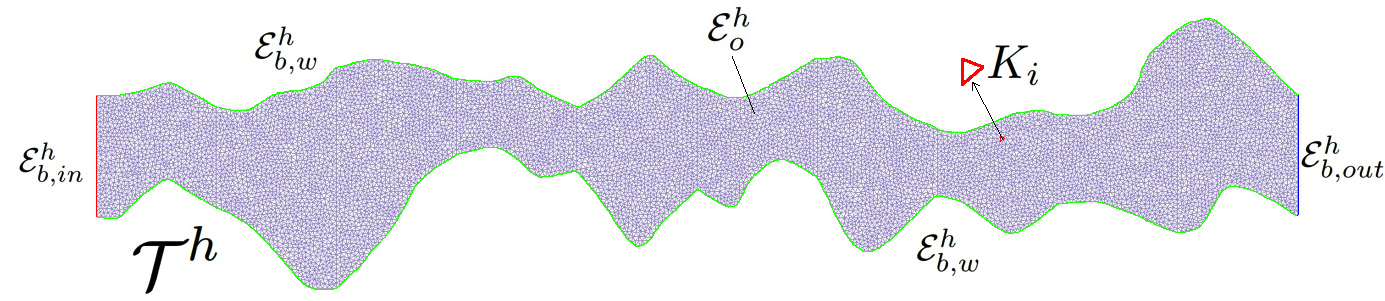}
\caption{Illustration of the fine grid for domain $\Omega$}
\label{gfine}
\end{figure}

In order to solve the problem \eqref{eq:mm1}, we generate an unstructured grid (fine grid) and use a finite element approximation for the spatial discretization. 
Let 
$\mathcal{T}^h$ be a fine-grid partition of the domain $\Omega$ given by
\[
\mathcal{T}^h = \bigcup_{i = 1}^{N_{cell}^h} K_i, 
\]
where $N_{cell}^h$ is the number of fine grid cells. We use
$\mathcal{E}^h$ to denote the set of facets in $\mathcal{T}^h$ with  
$\mathcal{E}^h = \mathcal{E}^h_{o} \cup \mathcal{E}^h_{b}$, 
where 
$\mathcal{E}^h_{o}$ is the set of interior facets and
$\mathcal{E}^h_{b}$ is the set of boundary facets with
$\mathcal{E}^h_{b} = \mathcal{E}^h_{b,in} \cup \mathcal{E}^h_{b, w} \cup \mathcal{E}^h_{b, out}$ (see Figure \ref{gfine}).
We use the notations $K$ and $E$ to denote a generic cell and facet in the fine grid $\mathcal{T}^h$ (see Figure \ref{gfine}). We define the jump $[u]$ and the average $\{u\}$ of a function $u$ on interior facet by
\[
[u] = u_+ - u_-, \quad 
\{u\} = \frac{u_+ + u_-}{2},
\]
where $u_+ = u|_{K^{+}}$, $u_- =  u|_{K^{-}}$, $K^+$ and $K^-$ are the two cells sharing the facet $E$. 
Note that, we define $[u] = u|_E$ and  $\{u\} = u|_E$ for $E \in \mathcal{E}^h_{b}$.

We define the fine-scale velocity space
\[
V_h = \{
v \in L^2(\Omega): \,
v|_{K} \in (\mathbb{P}_1(K))^2, \, \forall K \in \mathcal{T}^h\},
\]
which contains functions that are piecewise linear in each fine-grid element $K$ and are discontinuous across coarse grid edges.  
For the pressure, we use the space of piecewise constant functions
\[
Q_h = \{ q \in L^2(\Omega) : \,  q |_K \in \mathbb{P}_0(K), \, \forall K \in \mathcal{T}^h \}.
\]  
The  fine-scale  space for concentration is the following
\[
P_h =\{ 
v \in L^2(\Omega): \, 
v|_{K} \in \mathbb{P}_1(K), \, \forall K \in \mathcal{T}^h \}.
\]
Using these spaces, we
have the following IPDG variational formulation with implicit time approximation for the approximation of \eqref{eq:mm1}:
\begin{itemize}

\item Flow problem: find $(u_h, p_h) \in V_h \times  Q_h$ such that
\begin{equation}
\label{eq:flow}
\begin{split}
\frac{1}{\tau }m^u(u_h - \check{u}_h, v) 
+ a^u_{\text{DG}}(u_h, v) 
+ b_{\text{DG}}(p_h, v) &= l^u(v), \quad \forall v \in V_h, \\
b_{\text{DG}}(u_h, q) &= l^p(q), \quad \forall q \in Q_h,  
\end{split}
\end{equation} 
where
\[
\begin{split}
a^u_{\text{DG}}(u,v) 
&= \sum_{K \in \mathcal{T}^h} \int_K
\mu \nabla u \cdot \nabla v \, dx \\
&- \sum_{E \in  \mathcal{E}^h / \mathcal{E}^h_{b,out}} \int_E \Big( 
\{ \mu \nabla u \cdot n \}\cdot [v] + 
\{ \mu \nabla v \cdot  n\}  \cdot [u] -  
\frac{\gamma_u}{h} \{\mu\} [u]\cdot [v] \Big) \, ds,
\\
m^u(u,v)  &= \sum_{K \in \mathcal{T}^h} 
 \int_K \rho \, u \, v \, dx, 
 \\
b_{\text{DG}}(u, p) &= - \sum_{K \in \mathcal{T}^h}  
\int_K p \, \nabla u \, dx  +  
\sum_{E \in \mathcal{E}^h / \mathcal{E}^h_{b,out}} \int_E p \, [u] \cdot n \, ds,
\\
l^u(v) &=  \sum_{E \in \mathcal{E}^h_{b, in}} \int_E 
\left( \frac{\gamma_u}{h} \, \mu  v -  (\mu \nabla v \cdot n) \right) \cdot g  \, ds, 
\\
l^p(q) &= \sum_{E \in \mathcal{E}^h_{b, in}} \int_E  (g \cdot n) \, q \, dx,
\end{split}
\]
and $\gamma_u$ is the penalty perm. 

\item  Transport problem: find $c_h \in P_h$ such that
\begin{equation}
\label{eq:transport}
\begin{split}
\frac{1}{\tau} m^c(c_h - \check{c}_h , r) 
+ c^c_{\text{DG}}(c_h, r) 
+ a^c_{\text{DG}}(c_h, r) &= l^c(r), 
\quad \forall r \in P_h,
\end{split}
\end{equation} 
where
\[
\begin{split}
a^c_{DG}(c, r) 
& = \sum_{K \in \mathcal{T}^h} 
\int_K D \nabla c \cdot \nabla r \, dx\\
&-\sum_{E \in \mathcal{E}^h_0 \cup \mathcal{E}^h_{b,in}} \int_E \Big(
 \{ D \nabla c \cdot n \} [r] + \{ D \nabla r \cdot n \} [c] - 
\frac{\gamma_c}{h} \{D\} [c][r] \Big) \, ds\\
& + 
\sum_{E \in \mathcal{E}^h_{b, w}} \int_E \alpha \, c \, r \, ds,
\\
m^c(c, r) &= \sum_{K \in \mathcal{T}^h} 
\int_K c \, r \, dx, 
\\ 
c^c_{DG}(c, r)  &=  \sum_{K \in \mathcal{T}^h} 
\int_K (u_h \, c) \cdot \nabla r \, dx + 
\sum_{E \in \mathcal{E}^h_0 \cup \mathcal{E}^h_{b,out}} 
\int_E (\tilde{u}_+ c_+ - \tilde{u}_- c_-) \, [r] \, ds,
\\
l^c(r) &=  \sum_{E \in \mathcal{E}^h_{b, in}} \int_E \left( 
\frac{\gamma_c}{h} \, D  \, r  -  \{ D \nabla r \cdot n \}  
\right) \, c_{in} \, dx 
+ 
\sum_{E \in \mathcal{E}^h_{w}} \int_E \alpha \, c_w \, r \, ds,
\end{split}
\]
and $\gamma_c$ is the penalty perm and 
$\tilde{u} = ( u_h \cdot n + |u_h \cdot n| )/2$. 
\end{itemize}
In the above systems (\ref{eq:flow}) and (\ref{eq:transport}), $\check{u}_h$ and $\check{c}_h$ are the solutions from the previous time step and $\tau$ is the given time step size.

We can write the above discrete systems (\ref{eq:flow}) and (\ref{eq:transport}) as follows. 
\begin{itemize}
\item Flow problem:
\begin{equation}
\begin{split}
\frac{1}{\tau}
\begin{pmatrix}
M^u_h & 0 \\
0 & 0
\end{pmatrix}
\begin{pmatrix}
u_h - \check{u}_h  \\
p_h - \check{p}_h
\end{pmatrix} + 
\begin{pmatrix}
A^u_h & B_h^T \\
B_h & 0
\end{pmatrix}
\begin{pmatrix}
u_h \\
p_h
\end{pmatrix} = 
\begin{pmatrix}
F_h^u \\
F_h^p
\end{pmatrix}.
\end{split}
\end{equation} 

\item Transport problem:
\begin{equation}
\frac{1}{\tau} M^c_h  (c_h - \check{c}_h) + (A^c_h + C^c_h(u_h)) c_h =  F^c_h.
\end{equation} 
Here, the matrix $C^c_h(u_h)$ depends on the function $u_h$.
\end{itemize}
In the next section, we will present the proposed multiscale method for the solution of the flow and transport problems that used to reduce the system size. 
In the multiscale method, we solve problems in local domains with various boundary conditions to form a snapshot space and use a spectral problem in the snapshot space
to perform the required dimension reduction.

\section{Multiscale method} \label{s-cg}

\begin{figure}[h!]
\centering
\includegraphics[width=0.9\textwidth]{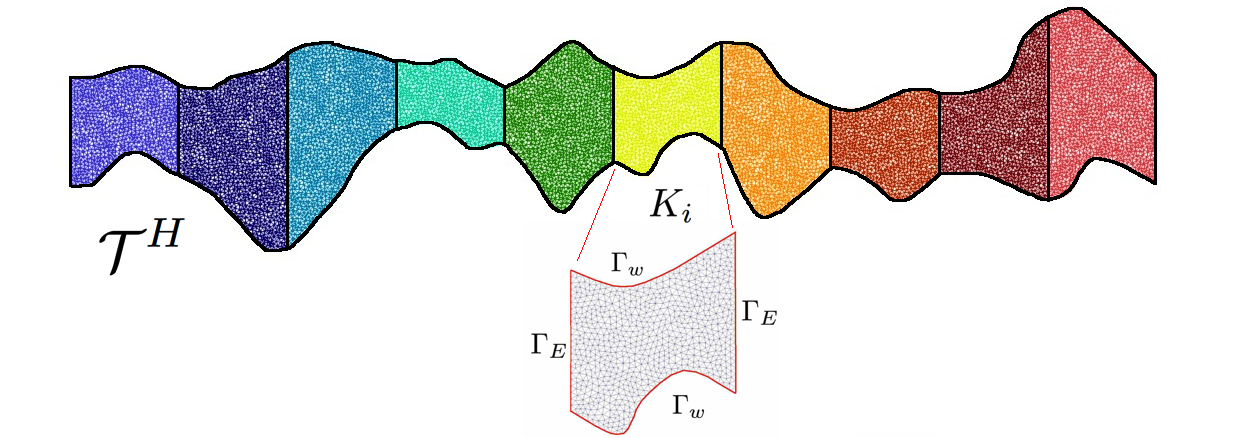}
\caption{Illustration of the coarse grid $\mathcal{T}^H$ with coarse cell $K_i$, $\Gamma_w \cup \Gamma_E = \partial K_i$}
\label{gcoarse}
\end{figure}

Let $\mathcal{T}^H$ be a coarse-grid partition of the domain $\Omega$ with mesh size $H$ (see Figure \ref{gcoarse})
\[
\mathcal{T}^H = \bigcup_{i = 1}^{N_{cell}^H} K_i, 
\]
where $N_{cell}^H$ is the number of coarse grid cells (local domains). We use
$\mathcal{E}^H$ to denote the set of facets in $\mathcal{T}^H$ with  $\mathcal{E}^H = \mathcal{E}^H_{o} \cup \mathcal{E}^H_{b}$. 
For the construction of the coarse grid approximation, we use the Discontinuous Galerkin Generalized Multiscale Finite Element Method (DG-GMsFEM). In  DG-GMsFEM, the multiscale basis functions are supported in each coarse cell $K_i$\cite{chung2016generalized, chung2018multiscale, chung2017conservative}. 
We define $V_{H}$ as the multiscale velocity space and  $P_{H}$ as the multiscale space for concentration. 
For the pressure approximation, we use the piecewise constant function space $Q_{H}$ over the coarse cells.

We denote the multiscale spaces for  concentration and  velocity as
\[
P_H = \text{span} \{ \phi_i \}_{i=1}^{N_c}, 
\quad 
V_H = \text{span} \{ \psi_i \}_{i=1}^{N_u},
\]
respectively,
where  
$N_u =\text{dim}(V_H)$ is the number of basis functions for velocity, 
$N_c =\text{dim}(P_H)$ is the number of basis functions for concentration.  
For the pressure, we use the space of piecewise constant functions over the coarse grid cells
\[
Q_H = \{ q \in L^2(\Omega) : \,  q |_K \in \mathbb{P}_0(K), \, \forall K \in \mathcal{T}^H \},
\]
where $N_p = \text{dim}(Q_H)$ is equal to the number of coarse grid cells ($N_p = N_{cell}^H$).

For the coarse grid approximation, we have the following variational formulation: 
\begin{itemize}
\item Flow problem: find $(u_H, p_H) \in V_H \times Q_H$ such that
\begin{equation}
\label{up-ms}
\begin{split}
\frac{1}{\tau}m^u(u_H - \check{u}_H, v) 
+  a^u_{\text{DG}}(u_H, v)  
+ b_{\text{DG}}(p_H, v)  & = l^u(v) 
\quad \forall v \in V_H, \\
b_{\text{DG}}(u_H, q) &= l^p(q), \quad \forall q \in Q_H.
\end{split}
\end{equation}

\item Transport problem:  find $c_H \in P_H$ such that
\begin{equation}
\label{c-ms}
\begin{split}
\frac{1}{\tau} m^c(c_H - \check{c}_H , r) 
+ a^c_{\text{DG}}(c_H, r) 
+ c^c_{\text{DG}}(c_H, r) &= l^c(r), 
\quad \forall r \in P_H.
\end{split}
\end{equation} 
\end{itemize}

Next, we consider the construction of the multiscale basis functions for velocity and concentration.

\subsection{Multiscale space for velocity}

For the construction of the multiscale space for the velocity, we start with the construction of the snapshot space. The snapshot space is formed by the solution of local problems with all possible boundary conditions up to the fine grid resolution in each coarse cell  $K_i, (i = 1, \cdots, N_{cell}^H)$, where $N_{cell}^H$ is the number of coarse blocks in $\Omega$.  
After that, we solve a local spectral problem to select dominant modes of the snapshot space. 

We consider two types of multiscale spaces: 
\begin{itemize}
\item \textit{Type 1}. Multiscale space for flow is defined so that snapshot spaces and spectral problems are constructed for all flow directions.
\item \textit{Type 2}. Multiscale space for flow is defined so that snapshot spaces and spectral problems are constructed separately for each flow direction.
\end{itemize}

We start with the construction of the \textit{Type 2} multiscale basis functions. 
The local snapshot space consists of functions  $u^{i,r}_l \in K_i$ which are solutions of the following problem
\begin{equation}
\label{eq:snap-u}
\begin{aligned}
-\mu \Delta u^{i,r}_l  + \nabla p &= 0, \quad & x \in K_i, \\
\nabla \cdot u^{i,r}_l  &= f^{i,r}_l, \quad & x \in  K_i,
\end{aligned}
\end{equation}
with boundary conditions
\[
u^{i,r}_l =\delta^{i,r}_l, \quad x \in \Gamma_E, \quad 
u^{i,r}_l = 0, \quad x \in \Gamma_w,
\]
for $l= 1, \cdots, J_i$,  $J_i$ is the number of fine grid facets on  $\Gamma_E$,
and $\delta^{i,r}_l$ is the vector discrete delta function defined on $\Gamma_E$,  $\Gamma_E = \partial K_i / \Gamma_w$ is the interface between local domains ($r = 1,..,d$ and $d = 2,3$ is dimension of problem). 
Here  constant $c$  is chosen by the compatibility condition, $f^{i,r}_l = \frac{1}{|K_i|}\int_{\Gamma_E} \delta^{i,r}_l \cdot n \, ds$. 
For 2D problem, we solve local problems \eqref{eq:snap-u} using the following 
\[
\delta_l^{i,1} = (\delta_l,0), \quad \text{ and } \quad 
\delta_l^{i,2} =(0,\delta_l).
\]  
For 3D problem, we solve the local problems \eqref{eq:snap-u} using the following
\[
\delta_l^{i,1}= (\delta_l,0,0), \quad 
\delta_l^{i,2} = (0,\delta_l,0), \quad \text{ and } \quad 
\delta_l^{i,3} = (0,0,\delta_l).
\]  
In the above, $\delta_l$ is the discrete delta function whose value is $1$ on the $l$-th fine grid node
and $0$ otherwise. 
We remark that the local problems \eqref{eq:snap-u} are solved on the fine mesh by using a standard numerical scheme.
Using these local solutions, we form the following local snapshot space 
\[
V^{i,\text{snap}}_r = \{ u_l^{i,r}: 1 \leq l \leq  J_i \},
\]
where the flow direction, $r = 1,..,d$, are considered separately. 

\begin{figure}[!h]
  \centering
  \includegraphics[width=1.0\linewidth]{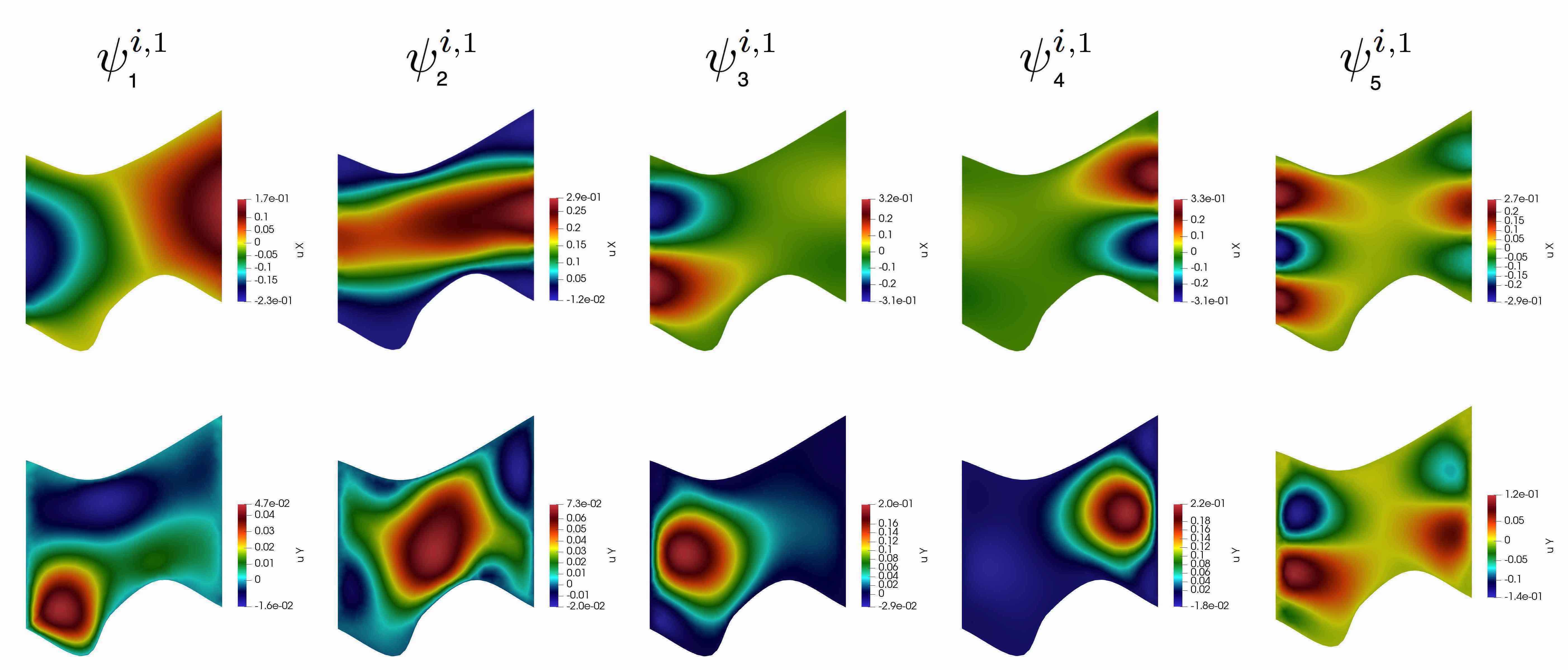}
\caption{Illustration of the multiscale basis functions for velocity 
$\psi_k^{i,r} = (\psi_{k, x}^{i,r}, \psi_{k, y}^{i,r})$ for  $r = 1$ and $k = 1,...,5$ (from left to right).  
\textit{Type 2} ($V_H^1$). 
First row: $\psi_{k, x}^{i,1}$. 
Second row: $\psi_{k, y}^{i,1}$ }
\label{ms-basis-u-t2a}
\end{figure}

\begin{figure}[!h]
  \centering
  \includegraphics[width=1.0\linewidth]{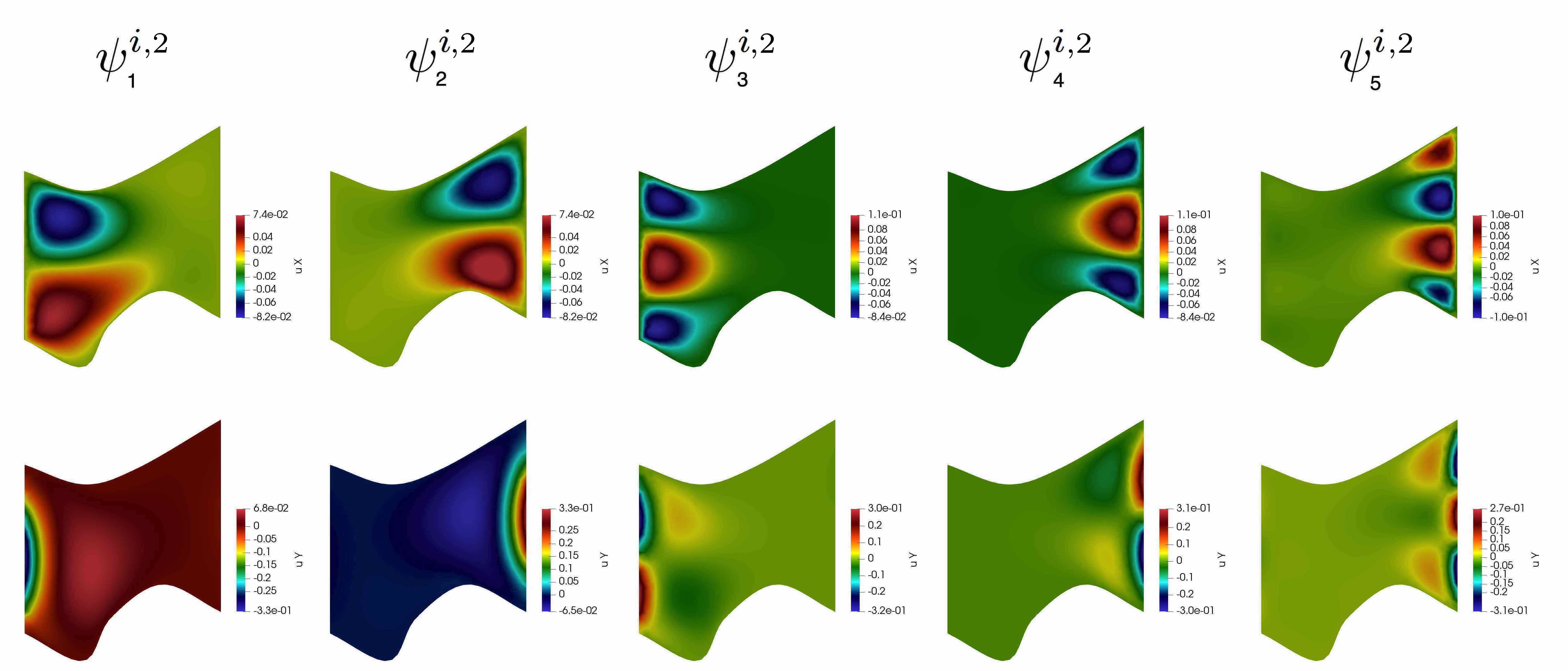}  
\caption{Illustration of the multiscale basis functions for velocity 
$\psi_k^{i,r} = (\psi_{k, x}^{i,r}, \psi_{k, y}^{i,r})$ for  $r = 2$ and $k = 1,...,5$ (from left to right).  
\textit{Type 2} ($V_H^2$). 
First row: $\psi_{k, x}^{i,2}$. 
Second row: $\psi_{k, y}^{i,2}$ }
\label{ms-basis-u-t2b}
\end{figure}

To reduce the size of direction based snapshot spaces, we solve following local spectral problem in each snapshot space $V^{i,\text{snap}}_r $ in $K_i$
\begin{equation}
\label{eq:eigen-u}
\tilde{A}^{u, K_i}_r \tilde{\psi}^{i,r} =  
\lambda  \tilde{S}^{u,K_i}_r \tilde{\psi}^{i,r}, 
\end{equation}
where 
\[
\tilde{A}^{u, K_i}_r = R^{u,r}_{i, \text{snap}} A_h^{u, K_i} (R^{u,r}_{i, \text{snap}})^T,
\quad
\tilde{S}^{u, K_i}_r = R^{u,r}_{i, \text{snap}} S_h^{u, K_i} (R^{u,r}_{i, \text{snap}})^T,  
\]
and $A_h^{u, K_i}$ is the matrix representation of the bilinear form $a^{u, K_i}_{DG}(u, v)$ and $S_h^{u, K_i}$ is the matrix representation of the bilinear form $s^{u, K_i}(u, v)$
\begin{equation}
\label{eq:eigen-u-as}
\begin{split}
a^{u, K_i}_{\text{DG}}(u,v) 
&= \sum_{K \in \mathcal{T}^h(K_i)} \int_K
\mu \nabla u \cdot \nabla v \, dx \\
&- \sum_{E \in  \mathcal{E}_o^h(K_i)} \int_E \Big( 
\{ \mu \nabla u \cdot n \}\cdot [v] + 
\{ \mu \nabla v \cdot  n\}  \cdot [u] -  
\frac{\gamma_u}{h} \{\mu\} [u]\cdot [v] \Big) \, ds,
\\
s^{u, K_i} (u,v) &= \sum_{E \in  \mathcal{E}_b^h(K_i)} \int_E  u \cdot v \, dx.
\end{split}
\end{equation}
The above 
snapshot space projection matrix is defined by collecting all local solutions  
\[
R^{u,r}_{i, \text{snap}} = \left[ 
u_{1}^{i,r}, \ldots,  u_{J_i}^{i,r} 
\right]^T, \quad r = 1,..,d.
\]
We remark that the integral in $s^{u, K_i}(u,v)$ is defined on the boundary of the coarse block, and $\mathcal{T}^h(K_i)$ is the fine grid for local domain $K_i$. 

To form a direction based multiscale space for velocity, we arrange the eigenvalues in increasing order and choose eigenvectors corresponding to the first $M^u_i$ the smallest eigenvalues for each flow direction
\[
V_H^{r} = \text{span} \{ \psi_k^{i,r}: \, 1 \leq i \leq N_{cell}^H, \,  1 \leq k \leq M^u_i \},
\]
where $\psi^{i,r}_k = R^{u,r}_{i, \text{snap}} \tilde{\psi}^{i,r}_k$ $k = 1,...,M^u_i$. In Figures \ref{ms-basis-u-t2a} and \ref{ms-basis-u-t2b}, we depicted the first five multiscale basis functions in $V_H^1$ and $V_H^2$, respectively.
Finally, the multiscale space of \textit{Type 2} is defined as follows
\[
V_H =  V_H^{1} +  V_H^{2} \text{ for 2D}, \quad 
V_H =  V_H^{1} +  V_H^{2} +  V_H^{3} \text{ for 3D}.
\]

\begin{figure}[h!]
\centering
\includegraphics[width=1.0\textwidth]{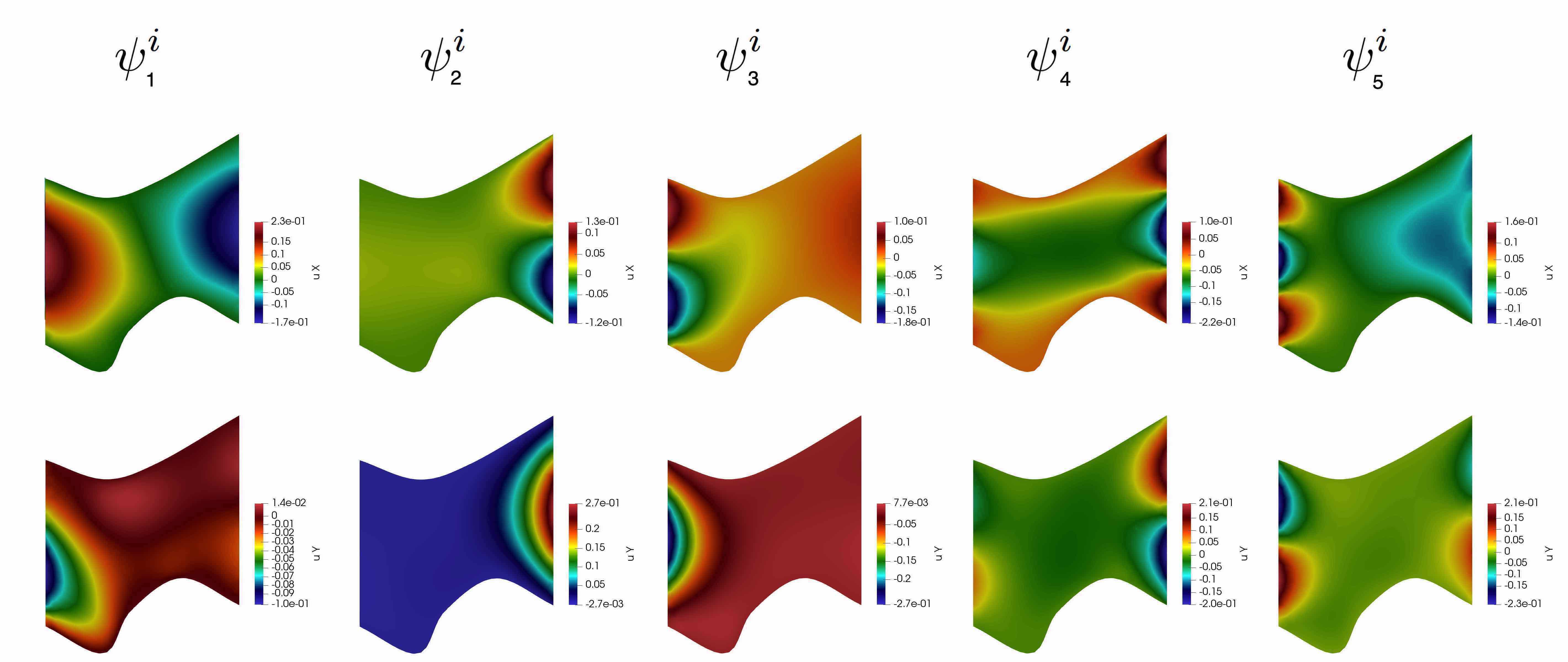}
\caption{Illustration of the multiscale basis functions for velocity 
$\psi_k^{i} = (\psi_{k, x}^{i}, \psi_{k, y}^{i})$ for $k = 1,...,5$ (from left to right).  \textit{Type 1}.  
First row: $\psi_{k, x}^{i}$. 
Second row: $\psi_{k, y}^{i}$ }
\label{ms-basis-u-t1}
\end{figure}

The \textit{Type 1} multiscale space is constructed in a similar way. Instead of flow separation by directions, we define the snapshot space as the collection of the local solutions for flow in all directions
\[
V^{i,\text{snap}} = \{ u_l^{i,r}: 1 \leq l \leq J_i, \, r = 1,..,d \}.
\]
where $d = 2,3$.  
To construct a multiscale space, we perform a dimension reduction by the solution of the local spectral problem on snapshot space $V^{i,\text{snap}}$ 
\begin{equation}
\label{eq:eigen-u-t1}
\tilde{A}^{u, K_i} \tilde{\psi}^{i} =  
\lambda  \tilde{S}^{u,K_i} \tilde{\psi}^{i}, 
\end{equation}
with
$\tilde{A}^{u, K_i} = R^{u}_{i, \text{snap}} A_h^{u, K_i} (R^{u}_{i, \text{snap}})^T$, 
$\tilde{S}^{u, K_i} = R^{u}_{i, \text{snap}} S_h^{u, K_i} (R^{u}_{i, \text{snap}})^T$
and   
\[
R^u_{i, \text{snap}} = \left[ 
u_{1}^{i,1}, \ldots u_{1}^{i,d}, \ldots, 
u_{J_i}^{i,1}, \ldots  u_{J_i}^{i,d} 
\right]^T.
\]
where $A_h^{u, K_i}$ and $S_h^{u, K_i}$ are given in \eqref{eq:eigen-u-as}.
We arrange the eigenvalues in increasing order and choose the first  eigenvectors corresponding to the first $M^u_i$ the smallest eigenvalues  as the basis functions
\[
V_H = \text{span} \{ \psi_k^i: \, 1 \leq i \leq N_{cell}^H, \,  1 \leq k \leq M^u_i \}. 
\]
where $\psi_k^i = R^u_{i, \text{snap}} \tilde{\psi}^i_k$,  $k = 1,...,M^u_i$. See Figure \ref{ms-basis-u-t1} for an illustration of \textit{Type 1} multiscale basis functions.


\subsection{Multiscale space for concentration}

The construction of the multiscale space for the concentration has a similar concept. The space is constructed for each coarse cell  $K_i$ (local domain) for $i = 1, \cdots, N_{cell}^H$.  
We consider two types of multiscale spaces: 
\begin{itemize}
\item \textit{Type 1}. Multiscale space for transport from all boundaries. Snapshots are constructed for flow in all directions with corresponded spectral problems to dominant mode extraction.
\item \textit{Type 2}. Multiscale space for  each boundary transport
 separately. Snapshot spaces are constructed separately with corresponded spectral problems for each of them.
\end{itemize}
We will consider two types of boundaries: (1) interface between local domains $\Gamma_E$ and (2) reactive wall boundaries $\Gamma_w$.  
This concept is based on the definition of the coarse grid variables and is similar to the approach that we used in \cite{chung2018non, vasilyeva2019nonlocal}. 


We will present the construction of the multiscale basis functions for three types of wall boundary conditions for the transport problem 
\begin{itemize}
\item Dirichlet boundary conditions (\textit{DBC}):
\begin{equation}
\label{eq:mm-dbc}
c = g, \quad x \in \Gamma_w.
\end{equation}

\item Neumann boundary conditions (\textit{NBC}):
\begin{equation}
\label{eq:mm-dbc}
- D \nabla c  \cdot n = g, \quad x \in \Gamma_w.
\end{equation}

\item Robin boundary conditions (\textit{RBC}):
\begin{equation}
\label{eq:mm-dbc}
- D \nabla c  \cdot n = \alpha (c  - g), \quad x \in \Gamma_w.
\end{equation}

\end{itemize}

We start with the construction of the \textit{Type 2} multiscale basis function. 
The local snapshot space consists of functions  $c^{i,r}_l$ which are solutions of the following local problem
\begin{equation}
\label{eq:snap-c}
- \nabla (D \nabla c^{i,r}_l)   = 0, \quad x \in  K_i.
\end{equation}
where boundary conditions depend on the type of non-homogeneous  boundary conditions on wall boundary
\begin{itemize}
\item  $r = 1$ for transport between local domains on interface $\Gamma_E$, we set 
\[
c^{i,r}_l =\delta_i^l, \quad x \in \Gamma_E,  
\]
and 
\[
c^{i,r}_l = 0, \quad x \in \Gamma_w \text{ (DBC)},
\] \[
- D \nabla c^{i,r}_l  \cdot n = 0, \quad x \in \Gamma_w \text{ (NBC)},
\] \[
- D \nabla c^{i,r}_l  \cdot n = \alpha c, \quad x \in \Gamma_w \text{ (RBC)}.
\]
\item  $r = 2$ for transport from walls boundary
\[
c^{i,r}_l = \delta_i^l, \quad x \in \Gamma_w \text{ (DBC)},
\] \[
- D \nabla c^{i,r}_l  \cdot n = \delta_i^l, \quad x \in \Gamma_w \text{ (NBC)},
\] \[
- D \nabla c^{i,r}_l  \cdot n = \alpha (c - \delta_i^l), \quad x \in \Gamma_w \text{ (RBC)}.
\]
and
\[
- D \nabla c^i_l \cdot n = 0, \quad x \in \Gamma_E,  
\]
\end{itemize}

Here $l= 1, \cdots,  J_i$, where $J_i$ is the number of fine grid facets on the boundary of $K_i$,
and $\delta_i^l$ is the discrete delta function defined on $\partial K_i$ and equal to 1 if $i = l$ and zero otherwise.  
This problem is solved on the fine mesh using an appropriate numerical scheme.

Using these local solutions, we form a local snapshot space  for concentration in $K_i$
\[
P^{i,\text{snap}}_r = \{ c^{i,r}_l: 1 \leq l \leq J_i \}, 
\]
and define the snapshot space projection matrix
\[
R^{c,r}_{i, \text{snap}} = \left[ c_1^{i,r}, \ldots, c_{J_i}^{i,r} \right]^T.
\] 
where $r = 1,2$.

\begin{figure}[h!]
\centering
\includegraphics[width=0.9\textwidth]{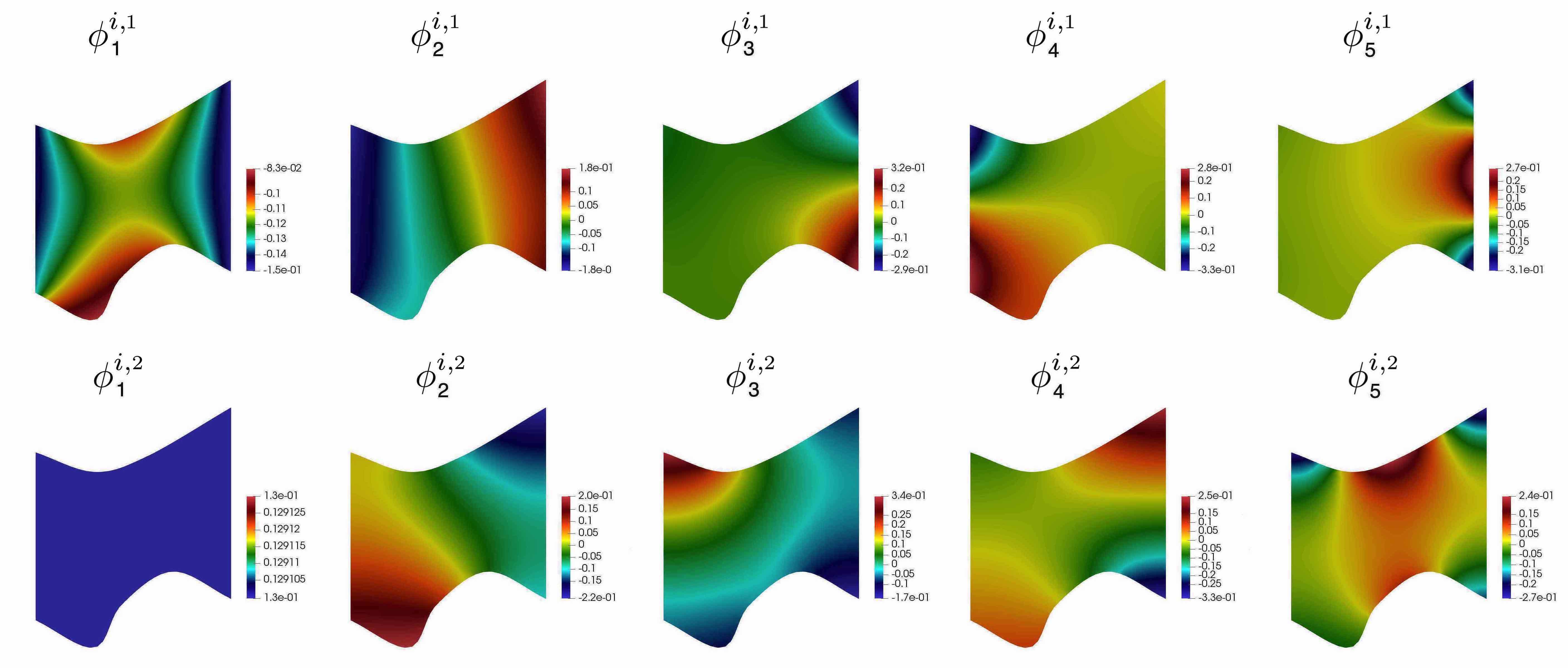}
\caption{Illustration of the multiscale basis functions for concentration 
$\phi_k^{i,r}$ for  $r = 2$ and $k = 1,...,5$ (from left to right).  
\textit{Type 2} for Robin boundary conditions on $\Gamma_w$ (\textit{RBC}). 
First row: $\phi_{k}^{i,1}$. 
Second row: $\phi_{k}^{i,2}$ }
\label{ms-basis-c-t2}
\end{figure}

To reduce the size of the snapshot space, we solve the following local spectral problem in the snapshot space $P^{i,\text{snap}}_r$
\begin{equation}
\label{eq:off-eq}
\tilde{A}^{c, K_i}_r \tilde{\phi}^{i,r} =  
\eta \tilde{S}^{u,K_i}_r \tilde{\phi}^{i,r}, 
\end{equation}
where 
\[
\tilde{A}^{c, K_i}_r = R^{c,r}_{i, \text{snap}} A_h^{c, K_i} (R^{c,r}_{i, \text{snap}})^T,
\quad
\tilde{S}^{c, K_i}_r = R^{c,r}_{i, \text{snap}} S_h^{c, K_i} (R^{c,r}_{i, \text{snap}})^T,  
\]
and $A_h^{c, K_i}$ is the matrix representation of the bilinear form $a^{c, K_i}_{DG}(c, r)$ and $S_h^{c, K_i}$ is the matrix representation of the bilinear form $s^{u, K_i}(c, r)$
\begin{equation}
\label{eq:off-c-as}
\begin{split}
a^{c, K_i}_{DG}(c, z) 
& = \sum_{K \in \mathcal{T}^h(K_i)} 
\int_K D \nabla c \cdot \nabla z \, dx\\
&-\sum_{E \in \mathcal{E}^h_0(K_i)} \int_E \Big(
 \{ D \nabla c \cdot n \} [z] + \{ D \nabla z \cdot n \} [c] - 
\frac{\gamma_c}{h} \{D\} [c][z] \Big) \, ds
\\
s^{c, K_i}(c, z) & = \sum_{E \in  \mathcal{E}_b^h(K_i)} \int_E c \, z \, ds.
\end{split}
\end{equation}
We remark that the integral in $s^{c, K_i}(c,r)$ is defined on the boundary of the coarse block and $\mathcal{T}^h(K_i)$ is the fine grid for local domain $K_i$. 

We arrange the eigenvalues in increasing order and choose the first  eigenvectors corresponding to the first $M^c_i$ the smallest eigenvalues  as the basis functions 
\[
P^r_H 
= \text{span} \{ 
\phi_k^{i,r}: \, 1 \leq i \leq N_{cell}^H, \,  1 \leq k \leq M^c_i \}. 
\]
where $\phi_k^{i,r} = R^{c,r}_{i, \text{snap}} \tilde{\phi}^{i,r}_k$ for $k = 1,...,M^c_i$ and $r =1,2$.  
Finally, the multiscale space of \textit{Type 2} is defined as follows
\[
P_H =  P_H^{1} \times  P_H^{2}.
\] 
See Figure \ref{ms-basis-c-t2} for an illustration of the multiscale basis functions.

\begin{figure}[h!]
\centering
\includegraphics[width=0.9\textwidth]{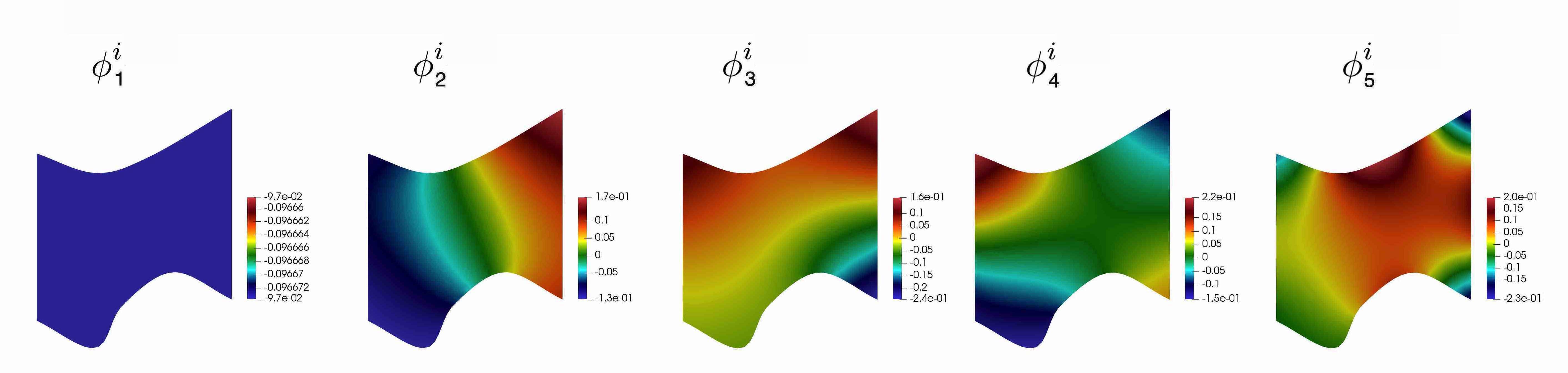}
\caption{Illustration of the multiscale basis functions for concentration. Type 1}
\caption{Illustration of the multiscale basis functions for concentration 
$\phi_k^{i}$ for $k = 1,...,5$ (from left to right).  
\textit{Type 1}.  }
\label{ms-basis-c-t1}
\end{figure}

The \textit{Type 1} multiscale space is constructed similarly. In the snapshot space, we collect all possible boundary conditions on $\partial K_i$ and solve the following problem
\begin{equation}
\label{eq:snap-c-all}
\begin{split}
- \nabla (D \nabla c^{i}_l) = 0, \quad x \in  K_i, \\
c^{i}_l =\delta_i^l, \quad x \in \partial K_i. 
\end{split}
\end{equation}
Local solutions are collected as a snapshot space
\[
V^{i,\text{snap}} = \{ c_l^{i}: 1 \leq l \leq J_i\}.
\]
Dimension reduction of the snapshot space is performed by the solution of the local spectral problem on snapshot space
\begin{equation}
\label{eq:off-eq}
\tilde{A}^{c, K_i} \tilde{\phi}^{i} =  
\eta \tilde{S}^{u,K_i} \tilde{\phi}^{i}, 
\end{equation}
where 
$\tilde{A}^{c, K_i} = R^{c}_{i, \text{snap}} A_h^{c, K_i} (R^{c}_{i, \text{snap}})^T$, 
$\tilde{S}^{c, K_i} = R^{c}_{i, \text{snap}} S_h^{c, K_i} (R^{c}_{i, \text{snap}})^T$ 
and   
\[
R^{c}_{i, \text{snap}} = \left[ c_1^i, \ldots, c_{J_i}^i \right]^T.
\]
where $A_h^{u, K_i}$ and $S_h^{u, K_i}$ are given in \eqref{eq:off-c-as}. 
We arrange the eigenvalues in increasing order and choose the first  eigenvectors corresponding to the first $M^c_i$ the smallest eigenvalues  as the basis functions
\[
P_H 
= \text{span} \{ 
\phi_k^{i}: \, 1 \leq i \leq N_{cell}^H, \,  1 \leq k \leq M^c_i \}. 
\]
where $\phi_k^{i} = R^{c}_{i, \text{snap}} \tilde{\phi}^{i}_k$ for $k = 1,...,M^c_i$.  Multiscale basis functions of \textit{Type 1} are presented in  Figure \ref{ms-basis-c-t1}. We note that, the \textit{Type 1} basis functions do not depend on type of boundary conditions on $\Gamma_w$.

\begin{figure}[h!]
\centering
\includegraphics[width=0.9\textwidth]{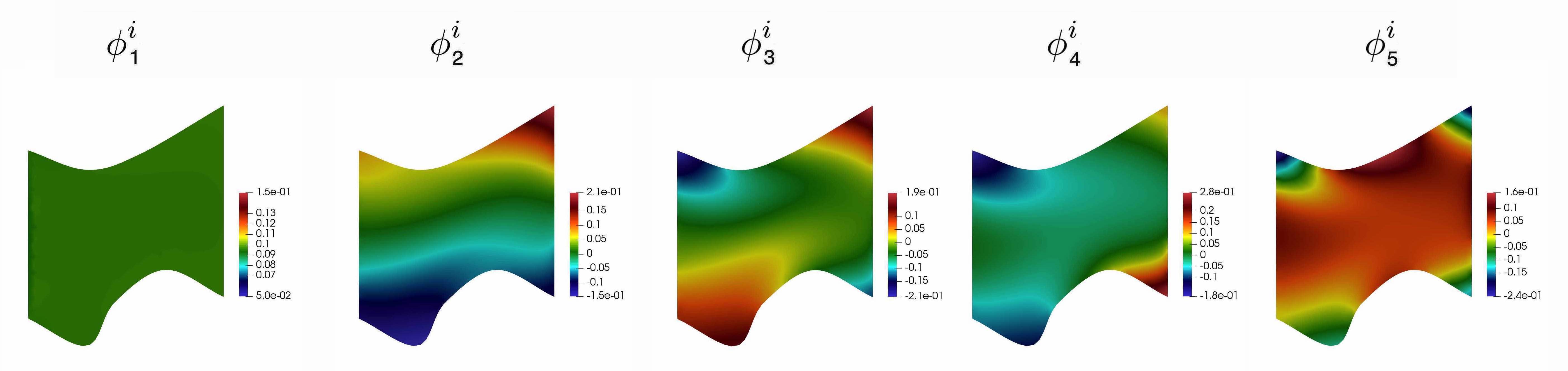}
\caption{Illustration of the multiscale basis functions for concentration 
$\phi_k^{i}$ for $k = 1,...,5$ (from left to right).  
\textit{Type 1} with velocity and time  }
\label{ms-basis-c-t1b}
\end{figure}

\begin{figure}[h!]
\centering
\includegraphics[width=0.9\textwidth]{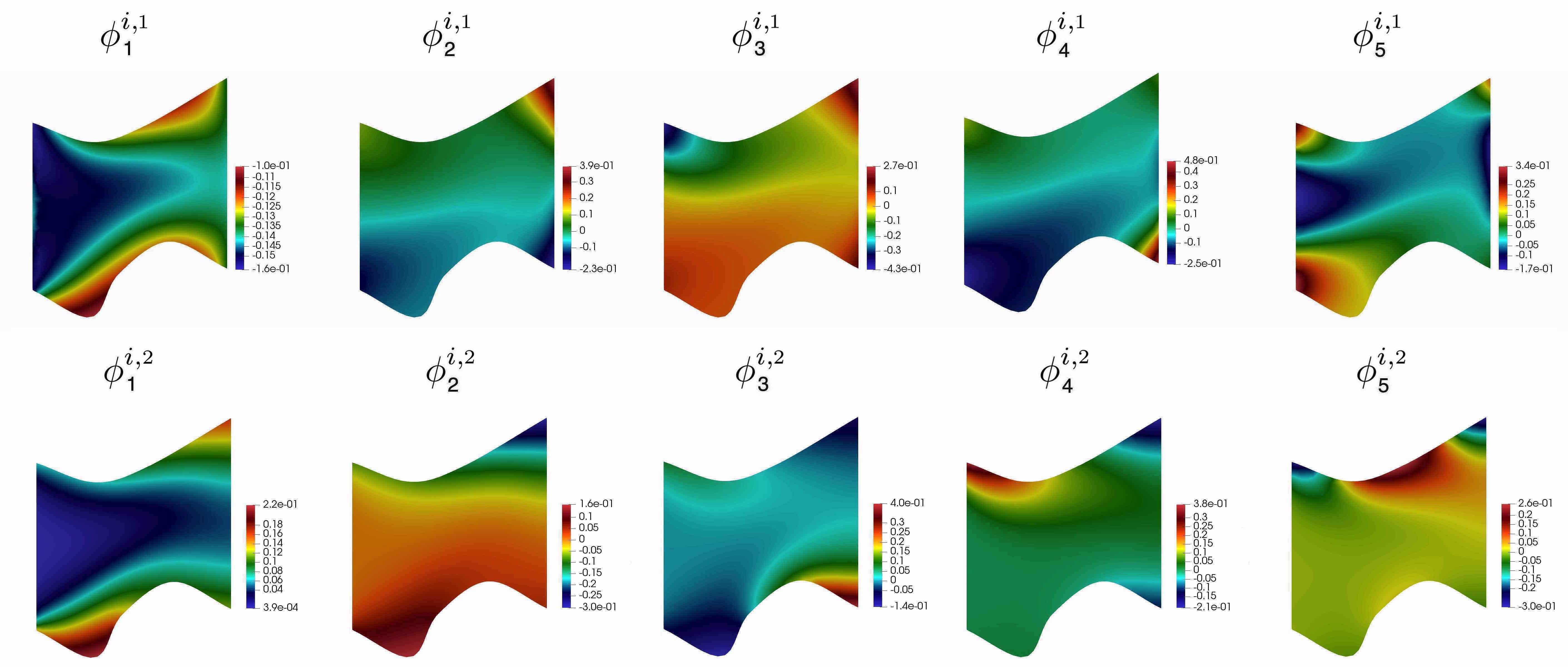}
\caption{Illustration of the multiscale basis functions for concentration 
$\phi_k^{i,r}$ for  $r = 2$ and $k = 1,...,5$ (from left to right).  
\textit{Type 2} with velocity and time for Robin boundary conditions on $\Gamma_w$ (\textit{RBC}) }
\label{ms-basis-c-t2b}
\end{figure}

\begin{remark}
For the use of the unsteady convection-diffusion problem instead of the elliptic problem \eqref{eq:snap-c}  for the construction of the snapshot space, we can use a local problem formulation that is equivalent to the global problem \cite{akkutlu2016multiscale, mehmani2019multiscale}
\begin{equation}
\label{eq:snap-ct}
\frac{1}{\tau}  c^{i}_l + u \nabla c^{i}_l - \nabla (D \nabla c^{i}_l) = 0, \quad x \in  K_i.
\end{equation}
We will use notations "Type 1 with time and velocity"  and "Type 2 with time and velocity"   in the numerical results section. 
Illustration of the multiscale basis functions for \textit{Type 1} and \textit{Type 2} with time and velocity are presented in Figures \ref{ms-basis-c-t1b} and  \ref{ms-basis-c-t2b}, respectively. 
\end{remark}

\begin{figure}[h!]
\centering
\includegraphics[width=0.5\textwidth]{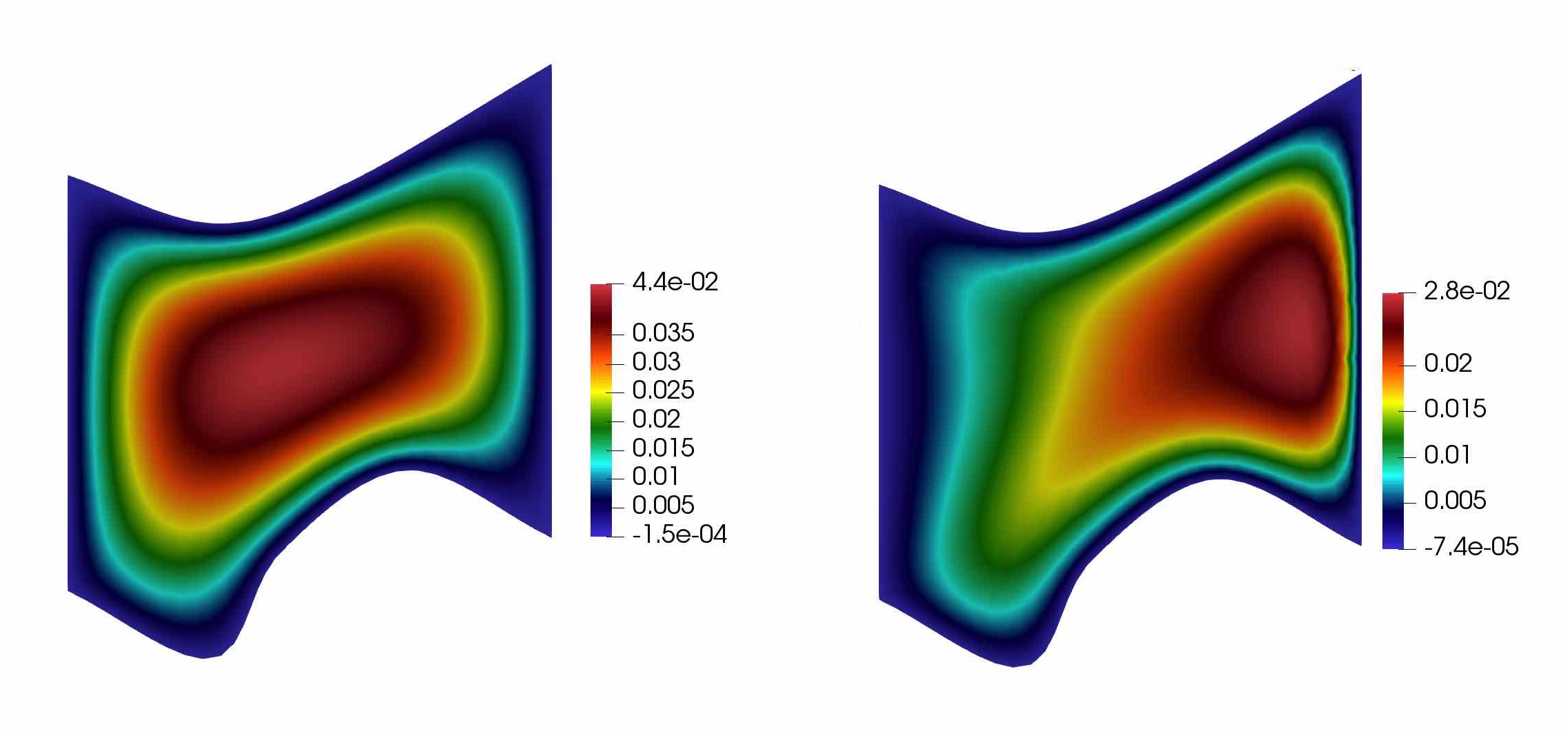}
\caption{Additional basis functions. 
Left: local problem  \eqref{eq:snap-c} with  zero Dirichlet boundary conditions. 
Right: local problem  \eqref{eq:snap-ct} with  zero Dirichlet boundary conditions. }
\label{ms-basis-c-int}
\end{figure}

The presented basis functions are associated with the boundary basis functions  \cite{chung2014adaptive, chung2016multiscale}.     
In addition to the presented multiscale basis functions, we add one interior basis functions, that is constructed by the solution of the corresponding local  problem \eqref{eq:snap-c} or \eqref{eq:snap-ct} with zero Dirichlet boundary conditions on  $x \in \partial K_i$ (see Figure \ref{ms-basis-c-int}).

\subsection{Coarse scale system}

To construct the coarse grid system, we construct projection  matrices using the computed multiscale basis functions for velocity and concentration
\[
R_u = \left[ \psi_1 , \ldots, \psi_{N_u} \right]^T,
\quad 
R_p = \left[ \eta_1 , \ldots, \eta_{N_p} \right]^T, 
\quad 
R_c= \left[ \phi_1 , \ldots, \phi_{N_c} \right]^T,
\]
where we used a single index notation. 
For \textit{Type 1}, we have 
$N_c = \sum_{i = 1}^{N_{cell}^H} M^c_i$  and  
$N_u = \sum_{i = 1}^{N_{cell}^H} M^u_i$. 
For \textit{Type 2}, we have 
$N_c = 2 \cdot \sum_{i = 1}^{N_{cell}^H} M^c_i$  and  
$N_u = d \cdot  \sum_{i = 1}^{N_{cell}^H} M^u_i$. 
Note that, we use the space of piecewise constant functions for $R_p$ over the coarse grid, and set $\eta_i(x)$ equal to 1 if $x \in K_i$ and zero otherwise  ($N_p = N_{cell}^H$). 

Using these matrices, we have the following computational systems in matrix form:
\begin{itemize}
\item Flow problem:
\begin{equation}
\begin{split}
\frac{1}{\tau}
\begin{pmatrix}
M^u_H & 0 \\
0 & 0
\end{pmatrix}
\begin{pmatrix}
u_H - \check{u}_H  \\
p_H - \check{p}_H
\end{pmatrix} + 
\begin{pmatrix}
A^u_H & B_H^T \\
B_H & 0
\end{pmatrix}
\begin{pmatrix}
u_H \\
p_H
\end{pmatrix} = 
\begin{pmatrix}
F_H^u \\
F_H^p
\end{pmatrix}
\end{split}
\end{equation} 
where 
\[
M^u_H = R_u M^u_h R_u^T, \quad 
A^u_H = R_u A^u_h R_u^T, \quad 
B_H = R_u B_h R^T_p, \quad 
F_H^u = R_u F_h^u, \quad
F_H^u = R_p F_h^p,
\] 
and after the solution of the coarse-scale approximation, we reconstruct velocity on a fine grid  $u_{ms}= R_u^T u_H$.

\item Transport problem:
\begin{equation}
\frac{1}{\tau} M^c_H  (c_H - \check{c}_H) + (A^c_H + C^c_H(u_{ms})) c_H =  F^c_H.
\end{equation} 
where 
\[
M^c_H = R_c M^c_h R_c^T, \quad 
A^c_H = R_c A^c_h R_c^T, \quad 
C^c_H(u_{ms}) = R_c C^c_h(u_{ms}) R_c^T, \quad 
F_H^c = R_c F_h^c,
\] 
and reconstruct concentration on the fine grid  $c_{ms}= R_c^T c_H$.
\end{itemize}


\section{Numerical results}

\begin{figure}[h!]
\centering
\begin{subfigure}{1.0\textwidth}
\centering
\includegraphics[width=0.45\textwidth]{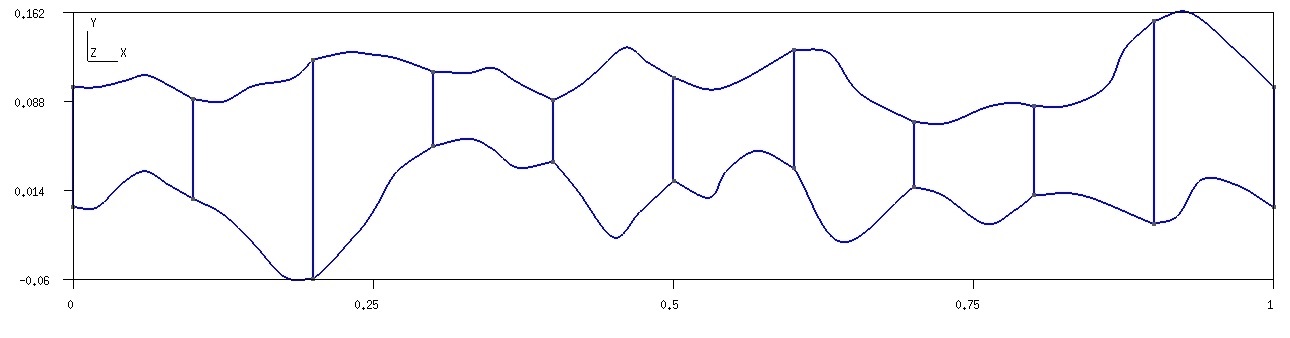}
\includegraphics[width=0.45\textwidth]{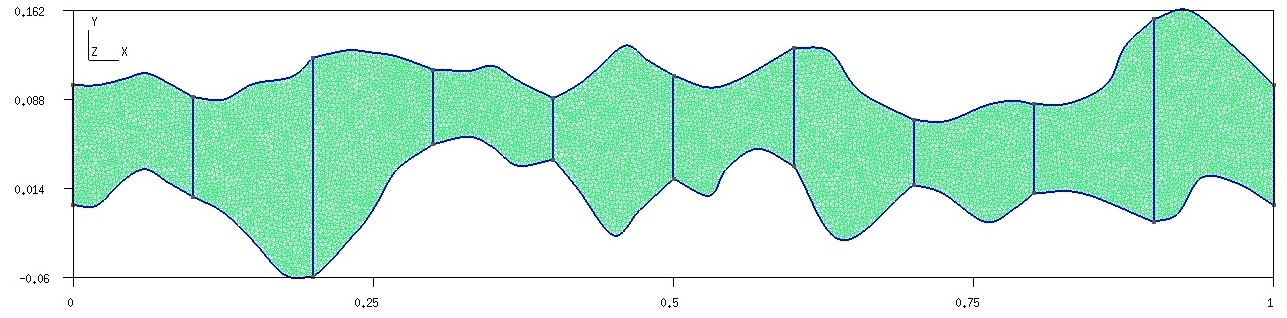}
\caption{Geometry 1 (fine grid with 17350 cells and coarse grid with 10 local domains)}
\end{subfigure}
\begin{subfigure}{1.0\textwidth}
\centering
\includegraphics[width=0.45\textwidth]{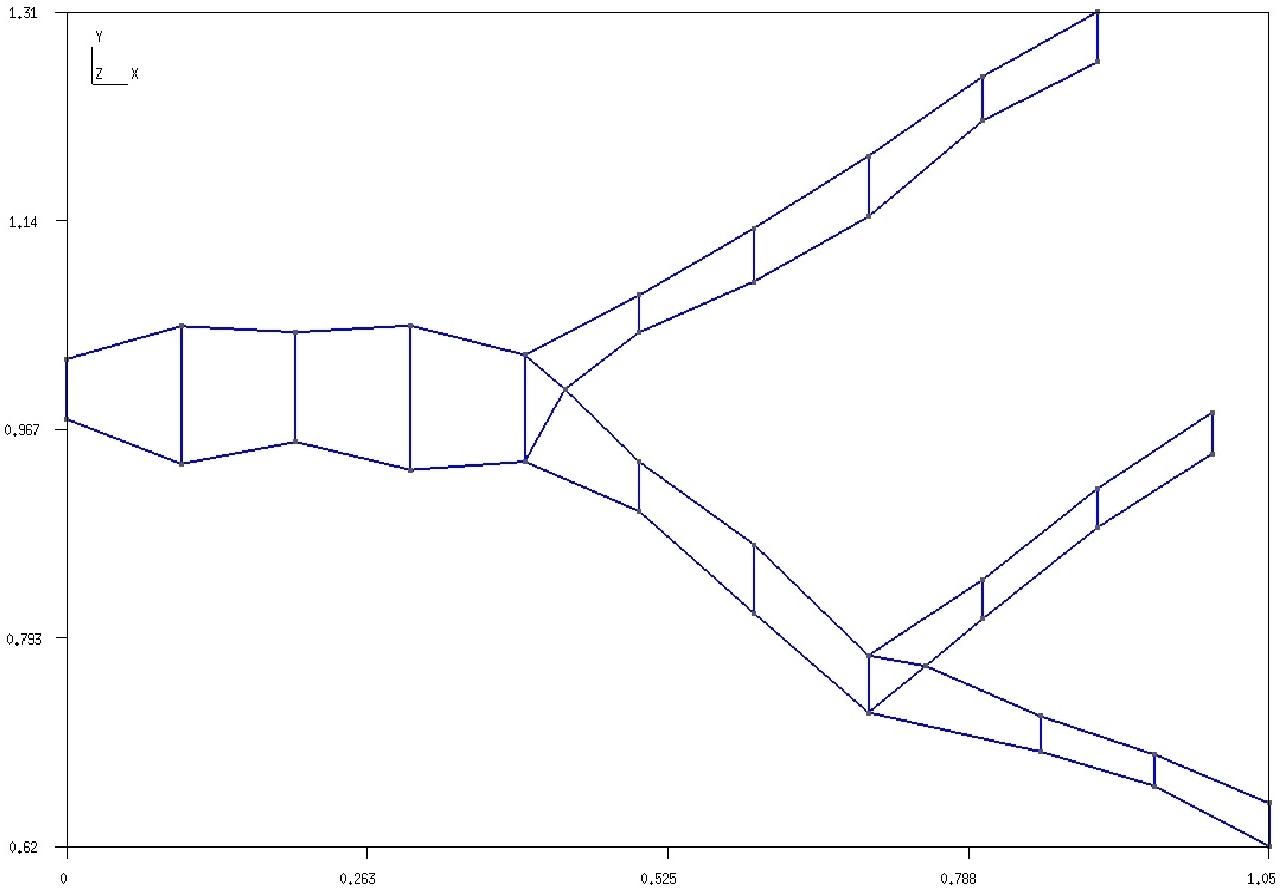}
\includegraphics[width=0.45\textwidth]{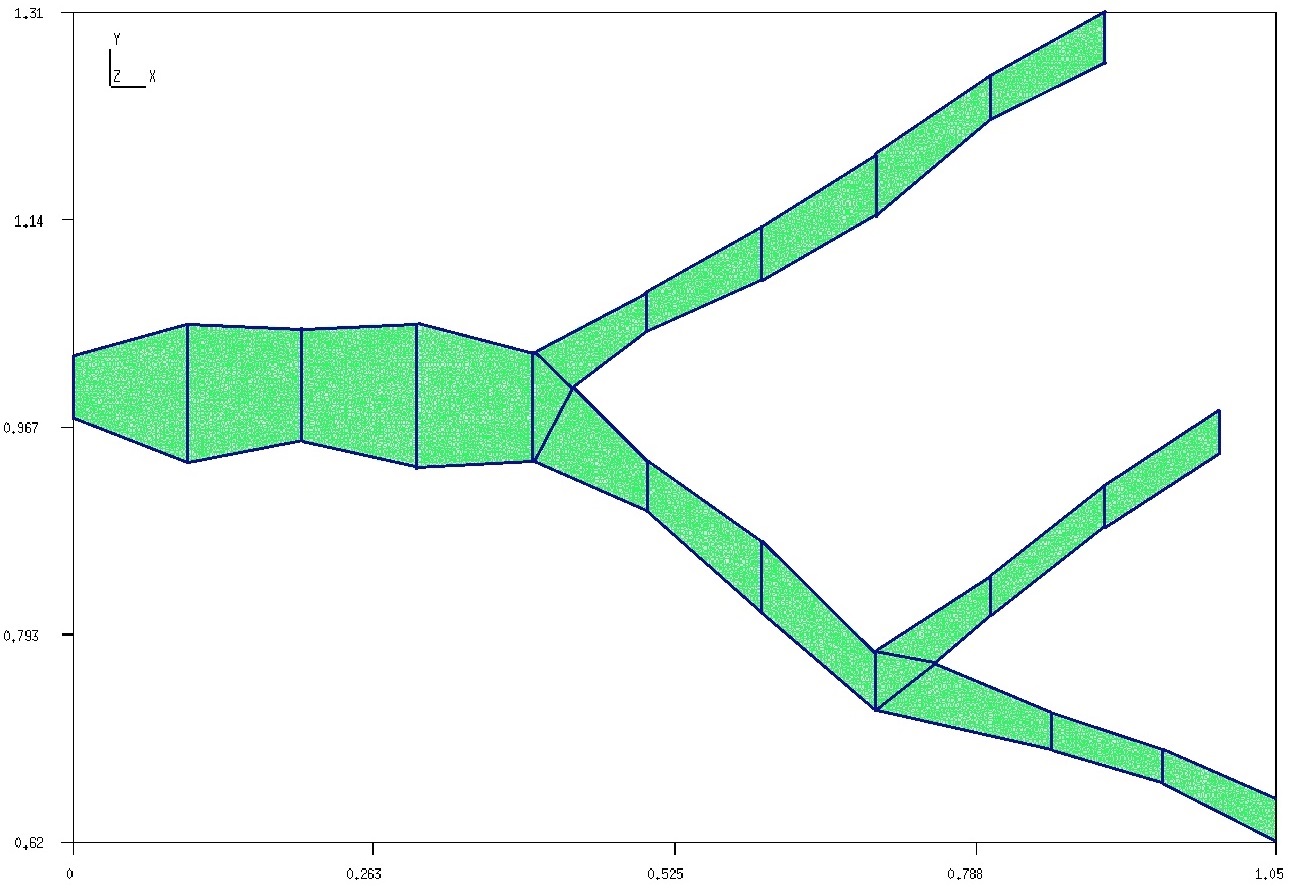}
\caption{Geometry 2 (fine grid with 18021 cells and coarse grid with 20  local domains)}
\end{subfigure}
\begin{subfigure}{1.0\textwidth}
\centering
\includegraphics[width=0.45\textwidth]{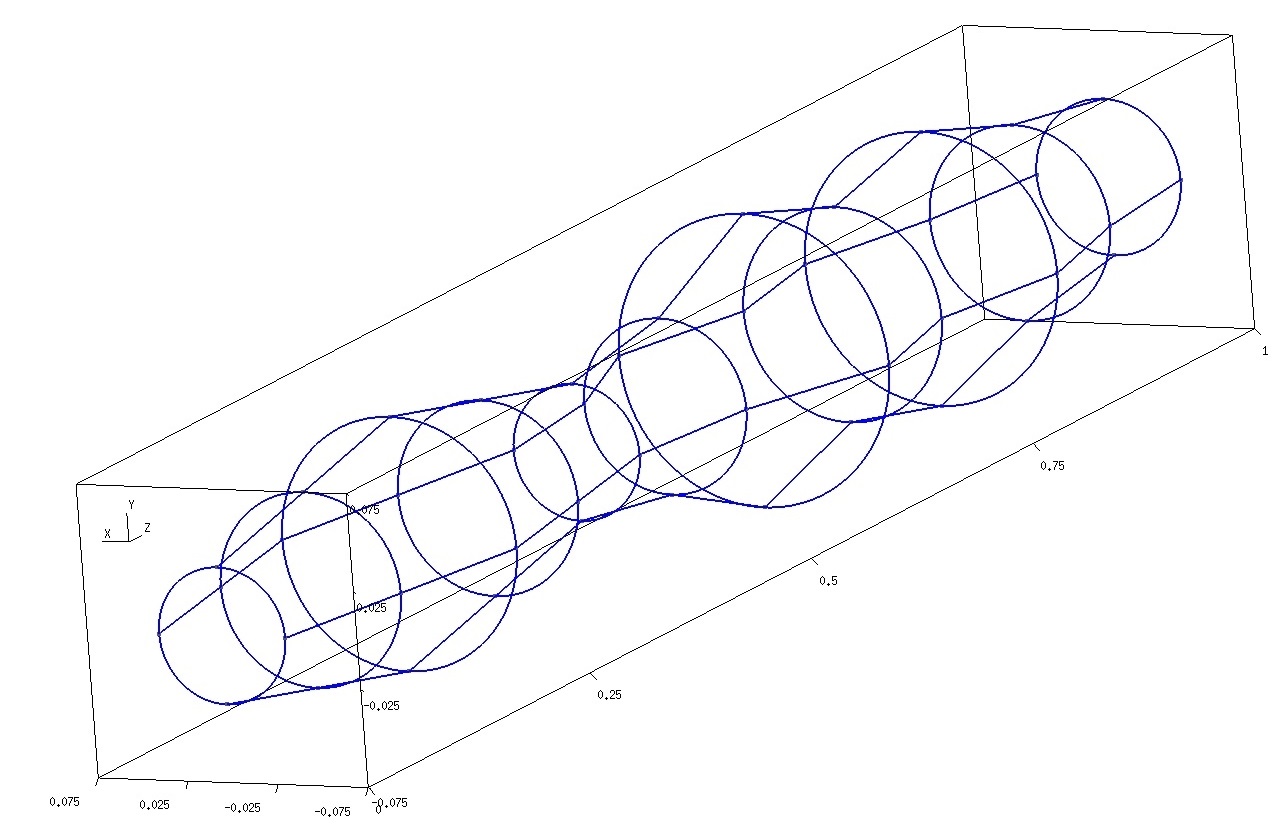}
\includegraphics[width=0.45\textwidth]{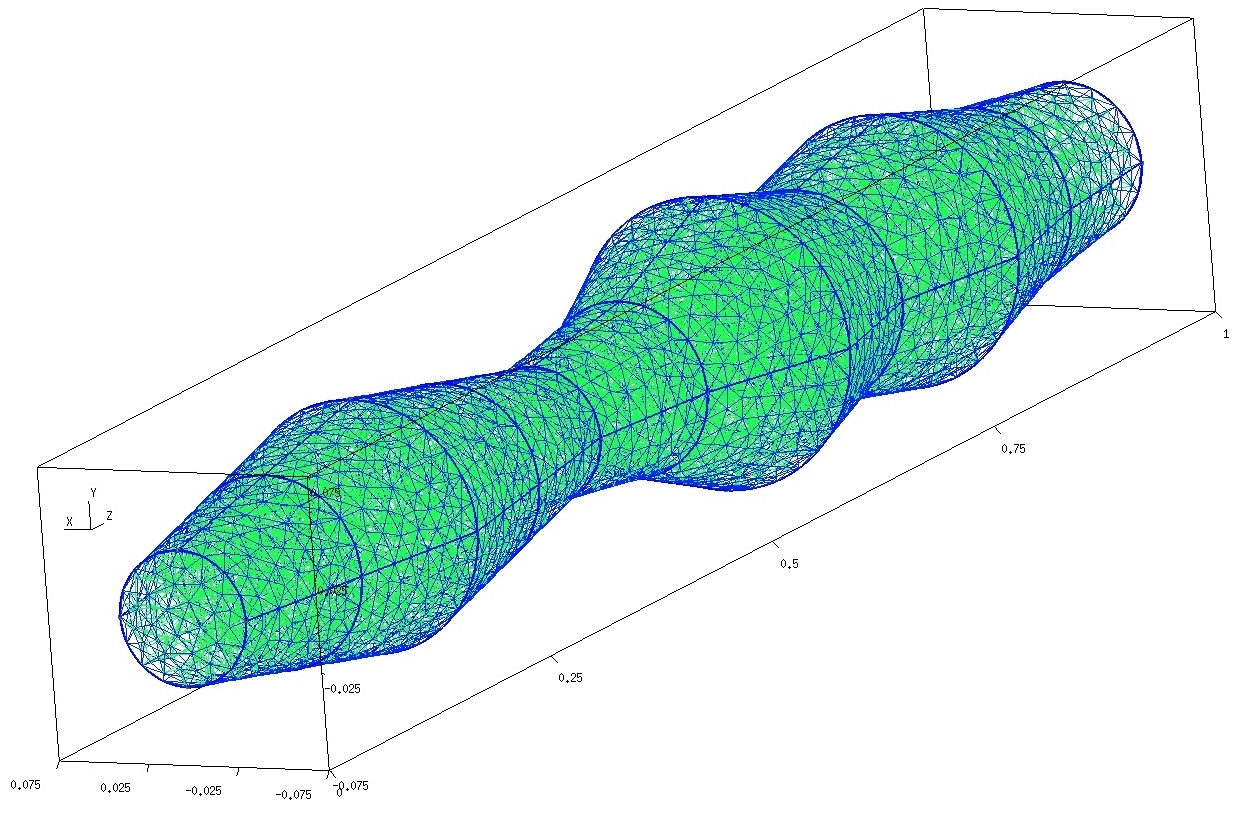}
\caption{Geometry 3 (fine grid with 15094 cells and coarse grid with 10 local domains)}
\end{subfigure}
\caption{Coarse grid and fine grid for Geometry 1, 2 and 3. 
Left: computation geometry and coarse grid.  Right: fine grid.}
\label{geom1}
\end{figure}

In this section, we will present some numerical results. 
We will use the following three computational domains (Figure \ref{geom1}) to demonstrate the performance of our method:
\begin{itemize}
\item[] \textit{Geometry 1} with fine grid that contains 17350 cells. Coarse grid contains 10  local domains.
\item[] \textit{Geometry 2} with fine grid that contains 18021 cells. Coarse grid contains 20 local domains.
\item[] \textit{Geometry 3} with fine grid that contains 15094 cells. Coarse grid contains 10 local domains.
\end{itemize}
In order to construct structured coarse grids, we explicitly add lines between local domains (coarse grid cells) in geometry construction. 
For unstructured coarse grid, local domains have a  rough interface between local domains. 
Computational domains with fine and coarse grids are presented in Figure \ref{geom1} for Geometry 1, 2, and 3.   
We use Gmsh to construct computational geometries and unstructured fine grids \cite{geuzaine2009gmsh}. The numerical implementation is based on the FEniCS library \cite{logg2012automated}.

To investigate the presented multiscale method for solving problems in thin domains with non-homogeneous boundary conditions, we consider the following test:
\begin{itemize}
\item \textit{Test 1} for different geometries of computational domain. 
We consider Geometries 1, 2 and 3 (see Figure \ref{geom1}). We consider problems with non-homogeneous Robin boundary conditions for concentration. 
\item \textit{Test 2} for different boundary conditions and diffusion coefficients ($D = 0.01$, $0.1$ and $1$). We consider Dirichlet and Neumann non-homogeneous boundary conditions on the wall boundary $\Gamma_w$. In this test, we also investigate different types of multiscale basis functions in detail.  
\item \textit{Test 3} for structured and unstructured coarse grids. We investigate the performance of the multiscale method with a rough interface between local domains (unstructured coarse grid). We also investigate the influence of the multiscale velocity accuracy to the concentration errors.
\end{itemize}

For flow problem, we set $\mu=1$, $\rho = 1$ and $u_0 = 0$. 
We set $g = (\tilde{g},0)$ for Geometry 1,2 and $g = (0, 0, \tilde{g})$ for Geometry 3 as inflow boundary condition $\Gamma_{in}$, 
where $\tilde{g}(x) = u_{in} n^{-1} (n+2) (1 - (r/r_{max})^n)$ with $n = 2$ and $u_{in} = 1$ \cite{dobroserdova2019multiscale, oshima2001finite}. Here $r$ is the distance to center point $x_0$ and $r_{max}$ is the radius of left boundary $\Gamma_{in}$, where 
$x_0 = (0, 0.05)$, $r_{max} = 0.05$ for Geometry 1, 
$x_0 = (0, 1)$, $r_{max} = 0.025$ for Geometry 2, and  
$x_0 = (0, 0, 1)$, $r_{max} = 0.035$ for Geometry 3. 

We calculate relative errors in $L^2$ norm in percentage
\[
e(c) =  \sqrt{ \frac{\int_{\Omega}  
(c_{ms} - c)^2 \, dx}{\int_{\Omega} c^2 \,dx} } \cdot 100 \%, 
\quad 
e(u) = \sqrt{ \frac{\int_{\Omega}  (u_{ms} - u, u_{ms} - u) \,dx}{\int_{\Omega} (u, u) \,dx} }  \cdot 100 \%, 
\]
where $c_{ms}$ and $u_{ms}$ are multiscale solutions, $c$ and $u$ are reference solutions. 


\subsection{Test 1 (different geometries)}

\begin{figure}[h!]
\centering
\begin{subfigure}{1.0\textwidth}
\centering
\includegraphics[width=0.5\textwidth]{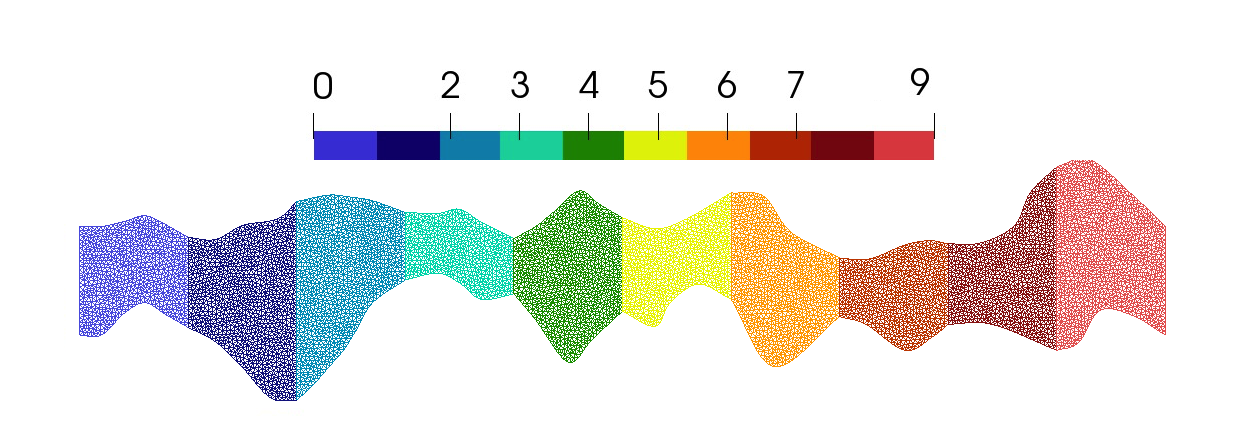}
\caption{Geometry 1 with 10 local domains}
\end{subfigure}
\begin{subfigure}{1.0\textwidth}
\centering
\includegraphics[width=0.5\textwidth]{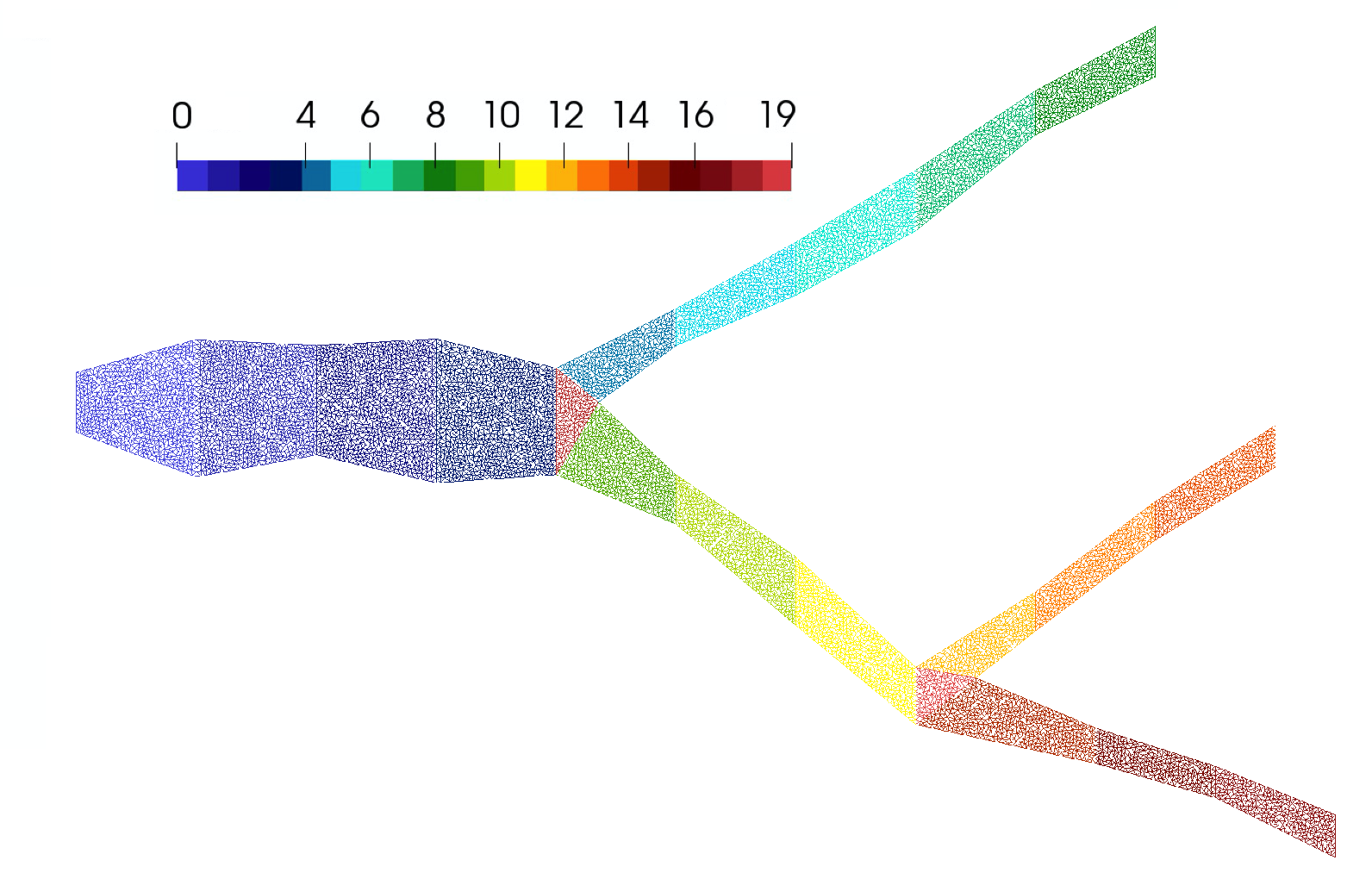}
\caption{Geometry 2 with 20 local domains}
\end{subfigure}
\begin{subfigure}{1.0\textwidth}
\centering
\includegraphics[width=0.45\textwidth]{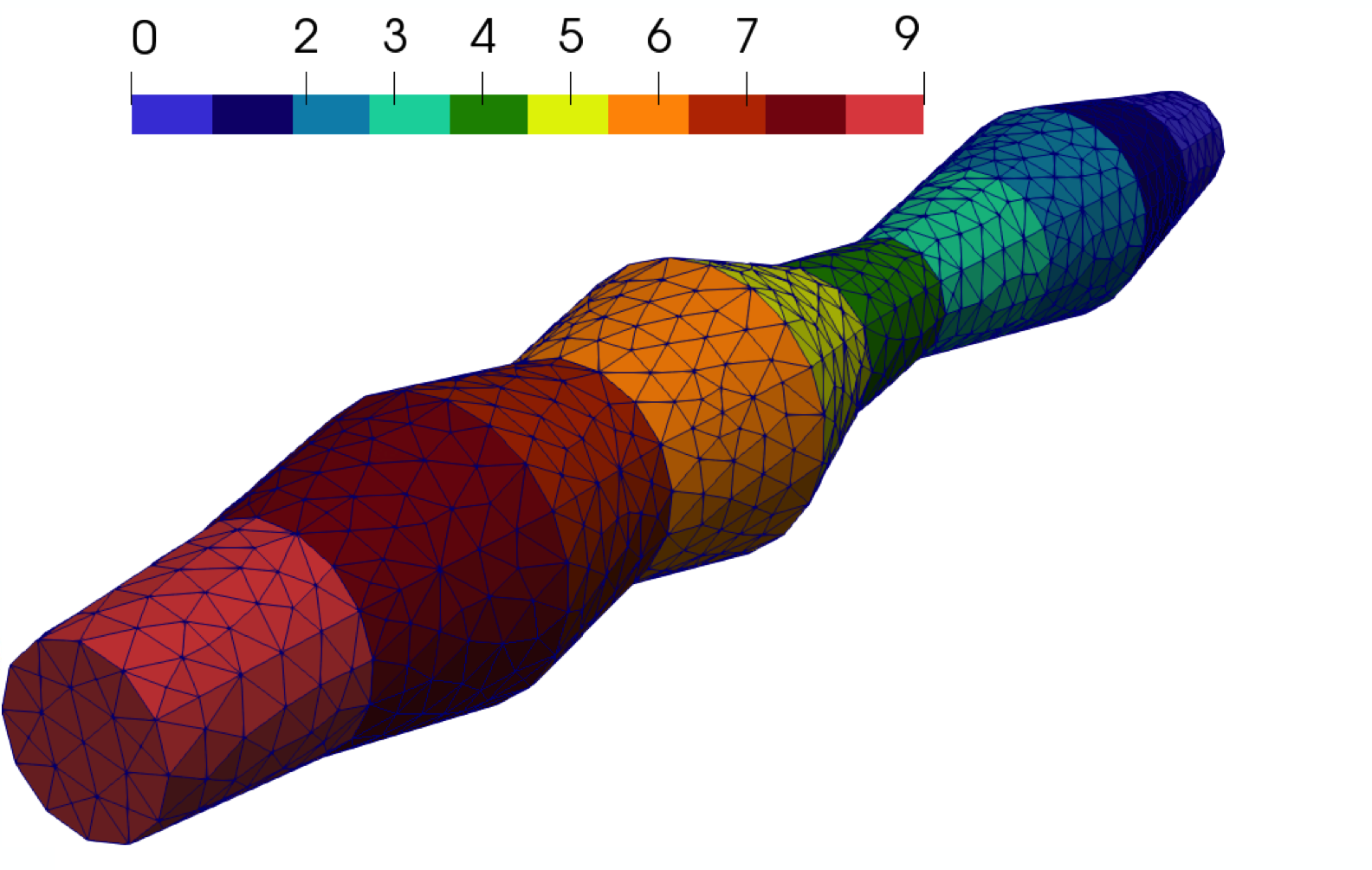}
\caption{Geometry 3 with 10 local domains}
\end{subfigure}
\caption{Fine grid with subdomain markers for Geometry 1, 2, and 3. }
\label{geom2}
\end{figure}

We consider a test problem with non-homogeneous Robin boundary condition for concentration
\[
- D \nabla c  \cdot n = \alpha (c - c_w), \quad x \in \Gamma_w, 
\]
where $c_w = 1$, $\alpha = 0.01$ and $D = 0.01$. 
As initial conditions, we set $c_0 = 1$ and $u_0 = 0$. 
For the inflow (left) boundary, we set $c_{in} = 0$ for $\Gamma_{in}$. 
We perform simulations for $t_{max} = 0.7$ (Geometry 1), $t_{max} = 1$ (Geometry 2)  and $t_{max} = 2$  (Geometry 3) with 40 time iterations. 
The coarse grid is structured with 20 local domains for Geometry 2, and with 10 local domains for Geometry 1 and 3. In Figure \ref{geom2}, we show local domain markers. 


\begin{figure}[h!]
\centering
\begin{subfigure}{1.0\textwidth}
\centering
\includegraphics[width=1.0\textwidth]{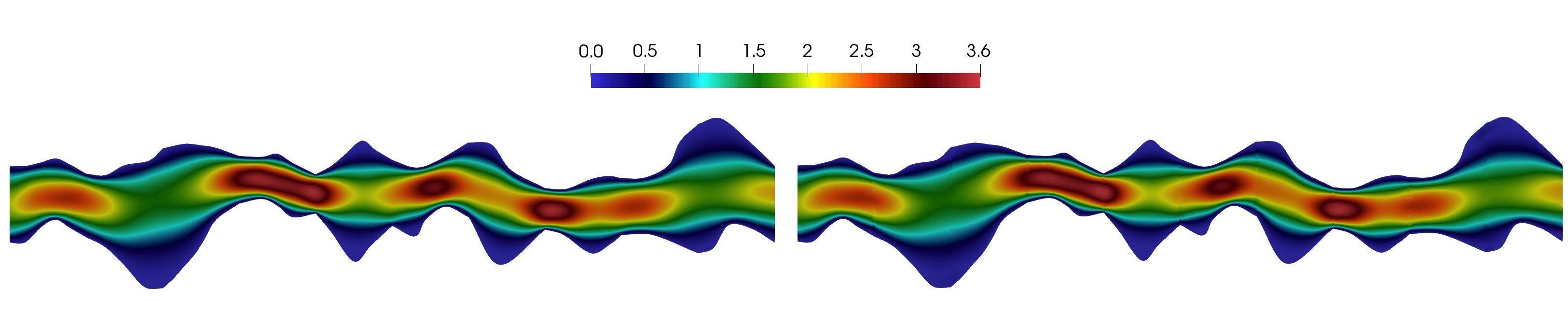}
\caption{Geometry 1. Left: reference solution, $DOF^u_h = 121 450$.  Right: multiscale solution, $DOF^u_H = 410$ }
\end{subfigure}
\begin{subfigure}{1.0\textwidth}
\centering
\includegraphics[width=1.0\textwidth]{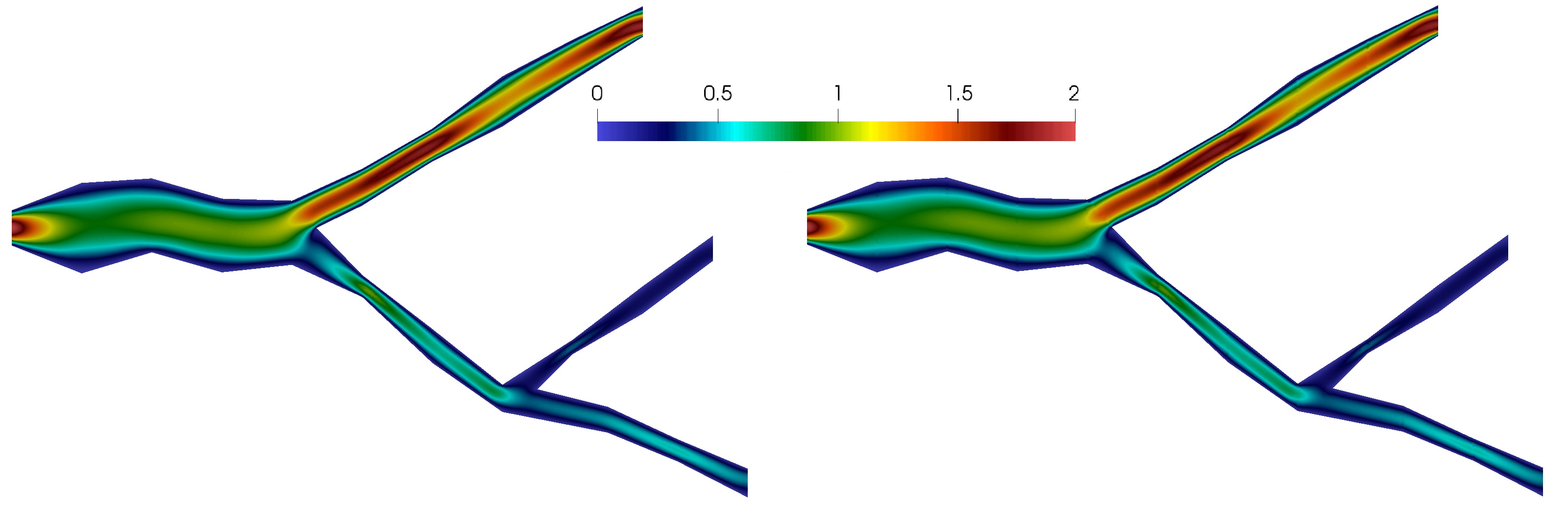}
\caption{Geometry 2. Left: reference solution, $DOF^u_h = 126 147$.  Right: multiscale solution, $DOF^u_H = 820$ }
\end{subfigure}
\begin{subfigure}{1.0\textwidth}
\centering
\includegraphics[width=1.0\textwidth]{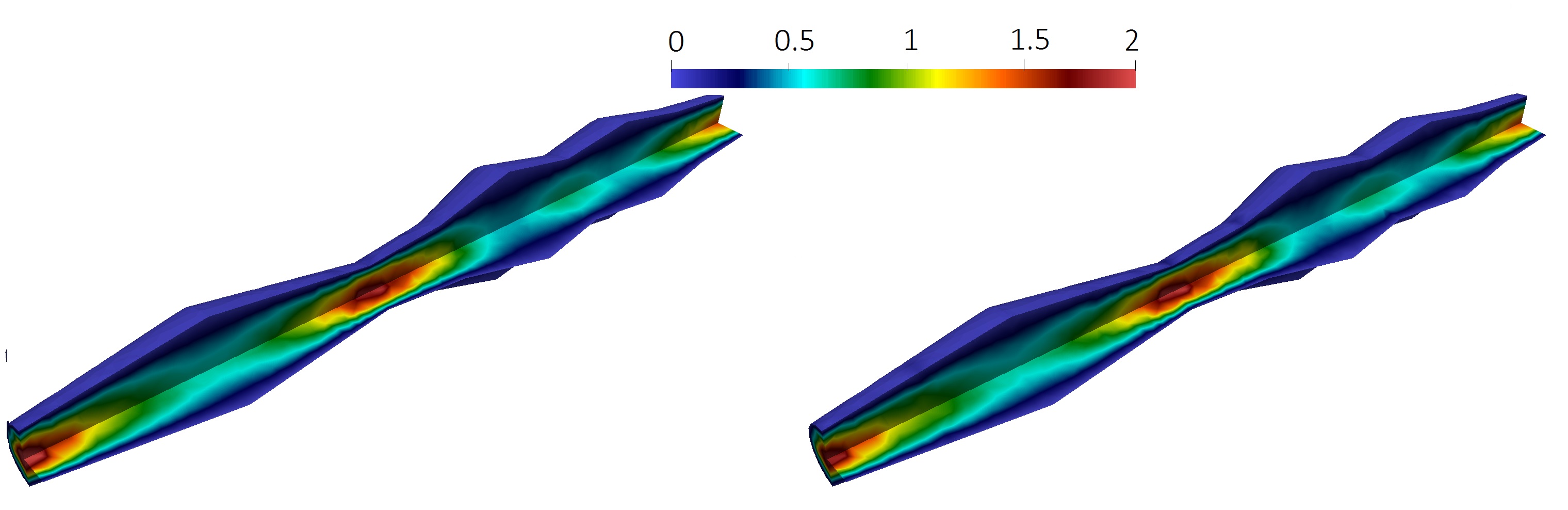}
\caption{Geometry 3. Left: reference solution, $DOF^u_h = 196 222$.  Right: multiscale solution, $DOF^u_H = 610$ }
\end{subfigure}
\caption{Reference and multiscale solutions of velocity (magnitude) at final time.  Geometry 1, 2 and 3 (Test 1). 
Left: reference solution.  Right: multiscale solution with 20 multiscale basis functions of Type 2. }
\label{sol-u}
\end{figure}

\begin{table}[h!]
\center
\begin{tabular}{|c|c|}
\hline
\multicolumn{2}{|c|}{Type 1} \\
\hline
$DOF_H^u (M^u)$ 
& $e(u)$ (\%) \\
\hline
\multicolumn{2}{|c|}{Geometry 1}  \\ 
\hline
110 (10)	& 11.70 \\
210 (20)	& 6.013 \\
410 (40)	& 1.536 \\
610 (60)	& 1.086 \\
810 (80)	& 1.064 \\
\hline
\multicolumn{2}{|c|}{Geometry 2}  \\ 
\hline
220 (10)	& 27.20 \\
420 (20)	& 6.205 \\
820 (40)	& 1.485 \\
1220 (60)	& 1.369 \\
1620 (80)	& 1.374 \\
\hline
\multicolumn{2}{|c|}{Geometry 3}  \\ 
\hline
160 (15)		& 36.70 \\
310 (30)		& 10.76 \\
610 (60)		& 9.493 \\
910 (90)		& 8.548 \\
1210 (120)	& 5.657 \\
\hline
\end{tabular}
\,\,
\begin{tabular}{|c|c|}
\hline
\multicolumn{2}{|c|}{Type 2} \\
\hline
$DOF_c^u (M^u)$ & $e(u)$ (\%) \\
\hline
\multicolumn{2}{|c|}{Geometry 1}  \\ 
\hline
110 (5)		& 19.85 \\
210 (10)	& 10.52 \\
410 (20)	& 3.230 \\
610 (30)	& 1.756 \\
810 (40)	& 1.346 \\
\hline
\multicolumn{2}{|c|}{Geometry 2}  \\ 
\hline
220 (5)		& 31.51 \\
420 (10)	& 10.08  \\
820 (20)	& 3.570 \\
1220 (30)	& 1.906 \\
1620 (40)	& 1.521  \\
\hline
\multicolumn{2}{|c|}{Geometry 3}  \\ 
\hline
160 (5)			& 15.96 \\
310 (10)		& 10.63 \\
610 (20)		& 6.548 \\
910 (30)		& 4.167 \\
1210 (40)		& 4.082 \\
\hline
\end{tabular}
\caption{Relative $L_2$ error for the velocity at the final time. Geometry 1, 2, and 3 (Test 1).
Left: Type 1 multiscale basis functions. 
Right: Type 2 multiscale basis functions. 
Reference solution with  
$DOF^u_h = 121 450$ (Geometry 1), 
$DOF^u_h = 126 147$ (Geometry 2) and 
$DOF^u_h = 196 222$ (Geometry 3). }
\label{err-u-g}
\end{table}


\begin{figure}[h!]
\centering
\includegraphics[width=0.3\textwidth]{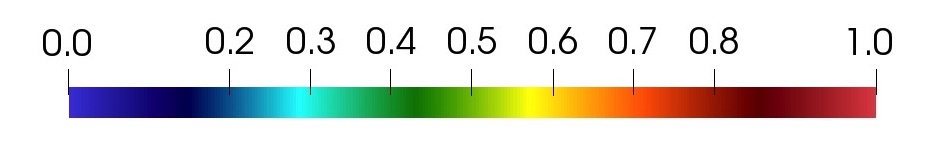}
\begin{subfigure}{1.0\textwidth}
\centering
\includegraphics[width=1.0\textwidth]{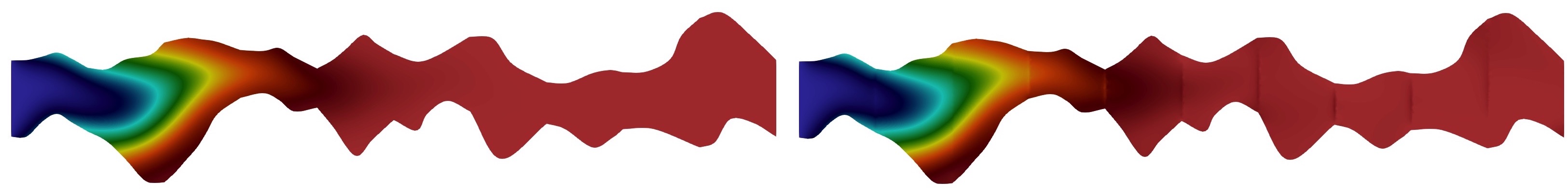}
\caption{Concentration, $c_{10}$. Left: reference solution. Right: multiscale solution}
\end{subfigure}
\begin{subfigure}{1.0\textwidth}
\centering
\includegraphics[width=1.0\textwidth]{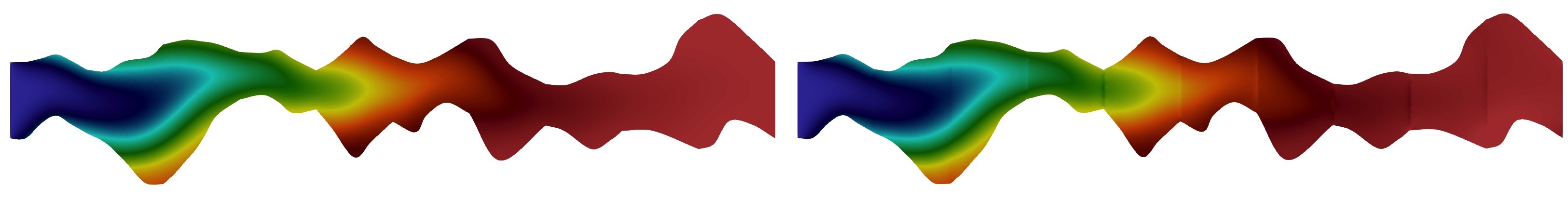}
\caption{Concentration, $c_{20}$. Left: reference solution. Right: multiscale solution}
\end{subfigure}
\begin{subfigure}{1.0\textwidth}
\centering
\includegraphics[width=1.0\textwidth]{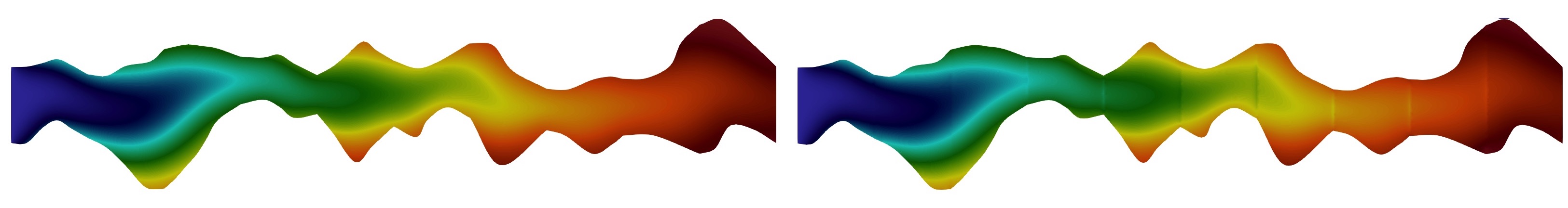}
\caption{Concentration, $c_{40}$. Left: reference solution. Right: multiscale solution}
\end{subfigure}
\caption{Reference and multiscale solutions of concentration at $t_m$ for $m = 10, 20$ and $40$.   
Geometry 1 (Test 1). 
Left: reference solution, $DOF^c_h = 52 050$. 
Right: multiscale solution with 20 multiscale basis functions of Type 2,  $DOF^c_H = 410$.}
\label{sol-c-g1}
\end{figure}

\begin{figure}[h!]
\centering
\includegraphics[width=0.3\textwidth]{bar01}
\begin{subfigure}{1.0\textwidth}
\centering
\includegraphics[width=1.0\textwidth]{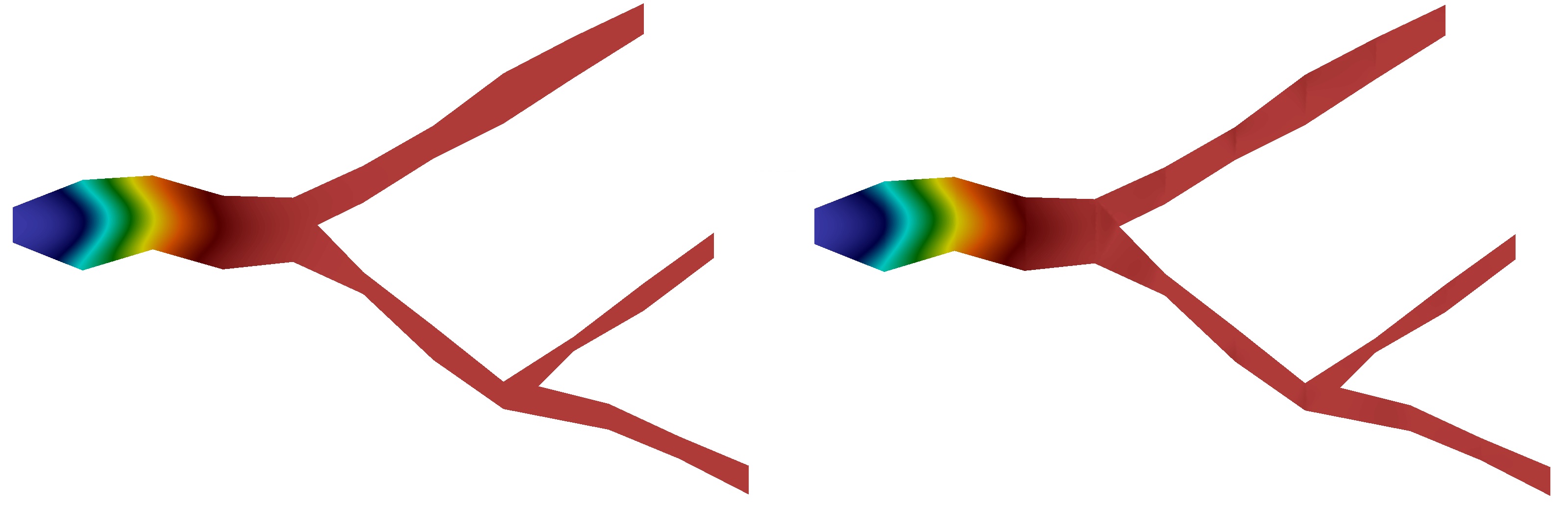}
\caption{Concentration, $c_{10}$. Left: reference solution. Right: multiscale solution}
\end{subfigure}
\begin{subfigure}{1.0\textwidth}
\centering
\includegraphics[width=1.0\textwidth]{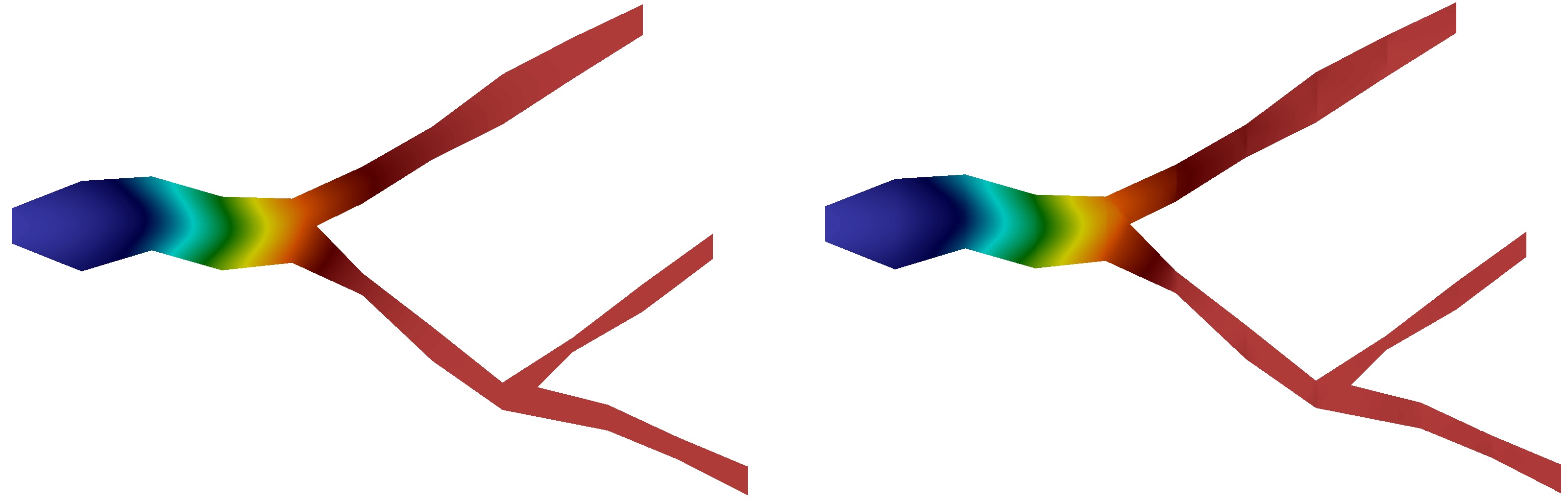}
\caption{Concentration, $c_{20}$. Left: reference solution. Right: multiscale solution}
\end{subfigure}
\begin{subfigure}{1.0\textwidth}
\centering
\includegraphics[width=1.0\textwidth]{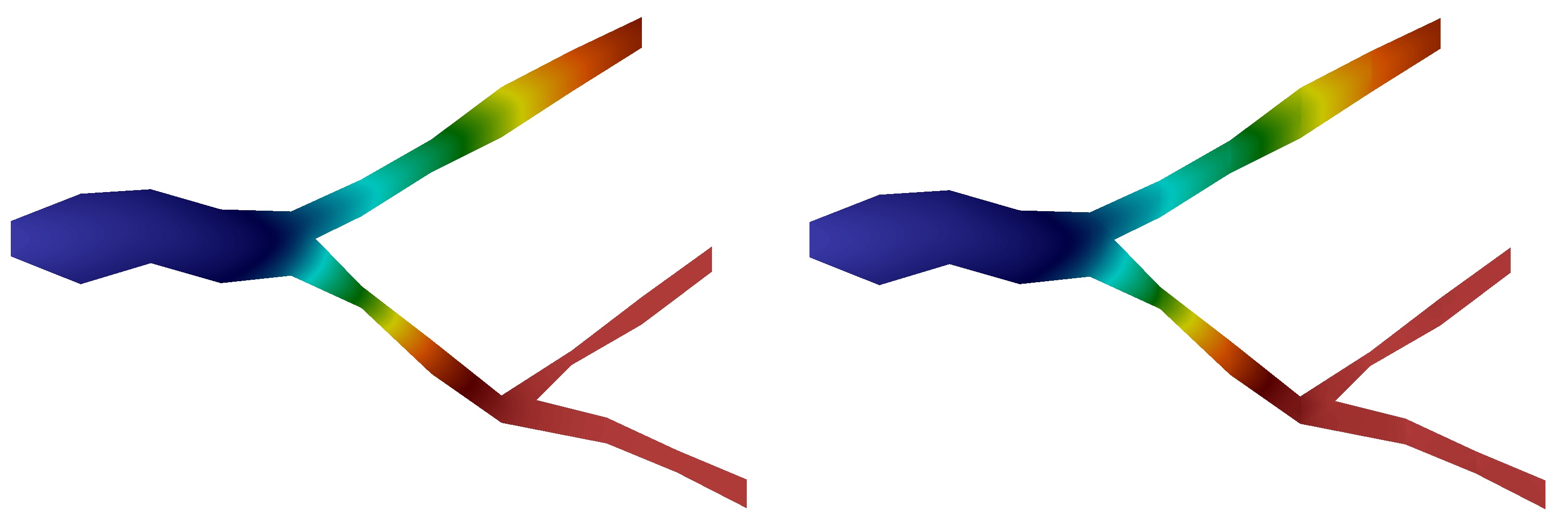}
\caption{Concentration, $c_{40}$. Left: reference solution. Right: multiscale solution}
\end{subfigure}
\caption{Reference and multiscale solutions of concentration at $t_m$ for $m = 10, 20$ and $40$. 
Geometry 2  (Test 1). 
Left: reference solution, $DOF^c_h = 54 063$. 
Right: multiscale solution with 20 multiscale basis functions of Type 2,  $DOF^c_H = 820$.}
\label{sol-c-g2}
\end{figure}

\begin{figure}[h!]
\centering
\includegraphics[width=0.3\textwidth]{bar01}
\begin{subfigure}{1.0\textwidth}
\centering
\includegraphics[width=1.0\textwidth]{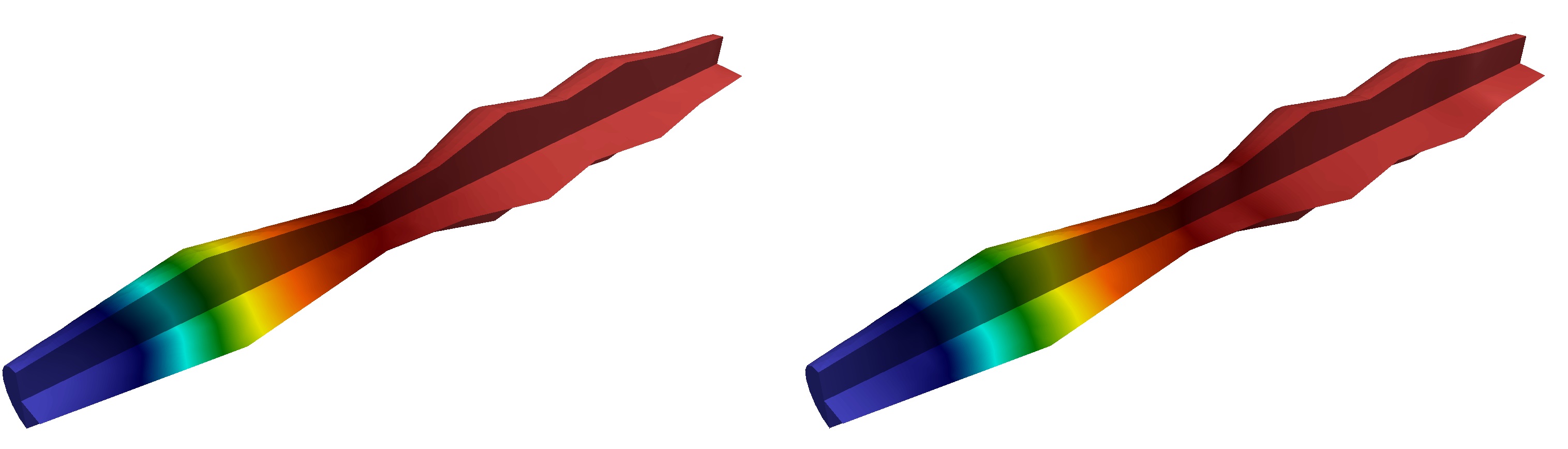}
\caption{Concentration, $c_{10}$. Left: reference solution. Right: multiscale solution}
\end{subfigure}
\begin{subfigure}{1.0\textwidth}
\centering
\includegraphics[width=1.0\textwidth]{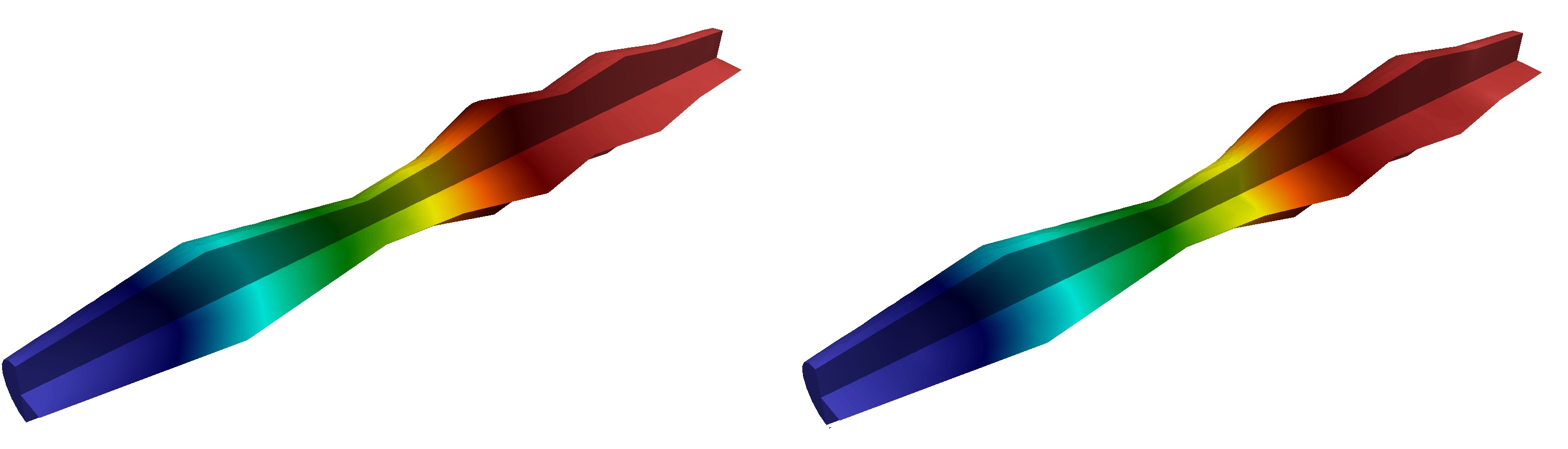}
\caption{Concentration, $c_{20}$. Left: reference solution. Right: multiscale solution}
\end{subfigure}
\begin{subfigure}{1.0\textwidth}
\centering
\includegraphics[width=1.0\textwidth]{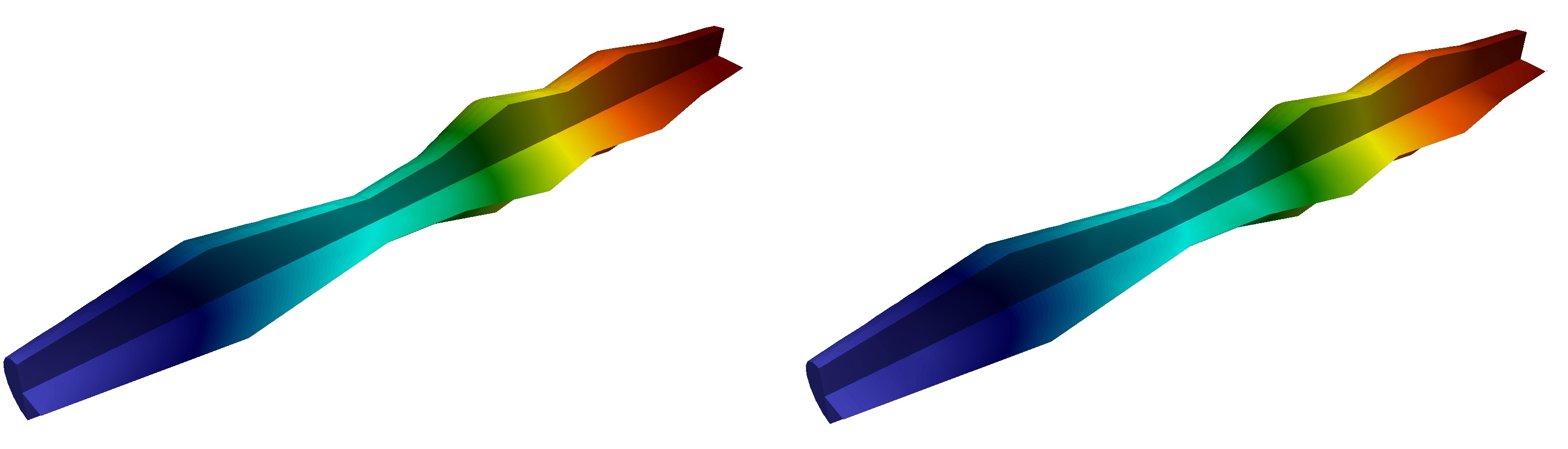}
\caption{Concentration, $c_{40}$. Left: reference solution. Right: multiscale solution}
\end{subfigure}
\caption{Reference and multiscale solutions of concentration at $t_m$ for $m = 10, 20$ and $40$.  
 Geometry 3 (Test 1). 
Left: reference solution, $DOF^c_h = 60 376$. 
Right: multiscale solution with 20 multiscale basis functions of Type 2,  $DOF^c_H = 410$.}
\label{sol-c-g3}
\end{figure}

\begin{table}[h!]
\center
\begin{tabular}{ | c | c  | c | c | c |}
\hline
\multicolumn{5}{|c|}{Type 1, without velocity} \\
\hline
$DOF_H^c (M^c)$  
& $e(c_{10}$) & $e(c_{20}$)  
& $e(c_{30}$) & $e(c_{40}$) \\
\hline
\multicolumn{5}{|c|}{Geometry 1} \\
\hline
30 (2)		& 42.27	& 38.93	& 37.68	& 37.64 \\
70 (6)		& 6.747	& 7.675	& 8.587	& 9.206 \\
110 (10)	& 3.160	& 3.818	& 4.269	& 4.549 \\
210 (20)	& 1.954	& 2.228	& 2.294	& 2.300 \\
410 (40)	& 1.924	& 2.189	& 2.252	& 2.259 \\ 
610 (60)	& 1.920	& 2.184	& 2.247	& 2.253 \\
810 (80)	& 1.916	& 2.177	& 2.236	& 2.240 \\
\hline
\multicolumn{5}{|c|}{Geometry 2} \\
\hline
60 (2)		& 60.30  	& 61.20  	& 55.67  	& 46.25 \\
140 (6)		& 1.931  	& 2.021  	& 2.668 	& 2.766 \\
220 (10)	& 1.784 	& 1.495  	& 1.471 	& 1.451 \\
420 (20)	& 1.041  	& 1.074   	& 1.059 	& 0.989 \\
820 (40)	& 1.048 	& 1.315 	& 1.069 	& 0.997 \\ 
1220 (60)	& 1.054   	& 1.085   	& 1.075   	& 1.004 \\
1620 (80)	& 1.058   	& 1.088  	& 1.079   	& 1.009 \\
\hline
\multicolumn{5}{|c|}{Geometry 3} \\
\hline
30 (2)		& 16.17  	& 16.07   	& 21.69  	& 23.32 \\
70 (6)		& 16.60  	& 14.86   	& 15.52 	& 15.85 \\
110 (10)	& 16.23 	& 14.50  	& 15.15 	& 15.47 \\
210 (20)	& 5.274  	& 4.940   	& 5.252 	& 5.648 \\
410 (40)	& 4.871 	& 4.607 	& 4.947 	& 5.338 \\ 
610 (60)	& 4.623   	& 4.391 	& 4.715   	& 5.094 \\
810 (80)	& 4.568   	& 4.337  	& 4.659   	& 5.037 \\
\hline
\end{tabular}
\,\,
\begin{tabular}{ | c | c  | c | c | c |}
\hline
\multicolumn{5}{|c|}{Type 2, without velocity} \\
\hline
$DOF_c^c (M^c)$  
& $e(c_{10}$) & $e(c_{20}$)  
& $e(c_{30}$) & $e(c_{40}$) \\
\hline
\multicolumn{5}{|c|}{Geometry 1} \\
\hline
30 (1)		& 46.99	& 53.03	& 52.54	 & 50.25 \\
70 (3)		& 4.111	& 4.070	& 4.746	& 5.154 \\
110 (5)		& 2.653	& 2.851	& 2.973	& 2.989 \\
210 (10)	& 1.950	& 2.196	& 2.236	& 2.224 \\
410 (20)	& 1.906	& 2.127	& 2.154	& 2.140 \\ 
610 (30)	& 1.902	& 2.121	& 2.148	& 2.135 \\
810 (40)	& 1.901	& 2.120	& 2.147	& 2.134 \\
\hline
\multicolumn{5}{|c|}{Geometry 2} \\
\hline
60 (1)		& 73.44  	& 84.91 	& 87.34   & 87.06 \\
140 (3)		& 2.256   	& 2.205   	& 2.349  	& 2.249 \\
220 (5)		& 1.341  	& 1.569   	& 1.692  	& 1.581 \\
420 (10)	& 1.142   	& 1.204  	& 1.228 	& 1.144 \\
820 (20)	& 1.145  	& 1.190 	& 1.215  	& 1.144 \\
1220 (30)	& 1.148   	& 1.191   	& 1.218   	& 1.149 \\
1620 (40)	& 1.149   	& 1.191  	& 1.219   	& 1.150 \\ 
\hline
\multicolumn{5}{|c|}{Geometry 3} \\
\hline
30 (1)		& 58.70   	& 71.74 	& 74.82   	& 73.37 \\
70 (3)		& 2.322  	& 1.866 	& 2.071   	& 2.229 \\
110 (5)		& 2.306   	& 1.895   	& 2.051  	& 2.214 \\
210 (10)	& 2.317  	& 2.142  	& 2.348 	& 2.586 \\
410 (20)	& 1.892  	& 1.773 	& 1.992  	& 2.210 \\ 
610 (30)	& 1.691   	& 1.595   & 1.792   	& 1.988 \\
810 (40)	& 1.651  	& 1.552  	& 1.746   	& 1.938 \\
\hline
\end{tabular}
\caption{Relative $L_2$ error for concentration at $t_m$ ($m=10, 20, 30$ and $40$). Geometry 1, 2 and 3 (Test 1). 
Multiscale basis functions without velocity. 
Left: Type 1 multiscale basis functions for concentration. 
Right: Type 2 multiscale basis functions for concentration. 
Reference solution with  $DOF^c_h = 52 050$ (Geometry 1), $DOF^c_h = 54 063$ (Geometry 2) and $DOF^c_h = 60 376$ (Geometry 3). }
\label{err-c-g-cu0}
\end{table}

\begin{table}[h!]
\center
\begin{tabular}{ | c | c  | c | c | c |}
\hline
\multicolumn{5}{|c|}{Type 1, with velocity} \\
\hline
$DOF_H^c (M^c)$  
& $e(c_{10}$) & $e(c_{20}$)  
& $e(c_{30}$) & $e(c_{40}$) \\
\hline
\multicolumn{5}{|c|}{Geometry 1} \\
\hline
30 (2)		& 60.48 	& 66.48	& 65.71	& 63.81 \\
70 (6)		& 41.95 	& 45.45	& 43.60	& 40.94 \\
110 (10)	& 31.93 	& 34.33	& 31.66	& 28.18 \\
210 (20)	& 6.114	& 7.652	& 6.285	& 4.931 \\
410 (40)	& 1.747	& 1.776	& 1.572	& 1.376 \\ 
610 (60)	& 1.742	& 1.762	& 1.555	& 1.359 \\
810 (80)	& 1.741	& 1.758	& 1.549	& 1.352 \\
\hline
\multicolumn{5}{|c|}{Geometry 2} \\
\hline
60 (2)		& 65.49  	& 64.93  	& 58.22  	& 48.53 \\
140 (6)		& 47.06  	& 53.16   	& 53.16 	& 42.65 \\
220 (10)	& 16.22 	& 15.71  	& 15.73 	& 14.42 \\
420 (20)	& 1.467  	& 1.491   	& 1.474 	& 1.275 \\
820 (40)	& 1.277 	& 1.315 	& 1.274 	& 1.089 \\ 
1220 (60)	& 1.262   	& 1.302 	& 1.259  	& 1.075 \\
1620 (80)	& 1.252   	& 1.293 	& 1.248  	& 1.065 \\
\hline
\multicolumn{5}{|c|}{Geometry 3} \\
\hline
30 (2)		& 43.35  	& 39.21 	& 40.01 	& 40.01 \\
70 (6)		& 22.71   	& 21.47   	& 21.47  	& 21.47 \\
110 (10)	& 21.90  	& 20.70  	& 21.68  	& 20.55 \\
210 (20)	& 15.56   	& 14.70  	& 15.34 	& 15.34 \\
410 (40)	& 12.48   	& 12.05 	& 12.96  	& 13.77 \\ 
610 (60)	& 11.50   	& 11.24   	& 12.14   	& 12.94 \\
810 (80)	& 11.20   	& 10.98  	& 11.88   	& 12.68 \\
\hline
\end{tabular}
\,\,
\begin{tabular}{ | c | c  | c | c | c |}
\hline
\multicolumn{5}{|c|}{Type 2, with velocity} \\
\hline
$DOF_c^c (M^c)$  
& $e(c_{10}$) & $e(c_{20}$)  
& $e(c_{30}$) & $e(c_{40}$) \\
\hline
\multicolumn{5}{|c|}{Geometry 1} \\
\hline
30 (1)		& 49.91	& 53.68	& 52.64	& 50.36 \\
70 (3)		& 26.22	& 24.27	& 21.11	& 19.47 \\
110 (5)		& 14.63	& 13.24	& 9.694	& 7.351 \\
210 (10)	& 6.909	& 5.760	& 3.943	& 2.372 \\
410 (20)	& 1.867	& 2.004	& 1.870	& 1.701 \\ 
610 (30)	& 1.850	& 1.979	& 1.842	& 1.674 \\
810 (40)	& 1.846	& 1.973	& 1.836	& 1.669 \\
\hline
\multicolumn{5}{|c|}{Geometry 2} \\
\hline
60 (1)		& 48.37  	& 53.67 	& 52.88   	& 49.64 \\
140 (3)		& 14.56  	& 19.78   	& 20.14  	& 15.63 \\
220 (5)		& 12.46  	& 16.05   	& 16.30  	& 12.98 \\
420 (10)	& 1.474   	& 1.517  	& 1.532 	& 1.350 \\
820 (20)	& 1.308  	& 1.362 	& 1.340  	& 1.159 \\
1220 (30)	& 1.282  	& 1.338 	& 1.310   	& 1.133 \\
1620 (40)	& 1.272  	& 1.331 	& 1.302   	& 1.126 \\ 
\hline
\multicolumn{5}{|c|}{Geometry 3} \\
\hline
30 (1)		& 32.10  	& 32.21 	& 34.24   	& 32.98 \\
70 (3)		& 32.09  	& 32.19 	& 34.22   	& 32.96 \\
110 (5)		& 10.03   	& 9.976   	& 10.90  	& 11.25 \\
210 (10)	& 7.491   	& 7.609  	& 8.511 	& 8.736 \\
410 (20)	& 6.934  	& 7.026 	& 7.702  	& 8.124 \\ 
610 (30)	& 6.189  	& 6.384   	& 7.065   	& 7.484 \\
810 (40)	& 6.028   	& 6.234  	& 6.909   	& 7.320 \\
\hline
\end{tabular}
\caption{Relative $L_2$ error for concentration at $t_m$ ($m=10, 20, 30$ and $40$). Geometry 1, 2 and 3 (Test 1). 
Multiscale basis functions with velocity. 
Left: Type 1 multiscale basis functions for concentration. 
Right: Type 2 multiscale basis functions for concentration. 
Reference solution with  $DOF^c_h = 52 050$ (Geometry 1), $DOF^c_h = 54 063$ (Geometry 2) and $DOF^c_h = 60 376$ (Geometry 3). }
\label{err-c-g-cu1}
\end{table}


We consider the performance of the presented method for the solution of transport and flow problems in different computational domains. We consider two types of multiscale basis functions (see Section \ref{s-cg}). 
$M^c$ and $M^u$ are the number of the multiscale basis functions for concentration and velocity, respectively. 
$DOF^u_h$  and $DOF^c_h$  are the number of degrees of freedom for reference (fine grid) solution. 
$DOF^u_H$ and $DOF^c_H$ are the numbers of degrees of freedom for a multiscale solution. 
The reference solution is performed on the fine grid using IPDG approximation presented in Section \ref{s-mmfg}. We used linear basis functions for both velocity and concentration fields on the fine grid. Therefore, 
$DOF^u_h = N_{cell}^h \cdot (d \, (d+1)  + 1) $ and  
$DOF^c_h = N_{cell}^h  \cdot (d+1)$. 
For multiscale solver, we have 
$DOF^u_H = N_{cell}^H (M^u + 1)$,   
$DOF^c_H = N_{cell}^H (M^c + 1)$ for Type 1 and 
$DOF^u_H = N_{cell}^H (d \cdot M^u + 1) $,   
$DOF^c_H = N_{cell}^H (2 \cdot M^c + 1)$ for Type 2 multiscale basis functions.

Because velocity does not depend on the concentration in our test problems. We start with results for flow problems.
In Figure \ref{sol-u}, we present reference and multiscale solutions for Geometries 1, 2, and 3. We depicted the magnitude of the velocity field at the final time.  In a multiscale solver, we used $M^u = 20$ multiscale basis functions for each direction (Type 2). 
For Geometry 1, we have $DOF^u_h = 121 450$ and $DOF^u_H = 410$ (0.33 \% from $DOF^u_h$). 
For Geometry 2, we have $DOF^u_h = 126 147$ and $DOF^u_H = 820$ (0.65 \% from $DOF^u_h$). 
For Geometry 3, we have $DOF^u_h = 196 222$ and $DOF^u_H = 610$ (0.31 \% from $DOF^u_h$).
We observe a good accuracy of the presented method with huge reduction of the discrete system size. 

In Table \ref{err-u-g}, we present relative errors in \% for velocity between reference solution and multiscale solution with different numbers of the multiscale basis functions at the final times.
We observe a reduction of the error with an increasing number of multiscale basis functions for all geometries. 
For example for 40 multiscale basis functions of Type 2, we have 
$1.3$ \% of error for Geometry 1,  
$1.5$ \% of error for Geometry 2, and 
$4.0$ \% of error for Geometry 3. 
For Type 1 and 2, we observe similar errors for two-dimensional problems (Geometry 1 and 2). For three - dimensional domain in Geometry 3, we obtain better results with Type 2 multiscale basis functions.   

In Figures \ref{sol-c-g1}, \ref{sol-c-g2} and \ref{sol-c-g3}, we present concentration distributions for the reference (fine scale) and multiscale solutions at different time layers $t_m$ for $m = 10, 20$ and $40$. 
In these calculations, we used a fixed number of multiscale basis functions for the velocity field ($M^u = 20$ of Type 2).  We will investigate the influence of the velocity accuracy on the concentration solution later in Test 3. 
In figures, we presented results for $M^u = 20$ multiscale basis functions for each type of local domain boundaries (Type 2). 
For Geometry 1 (Figure \ref{sol-c-g1}), we have $DOF^c_h = 52 050$ and $DOF^c_H = 410$ (0.78 \% from $DOF^c_h$). 
For Geometry 2 Figure \ref{sol-c-g2}), we have $DOF^c_h = 54 063$ and $DOF^c_H = 820$ (1.5 \% from $DOF^c_h$). 
For Geometry 3 Figure \ref{sol-c-g3}), we have $DOF^c_h = 60 376$ and $DOF^c_H = 410$ (0.67 \% from $DOF^c_h$). 
We observe good results of the presented method for solving transport problems for all three geometries. 

We present numerical results in Tables \ref{err-c-g-cu0}  and  \ref{err-c-g-cu1}  for different number of multiscale basis functions for concentration (Type 1 and 2) for a fixed number of multiscale basis functions for velocity field ($M^u = 20$ of Type 2). Relative errors for concentrations are presented for four time layers $t_m$ with $m = 10, 20, 30$ and $40$. 
We can obtain good multiscale solution when we take a sufficient number of multiscale basis functions for pressure and for displacements. 
For example, we obtain near 40 \% of concentration error, when we take 3 multiscale basis functions of Type 1 and near 20 \% for Type 2. 
For 40 multiscale basis functions of Type 2, the error reduce to
$1.6$ \% for Geometry 1,  
$1.1$ \%  for Geometry 2, and 
$1.9$ \%  for Geometry 3. 
For two-dimensional problems (Geometry 1 and 2), we obtain a similar error for Type 1 and 2 multiscale basis functions. However, for three - dimensional problem (Geometry 3), the multiscale basis of Type 2 is better.



We observe that the presented multiscale method provides good results with small errors and huge reduction of the system size for all three test geometries. 
For the presented test problem in the 3D case is better to use basis functions without information on the velocity field. For the 2D result, we observe similar errors for basis with and without velocity. 
For three-dimensional problem Type 2 multiscale basis functions are better than Type 1 due to direct definitions of the flow and transport directions. However, for 2D problem, we obtain similar results for Type 1 and 2 basis functions.

\subsection{Test 2  (different boundary conditions and diffusion coefficients)}

Next, we consider the efficiency of the presented method for the solution of the transport problem with different types of boundary conditions and different values of diffusion coefficient $D$. 

\begin{figure}[h!]
\centering
\includegraphics[width=0.3\textwidth]{bar01}
\begin{subfigure}{1.0\textwidth}
\centering
\includegraphics[width=1\textwidth]{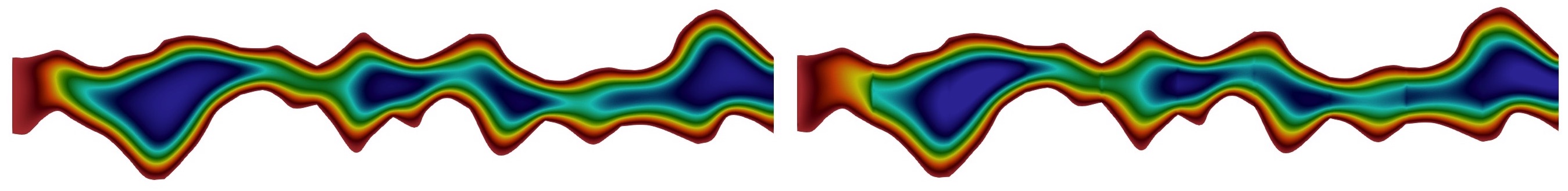}
\caption{Concentration, $c_{10}$. Left: reference solution. Right: multiscale solution}
\end{subfigure}
\begin{subfigure}{1.0\textwidth}
\centering
\includegraphics[width=1\textwidth]{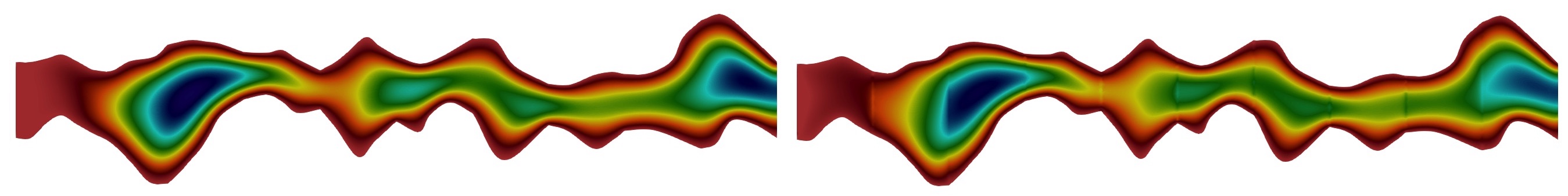}
\caption{Concentration, $c_{20}$. Left: reference solution. Right: multiscale solution}
\end{subfigure}
\begin{subfigure}{1.0\textwidth}
\centering
\includegraphics[width=1\textwidth]{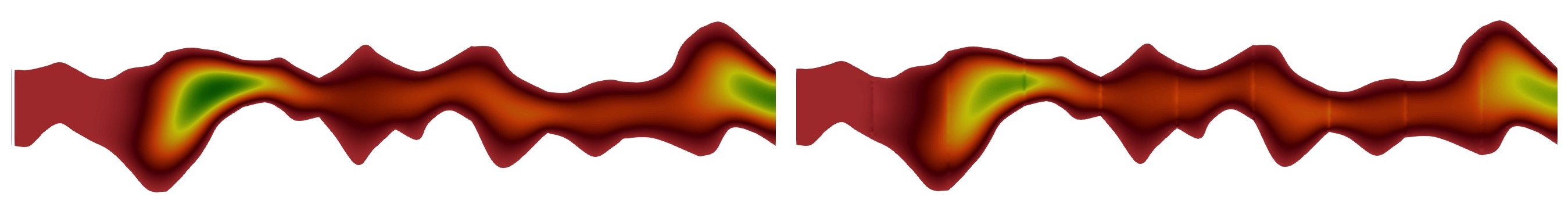}
\caption{Concentration, $c_{40}$. Left: reference solution. Right: multiscale solution}
\end{subfigure}
\caption{Reference and multiscale solutions of concentration at $t_m$ for $m = 10, 20$ and $40$.   
Dirichlet boundary conditions (Test 2).   
Left: reference solution. 
Right: multiscale solution with 20 multiscale basis functions of Type 2.}
\label{sol-dbc}
\end{figure}

\begin{figure}[h!]
\centering
\includegraphics[width=0.3\textwidth]{bar01}
\begin{subfigure}{1.0\textwidth}
\centering
\includegraphics[width=1\textwidth]{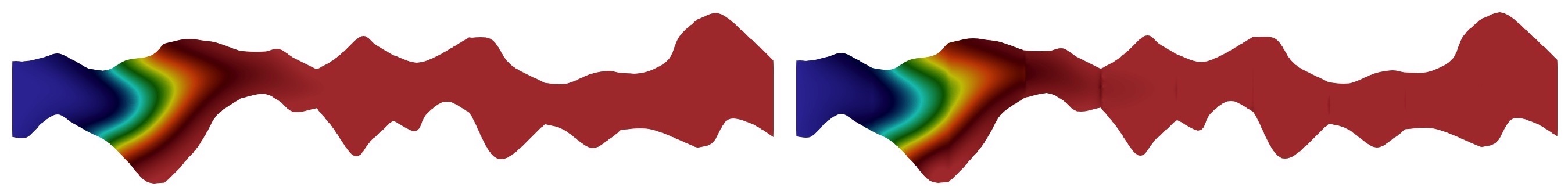}
\caption{Concentration, $c_{10}$. Left: reference solution. Right: multiscale solution}
\end{subfigure}
\begin{subfigure}{1.0\textwidth}
\centering
\includegraphics[width=1\textwidth]{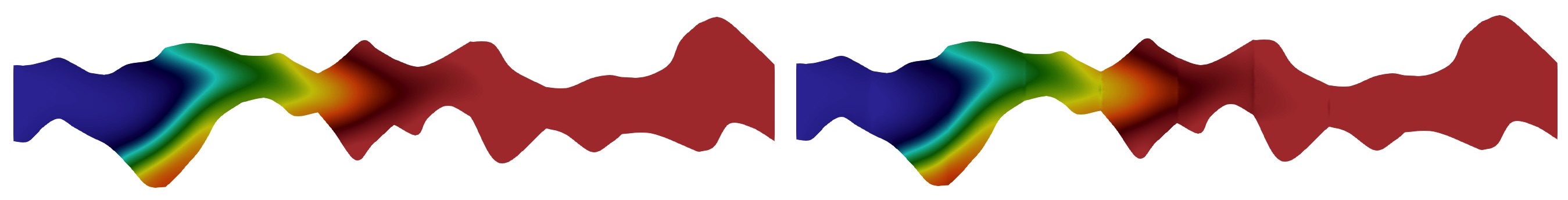}
\caption{Concentration, $c_{20}$. Left: reference solution. Right: multiscale solution}
\end{subfigure}
\begin{subfigure}{1.0\textwidth}
\centering
\includegraphics[width=1\textwidth]{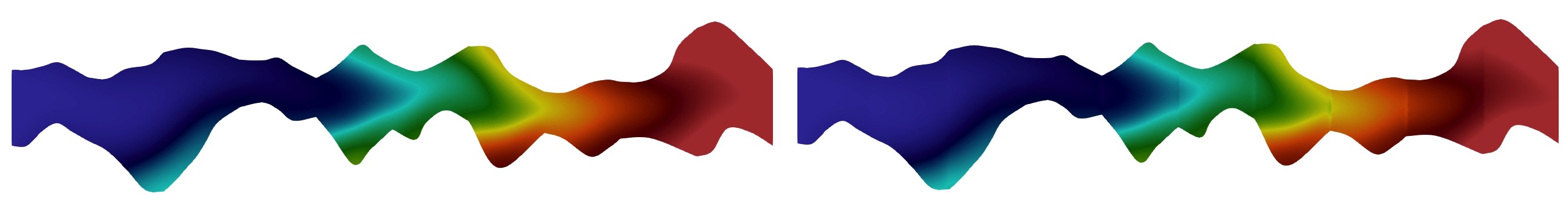}
\caption{Concentration, $c_{40}$. Left: reference solution. Right: multiscale solution}
\end{subfigure}
\caption{Reference and multiscale solutions of concentration at $t_m$ for $m = 10, 20$ and $40$.  
Neumann boundary conditions (Test 2).   
Left: reference solution. 
Right: multiscale solution with 20 multiscale basis functions of Type 2.}
\label{sol-nbc}
\end{figure}

\begin{table}[h!]
\center
\begin{tabular}{ | c | c  | c | c | c | }
\hline
\multicolumn{5}{|c|}{with velocity}  \\ 
\hline
$DOF^c_H (M^c)$  
& $e(c_{10}$) & $e(c_{20}$)  
& $e(c_{30}$) & $e(c_{40}$) \\
\hline
\multicolumn{5}{|c|}{Type 1, without time} \\
\hline
30 (2)		& 51.74	& 36.78	& 26.43	& 18.88 \\
70 (6)		& 51.75	& 36.82	& 26.45	& 18.90 \\
110 (10)	& 44.35	& 31.52	& 22.54	& 16.19 \\
210 (20)	& 18.04	& 11.93	& 8.893	& 6.877 \\
410 (40)	& 9.781	& 5.784	& 4.241	& 3.229 \\ 
610 (60)	& 9.205	& 5.434	& 3.993	& 3.036 \\
810 (80)	& 9.114	& 5.366	& 3.946	& 3.001 \\
\hline
\multicolumn{5}{|c|}{Type 2,  without time} \\
\hline
30 (1)		& 44.40	& 31.76	& 23.15	& 16.78 \\
70 (3)		& 17.92	& 11.85	& 8.891	& 6.782 \\
110 (5)		& 12.42	& 6.768	& 4.681	& 3.486 \\
210 (10)	& 9.150	& 5.210	& 3.815	& 2.906 \\
410 (20)	& 8.845	& 5.065	& 3.848	& 2.888 \\ 
610 (30)	& 8.831	& 5.132	& 3.787	& 2.887 \\
810 (40)	& 8.829	& 5.131	& 3.786	& 2.886  \\
\hline
\multicolumn{5}{|c|}{Type 2, with time} \\
\hline
30 (1)		& 44.40	& 31.76	& 23.15	& 16.78 \\
70 (3)		& 17.92	& 11.86	& 8.892	& 6.783 \\
110 (5)		& 12.42	& 6.768	& 4.681	& 3.487 \\
210 (10)	& 9.150	& 5.210	& 3.815	& 2.906  \\
410 (20)	& 8.844	& 5.135	& 3.788	& 2.888  \\ 
610 (30)	& 8.831	& 5.132	& 3.787	& 2.887 \\
810 (40)	& 8.829	& 5.131	& 3.786	& 2.886 \\
\hline
\end{tabular}
\,\,
\begin{tabular}{ | c | c  | c | c | c | }
\hline
\multicolumn{5}{|c|}{without velocity}  \\ 
\hline
$DOF^c_H (M^c)$  
& $e(c_{10}$) & $e(c_{20}$)  
& $e(c_{30}$) & $e(c_{40}$) \\
\hline
\multicolumn{5}{|c|}{Type 1, with time} \\
\hline
30 (2)		& 67.69	& 49.62	& 34.35	& 23.66 \\
70 (6)		& 67.83	& 49.75	& 34.47	& 23.78 \\
110 (10)	& 64.16	& 48.15	& 33.99	& 23.80 \\
210 (20)	& 15.70	& 13.12	& 11.58	& 9.439 \\
410 (40)	& 9.781	& 5.065	& 3.848	& 2.939 \\ 
610 (60)	& 6.883	& 4.778	& 3.650	& 2.779 \\
810 (80)	& 6.831	& 4.736	& 3.619	& 2.754 \\
\hline
\multicolumn{5}{|c|}{Type 2,  without time} \\
\hline
30 (1)		& 89.85	& 83.50	& 73.50	& 65.06 \\
70 (3)		& 62.18	& 62.39	& 59.84	& 56.96 \\
110 (5)		& 51.05	& 48.08	& 42.41	& 36.61 \\
210 (10)	& 7.397	& 4.893	& 3.972	& 3.309 \\
410 (20)	& 6.750	& 4.405	& 3.397	& 2.656  \\ 
610 (30)	& 6.728	& 4.392	& 3.381	& 2.648 \\
810 (40)	& 6.726	& 4.392	& 3.381	& 2.649 \\
\hline
\multicolumn{5}{|c|}{Type 2, with time} \\
\hline
30 (1)		& 89.85	& 83.51	& 73.50	& 65.07 \\
70 (3)		& 62.18	& 62.40	& 59.85	& 56.97 \\
110 (5)		& 51.03	& 48.06	& 42.40	& 36.60 \\
210 (10)	& 7.397	& 4.892	& 3.971	& 3.308 \\
410 (20)	& 6.750	& 4.405	& 3.398	& 2.656 \\ 
610 (30)	& 6.728	& 4.392	& 3.381	& 2.648 \\
810 (40)	& 6.726	& 4.393	& 3.381	& 2.649 \\
\hline
\end{tabular}
\caption{Relative $L_2$ error for concentration at $t_m$ ($m=10, 20, 30$ and $40$). 
Dirichlet boundary conditions (Test 2).   
Left: Type 1 multiscale basis functions for concentration. 
Right: Type 2 multiscale basis functions for concentration. 
Reference solution with  $DOF^c_h = 52 050$}
\label{table-t2-dbc}
\end{table}

\begin{table}[h!]
\center
\begin{tabular}{ | c | c  | c | c | c | }
\hline
\multicolumn{5}{|c|}{without velocity} \\
\hline
$DOF^c_H (M^c)$    
& $e(c_{10}$) & $e(c_{20}$)  
& $e(c_{30}$) & $e(c_{40}$) \\
\hline
\multicolumn{5}{|c|}{Type 1, without time} \\
\hline
30 (2)		& 45.87	& 48.36	& 50.73	& 48.27 \\
70 (6)		& 5.842	& 6.302	& 7.917	& 11.05 \\
110 (10)	& 2.860	& 3.474	& 4.482	& 6.459 \\
210 (20)	& 2.131	& 2.703	& 3.338	& 4.570 \\
410 (40)	& 2.104	& 2.670	& 3.284	& 4.465 \\ 
610 (60)	& 2.102	& 2.667	& 3.279	& 4.453 \\
810 (80)	& 2.099	& 2.664	& 3.274	& 4.442 \\
\hline
\multicolumn{5}{|c|}{Type 2, without time} \\
\hline
30 (1)		& 75.71	& 89.67	& 93.60	& 93.30 \\
70 (3)		& 6.233	& 6.465	& 8.005	& 11.93 \\
110 (5)		& 4.621	& 6.042	& 7.537	& 10.75 \\
210 (10)	& 3.385	& 4.431	& 5.718	& 8.192 \\
410 (20)	& 3.164	& 4.037	& 5.214	& 7.593 \\ 
610 (30)	& 3.088	& 3.928	& 5.065	& 7.375 \\
810 (40)	& 3.056	& 3.896	& 5.021	& 7.301 \\
\hline
\multicolumn{5}{|c|}{Type 2, with time} \\
\hline
30 (1)		& 67.10	& 82.64	& 88.56	& 88.73 \\
70 (3)		& 4.373	& 4.409	& 5.511	& 8.065 \\
110 (5)		& 2.680	& 3.212	& 4.097	& 5.803 \\
210 (10)	& 2.143	& 2.723	& 3.338	& 4.492 \\
410 (20)	& 2.123	& 2.691	& 3.290	& 4.409 \\ 
610 (30)	& 2.121	& 2.689	& 3.287	& 4.405 \\
810 (40)	& 2.120	& 2.688	& 3.286	& 4.402 \\
\hline
\end{tabular}
\,\,
\begin{tabular}{ | c | c  | c | c | c | }
\hline
\multicolumn{5}{|c|}{with velocity} \\
\hline
$DOF^c_H (M^c)$  
& $e(c_{10}$) & $e(c_{20}$)  
& $e(c_{30}$) & $e(c_{40}$) \\
\hline
\multicolumn{5}{|c|}{Type 1, with time} \\
\hline
30 (2)		& 69.34	& 84.22	& 89.44	& 89.28 \\
70 (6)		& 49.09	& 61.96	& 67.20	& 64.77 \\
110 (10)	& 36.89	& 47.61	& 51.47	& 46.71 \\
210 (20)	& 7.022	& 11.53	& 13.35	& 16.50 \\
410 (40)	& 2.080	& 2.587	& 2.744	& 2.371 \\ 
610 (60)	& 2.071	& 2.565	& 2.713	& 2.336 \\
810 (80)	& 2.069	& 2.559	& 2.705	& 2.327 \\
\hline
\multicolumn{5}{|c|}{Type 2, without time} \\
\hline
30 (1)		& 47.37	& 56.09	& 61.62	& 59.57 \\
70 (3)		& 28.27	& 31.79	& 31.76	& 24.81 \\
110 (5)		& 15.85	& 18.89	& 18.45	& 15.28 \\
210 (10)	& 8.635	& 8.903	& 7.371	& 9.659 \\
410 (20)	& 2.126	& 2.669	& 2.823	& 2.283 \\ 
610 (30)	& 2.112	& 2.651	& 2.801	& 2.272 \\
810 (40)	& 2.109	& 2.646	& 2.797	& 2.270 \\
\hline
\multicolumn{5}{|c|}{Type 2, with time} \\
\hline
30 (1)		& 47.38	& 56.10	& 61.64	& 59.59 \\
70 (3)		& 28.26	& 31.78	& 31.76	& 24.80 \\
110 (5)		& 15.85	& 18.89	& 18.45	& 15.29 \\
210 (10)	& 8.636	& 8.902	& 7.368	& 9.661 \\
410 (20)	& 2.126	& 2.669	& 2.822	& 2.283 \\ 
610 (30)	& 2.112	& 2.651	& 2.801	& 2.272 \\
810 (40)	& 2.109	& 2.646	& 2.797	& 2.270 \\
\hline
\end{tabular}
\caption{Relative $L_2$ error for concentration at $t_m$ ($m=10, 20, 30$ and $40$). 
Neumann boundary conditions (Test 2).   
Left: Type 1 multiscale basis functions for concentration. 
Right: Type 2 multiscale basis functions for concentration. 
Reference solution with  $DOF^c_h = 52 050$}
\label{table-t2-nbc}
\end{table}

\begin{figure}[h!]
\centering
\includegraphics[width=0.3\textwidth]{bar01}
\begin{subfigure}{1.0\textwidth}
\centering
\includegraphics[width=1\textwidth]{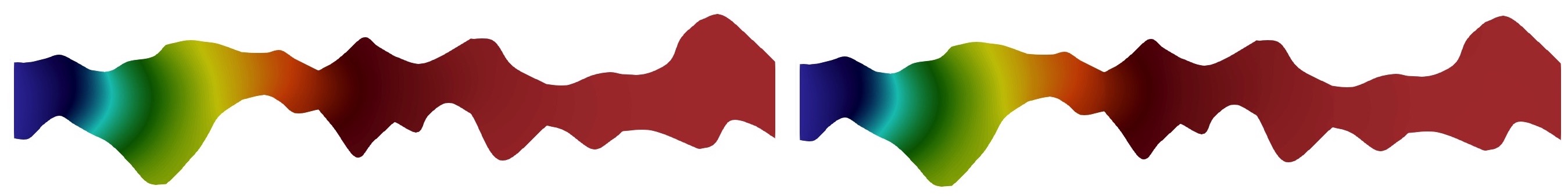}
\caption{Concentration, $c_{10}$. Left: reference solution. Right: multiscale solution}
\end{subfigure}
\begin{subfigure}{1.0\textwidth}
\centering
\includegraphics[width=1\textwidth]{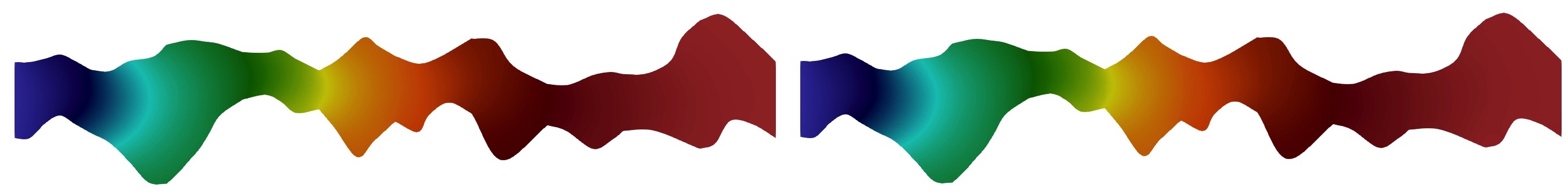}
\caption{Concentration, $c_{20}$. Left: reference solution. Right: multiscale solution}
\end{subfigure}
\begin{subfigure}{1.0\textwidth}
\centering
\includegraphics[width=1\textwidth]{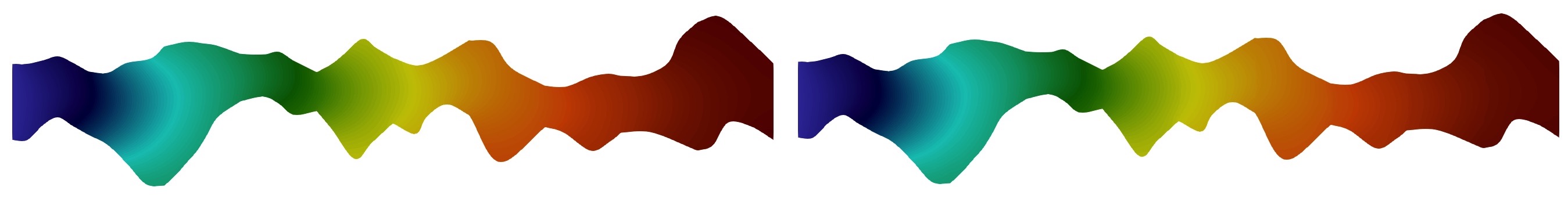}
\caption{Concentration, $c_{40}$. Left: reference solution. Right: multiscale solution}
\end{subfigure}
\caption{Relative $L_2$ error for concentration at $t_m$ ($m=10, 20, 30$ and $40$). 
Diffusion coefficient $D = 0.1$ (Test 2).   
Left: Type 1 multiscale basis functions for concentration. 
Right: Type 2 multiscale basis functions for concentration. 
Reference solution with  $DOF^c_h = 52 050$}
\label{sol-d0c1}
\end{figure}

\begin{figure}[h!]
\centering
\includegraphics[width=0.3\textwidth]{bar01}
\begin{subfigure}{1.0\textwidth}
\centering
\includegraphics[width=1\textwidth]{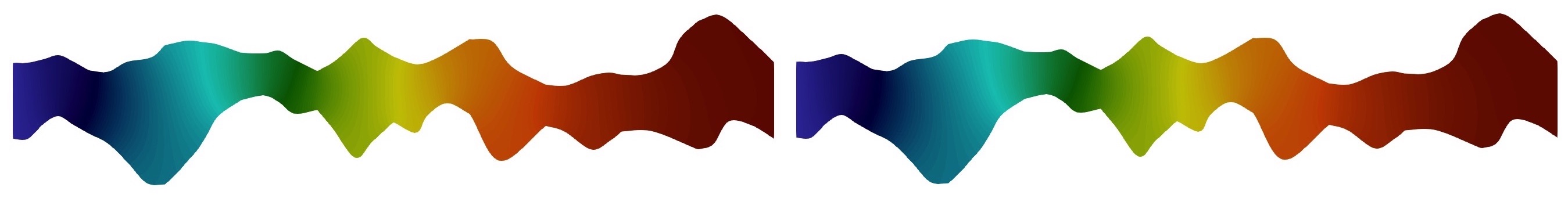}
\caption{Concentration, $c_{10}$. Left: reference solution. Right: multiscale solution}
\end{subfigure}
\begin{subfigure}{1.0\textwidth}
\centering
\includegraphics[width=1\textwidth]{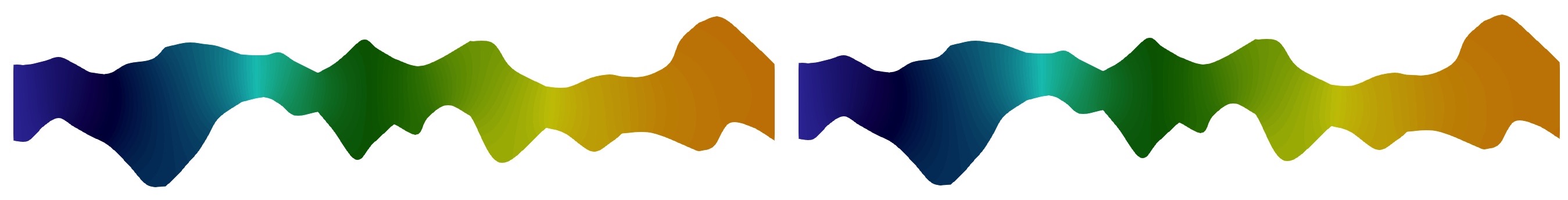}
\caption{Concentration, $c_{20}$. Left: reference solution. Right: multiscale solution}
\end{subfigure}
\begin{subfigure}{1.0\textwidth}
\centering
\includegraphics[width=1\textwidth]{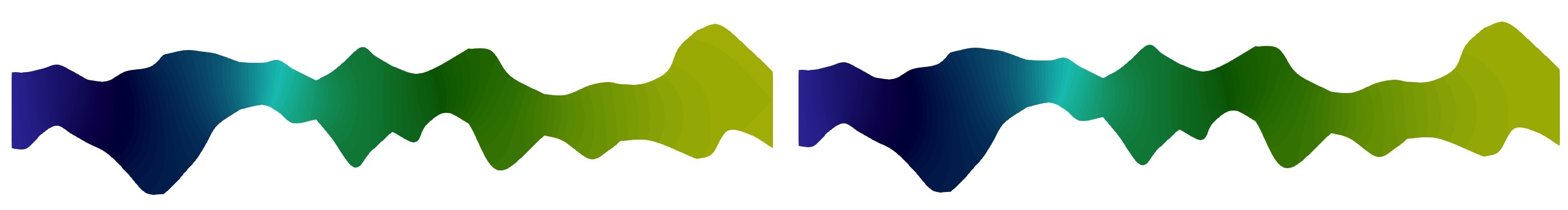}
\caption{Concentration, $c_{40}$. Left: reference solution. Right: multiscale solution}
\end{subfigure}
\caption{Relative $L_2$ error for concentration at $t_m$ ($m=10, 20, 30$ and $40$). 
Diffusion coefficient $D = 1$ (Test 2).   
Left: Type 1 multiscale basis functions for concentration. 
Right: Type 2 multiscale basis functions for concentration. 
Reference solution with  $DOF^c_h = 52 050$}
\label{sol-d1}
\end{figure}

\begin{table}[h!]
\center
\begin{tabular}{ | c | c | c | c | c |}
\hline
\multicolumn{5}{|c|}{without velocity} \\
\hline
$DOF^c_H (M^c)$  
& $e(c_{10}$) & $e(c_{20}$)  
& $e(c_{30}$) & $e(c_{40}$) \\
\hline
\multicolumn{5}{|c|}{Type 1, without time} \\
\hline
30 (2)		& 45.03	& 46.24	& 45.64	& 44.56 \\
70 (6)		& 3.870	& 3.959	& 3.910	& 3.813 \\
110 (10)	& 0.843	& 0.842	& 0.841	& 0.816 \\
210 (20)	& 0.505	& 0.502	& 0.505	& 0.492 \\
410 (40)	& 0.471	& 0.465	& 0.470	& 0.462 \\ 
610 (60)	& 0.469	& 0.462	& 0.468	& 0.460 \\
810 (80)	& 0.468	& 0.462	& 0.467	& 0.459 \\
\hline
\multicolumn{5}{|c|}{Type 2, without time} \\
\hline
30 (1)		& 93.96	& 93.31	& 92.55	& 92.02 \\
70 (3)		& 3.892	& 3.337	& 2.978	& 2.793 \\
110 (5)		& 0.486	& 0.488	& 0.503	& 0.509 \\
210 (10)	& 0.384	& 0.361	& 0.362	& 0.362 \\
410 (20)	& 0.380	& 0.357	& 0.356	& 0.357 \\ 
610 (30)	& 0.380	& 0.356	& 0.355	& 0.356 \\
810 (40)	& 0.380	& 0.356	& 0.356	& 0.356 \\
\hline
\multicolumn{5}{|c|}{Type 2, with time} \\
\hline
30 (1)		& 93.96	& 93.31	& 92.55	& 92.02 \\
70 (3)		& 3.893	& 3.339	& 2.980	& 2.795 \\
110 (5)		& 0.486	& 0.488	& 0.503	& 0.509 \\
210 (10)	& 0.384	& 0.361	& 0.362	& 0.363 \\
410 (20)	& 0.380	& 0.357	& 0.356	& 0.357 \\ 
610 (30)	& 0.380	& 0.356	& 0.355	& 0.356 \\
810 (40)	& 0.380	& 0.356	& 0.355	& 0.356 \\
\hline
\end{tabular}
\,\,
\begin{tabular}{ | c | c | c | c | c |}
\hline
\multicolumn{5}{|c|}{with velocity} \\
\hline
$DOF^c_H (M^c)$   
& $e(c_{10}$) & $e(c_{20}$)  
& $e(c_{30}$) & $e(c_{40}$) \\
\hline
\multicolumn{5}{|c|}{Type 1, with time} \\
\hline
30 (2)		& 48.32	& 51.01	& 50.58	& 49.42 \\
70 (6)		& 6.069	& 6.268	& 6.097	& 5.845 \\
110 (10)	& 0.827	& 0.875	& 0.870	& 0.826 \\
210 (20)	& 0.516	& 0.569	& 0.570	& 0.537 \\
410 (40)	& 0.475	& 0.524	& 0.528	& 0.502 \\ 
610 (60)	& 0.472	& 0.521	& 0.525	& 0.499 \\
810 (80)	& 0.472	& 0.520	& 0.524	& 0.498 \\
\hline
\multicolumn{5}{|c|}{Type 2, without time} \\
\hline
30 (1)		& 71.43	& 75.52	& 73.94	& 72.40 \\
70 (3)		& 8.661	& 8.862	& 8.827	& 8.631 \\
110 (5)		& 1.208	& 1.440	& 1.440	& 1.492 \\
210 (10)	& 0.509	& 0.619	& 0.656	& 0.648 \\
410 (20)	& 0.408	& 0.467	& 0.490	& 0.486 \\ 
610 (30)	& 0.400	& 0.454	& 0.475	& 0.472 \\
810 (40)	& 0.399	& 0.451	& 0.472	& 0.469 \\
\hline
\multicolumn{5}{|c|}{Type 2, with time} \\
\hline
30 (1)		& 71.43	& 75.52	& 73.94	& 72.40 \\
70 (3)		& 8.662	& 8.863 	& 8.828	& 8.632 \\
110 (5)		& 1.209	& 1.441	& 1.514	& 1.493 \\
210 (10)	& 0.509	& 0.619	& 0.656	& 0.649 \\
410 (20)	& 0.408	& 0.467	& 0.490	& 0.487 \\ 
610 (30)	& 0.400	& 0.454	& 0.475	& 0.472 \\
810 (40)	& 0.399	& 0.452	& 0.472	& 0.469 \\
\hline
\end{tabular}
\caption{Relative $L_2$ error for concentration at $t_m$ ($m=10, 20, 30$ and $40$). 
Diffusion coefficient $D = 0.1$ (Test 2).    
Left: Type 1 multiscale basis functions for concentration. 
Right: Type 2 multiscale basis functions for concentration. 
Reference solution with  $DOF^c_h = 52 050$}
\label{table-t2-d0c1}
\end{table}

\begin{table}[h!]
\center
\begin{tabular}{ | c | c | c | c |   c  | }
\hline
\multicolumn{5}{|c|}{without velocity} \\
\hline
$DOF_c^c (M)$  
& $e(c_{10}$) & $e(c_{20}$)  
& $e(c_{30}$) & $e(c_{40}$) \\\hline
\multicolumn{5}{|c|}{Type 1, without time} \\
\hline
30 (2)  		& 82.00 	& 80.87	& 78.91 	& 78.04 \\
70 (6)  		& 9.095 	& 10.93 	& 11.00 	& 10.74 \\
110 (10)  	& 2.086   	& 2.625 	& 2.636  	& 2.562 \\
210 (20) 	& 1.344   	& 1.707  	& 1.718 	& 1.671  \\
410 (40) 	& 1.149 	& 1.458  	& 1.468  	& 1.428 \\ 
610 (60) 	& 1.128 	& 1.431  	& 1.442   	& 1.402 \\
810 (80) 	& 1.120   & 1.421  	& 1.431   	& 1.392 \\
\hline
\multicolumn{5}{|c|}{Type 2, without time} \\
\hline
30 (1)  		& 98.76	& 98.43	& 98.25	& 98.17 \\
70  (3)   	& 7.119	& 7.728	& 7.684	& 7.536 \\
110 (5) 	& 0.699	& 0.872	& 0.870	& 0.846 \\
210 (10)	& 0.197	& 0.251	& 0.253	& 0.246 \\
410 (20)  	& 0.114	& 0.143	& 0.145	& 0.141 \\ 
610 (30)  	& 0.109	& 0.136	& 0.138	& 0.134 \\
810 (40)  	& 0.108	&  0.135	& 0.137	& 0.133 \\
\hline
\multicolumn{5}{|c|}{Type 2, with time} \\
\hline
30 (1)  		& 98.76	& 98.43	& 98.25	& 98.17 \\
70  (3)   	& 7.119	& 7.728	& 7.684	& 7.536 \\
110 (5) 	& 0.699	& 0.872	& 0.872	& 0.846 \\
210 (10)	& 0.197	& 0.251	& 0.251	& 0.246 \\
410 (20)  	& 0.114	& 0.143	& 0.145	& 0.141 \\ 
610 (30)  	& 0.109	& 0.136	& 0.138	& 0.134 \\
810 (40)  	& 0.108	 & 0.135	& 0.137	& 0.133 \\
\hline
\end{tabular}
\,\,
\begin{tabular}{ | c | c | c | c |   c  | }
\hline
\multicolumn{5}{|c|}{with velocity} \\
\hline
$DOF_c^c (M)$  
& $e(c_{10}$) & $e(c_{20}$)  
& $e(c_{30}$) & $e(c_{40}$) \\
\hline
\multicolumn{5}{|c|}{Type 1, with time} \\
\hline
30 (2)  		& 82.06	& 80.92	& 78.95	& 78.09 \\
70 (6)  		& 9.142	& 10.99	& 11.05	& 10.79 \\
110 (10)  	& 2.094	& 2.636	& 2.647	& 2.572 \\
210 (20) 	& 1.354	& 1.720	& 1.731	& 1.683 \\
410 (40) 	& 1.159	& 1.471	& 1.481	& 1.441 \\ 
610 (60) 	& 1.138	& 1.444	& 1.455	& 1.415 \\
810 (80) 	& 1.130	& 1.434	& 1.444	& 1.405 \\
\hline
\multicolumn{5}{|c|}{Type 2, without time} \\
\hline
30 (1)  		& 95.75	& 94.61	& 94.00	& 93.75 \\
70  (3)   	& 5.556	& 6.164	& 6.173	& 6.059 \\
110 (5) 	& 0.765	& 0.965	& 0.988	& 0.973 \\
210 (10)	& 0.299	& 0.387	& 0.398	& 0.391 \\
410 (20)  	& 0.155	& 0.204	& 0.210	& 0.205 \\ 
610 (30)  	& 0.144	& 0.189	& 0.194	& 0.190 \\
810 (40)  	& 0.143	& 0.187	& 0.192	& 0.188 \\
\hline
\multicolumn{5}{|c|}{Type 2, with time} \\
\hline
30 (1)  		& 95.75 	& 94.61 	& 94.00	& 93.75 \\
70  (3)   	& 5.556 	& 6.164 	& 6.173	& 6.059 \\
110 (5) 	& 0.765 	& 0.965 	& 0.988	& 0.974 \\
210 (10)	& 0.299 	& 0.387 	& 0.398	& 0.391 \\
410 (20)  	& 0.155 	& 0.204 	& 0.210	& 0.205 \\ 
610 (30)  	& 0.144 	& 0.189 	& 0.194	& 0.190 \\
810 (40)  	& 0.143 	& 0.187 	& 0.192	& 0.188 \\
\hline
\end{tabular}
\caption{Relative $L_2$ error for concentration at $t_m$ ($m=10, 20, 30$ and $40$). 
Diffusion coefficient $D = 1$ (Test 2).    
Left: Type 1 multiscale basis functions for concentration. 
Right: Type 2 multiscale basis functions for concentration. 
Reference solution with  $DOF^c_h = 52 050$}
\label{table-t2-d1}
\end{table}


To investigate the influence of the boundary conditions on the results of the presented multiscale method, we consider Geometry 1 and set diffusion coefficient $D = 0.01$.  
We consider following types of boundary conditions on $\Gamma_w$:
\begin{itemize}
\item[] \textit{DBC} - Dirichlet type boundary conditions:
\[
c = c_w, \quad x \in \Gamma_w, 
\]
with $c_w = 1$, 
initial condition $c_0 = 0$ and $c_{in} = 1$ on $\Gamma_{in}$. 
We perform simulations for $t_{max} = 0.1$  with 40 time iterations. 


\item[] \textit{NBC} - Neumann type boundary conditions: 
\[
- D \nabla c  \cdot n = \beta, \quad x \in \Gamma_w, 
\]
with $\beta = 0.01$, 
initial condition $c_0 = 1$ and $c_{in} = 0$ on $\Gamma_{in}$. 
We perform simulations for $t_{max} = 0.5$ with 40 time iterations. 
\end{itemize}
The coarse grid is structured with 10  local domains. 

To investigate the presented method for different diffusion coefficients, we consider transport problem in Geometry 1 with Robin boundary conditions (RBC) with $c_w = 1$ and  $\alpha = 0.1$, 
We set initial condition $c_0 = 1$ and $c_{in} = 0$ on $\Gamma_{in}$. 
Simulations are performed for $t_{max} = 0.7$ with 40 time iterations. 
We  consider $D = 0.1$ and $D = 1$. 
 
For each type of boundary condition and value of diffusion coefficient, we investigate the influence of the types of multiscale basis functions in detail. 
We consider \textit{Type 1} and \textit{Type 2} multiscale basis functions and present the results for multiscale basis functions with and without time and velocity in the construction of the snapshot space. Note that, spectral problems are similar for all types of snapshots and contain only information about diffusion.

In Figures \ref{sol-dbc} and \ref{sol-nbc}, we present results for Dirichlet and Neumann boundary conditions.  
Concentration distributions are depicted for the reference (fine scale) and multiscale solutions at different time layers $t_m$ for $m = 10, 20$ and $40$. 
In calculations, we used a fixed number of multiscale basis functions for the velocity field ($M^u = 20$ of Type 2 with time and velocity).  
We observe good results of the presented method for solving transport problems for all types of nonhomogeneous boundary conditions. 
  
In Tables \ref{table-t2-dbc}  and  \ref{table-t2-nbc}, we present relative errors for concentration in \%  for different number of multiscale basis functions for a fixed number of multiscale basis functions for velocity field ($M^u = 20$ of Type 2). 
Relative errors for concentrations are presented for four time layers $t_m$ with $m = 10, 20, 30$ and $40$. 
For a test with Dirichlet boundary conditions and snapshots with velocity information, we obtain near 2 \% of concentration error at the final time, when we take 20  multiscale basis functions of Type  2. For Type 1  basis functions, we obtain similar results for the same size of the coarse grid system ($DOF^c_H = 410$).  Numerical results are almost similar for any type of snapshot space due to large influence of the boundary conditions to the transport problem solution (see Figure {\ref{sol-dbc}). 
For Neumann boundary conditions, we observe that the multiscale basis functions are better when we take into account velocity into the local problems. For example, we obtain  4-7 \% of errors for basis functions without velocity and reduce errors to 2 \% for basis with velocity information. However, adding time information into basis construction does not affect the errors for tests with any typed of boundary conditions. 
We observe that Type 1 and 2 basis functions provide similar results for the 2D transport problem. Note that, for the homogeneous boundary conditions it is not needed to add second space for handling boundary conditions in Type 2 basis functions and the size of the coarse grid system will be two times smaller.

In Figures \ref{sol-d0c1} and \ref{sol-d1}, we present results for $D = 0.1$ and $1$ with Robin boundary conditions on the wall boundary.  
Concentration distributions are depicted for the reference (fine scale) and multiscale solutions at different time layers $t_m$ for $m = 10, 20$ and $40$. 
We used a fixed number of multiscale basis functions for the velocity field ($M^u = 20$ of Type 2 with time and velocity).

Relative errors for concentration for different types of multiscale basis functions are presented in Tables \ref{table-t2-d0c1}  and  \ref{table-t2-d1}. For velocity field we used $M^u = 20$ of Type 2 multiscale basis functions in all calculations.   
Concentration errors are presented for four time layers $t_m$ with $m = 10, 20, 30$ and $40$. 
We observe smaller errors for larger values of diffusion coefficients. For example, we have $0.4$  \% of errors for  $M^u = 20$ of Type 2 multiscale basis functions for $D = 0.1$ and $0.2$  \% of errors for $D = 0.1$. 
For larger diffusion, coefficients errors reduce faster with an increasing number of basis functions, and we can use 5 multiscale basis functions of Type 2 for obtaining results with less than one percent of errors. Numerical results show that the results are better for basis functions construction without time and velocity. However, they are almost the same and in general, it is better to construct a basis using all information in the snapshot space constructions.

For Test 2 with different boundary conditions and diffusion coefficients, we have the following conclusions for the transport problem:
\begin{itemize}
\item Multiscale basis functions are better to construct with all information (time and velocity). 
\item  Type 1 and 2 provide almost the same results with the same size of the coarse grid system for the 2D problem. However, Type 2 is preferable due to the exact definition of the coarse grid parameters and can be more usable for the general case with homogeneous boundary conditions.  
\item The presented multiscale method provides good results with small errors and gives a huge reduction of the system size. 
\end{itemize}

\subsection{Test 3  (unstructured coarse grid)}

Finally, we consider the test for unstructured coarse grids and investigate the influence of the velocity accuracy to the concentration errors for both structured and unstructured coarse grids.

\begin{figure}[h!]
\centering
\includegraphics[width=0.75\textwidth]{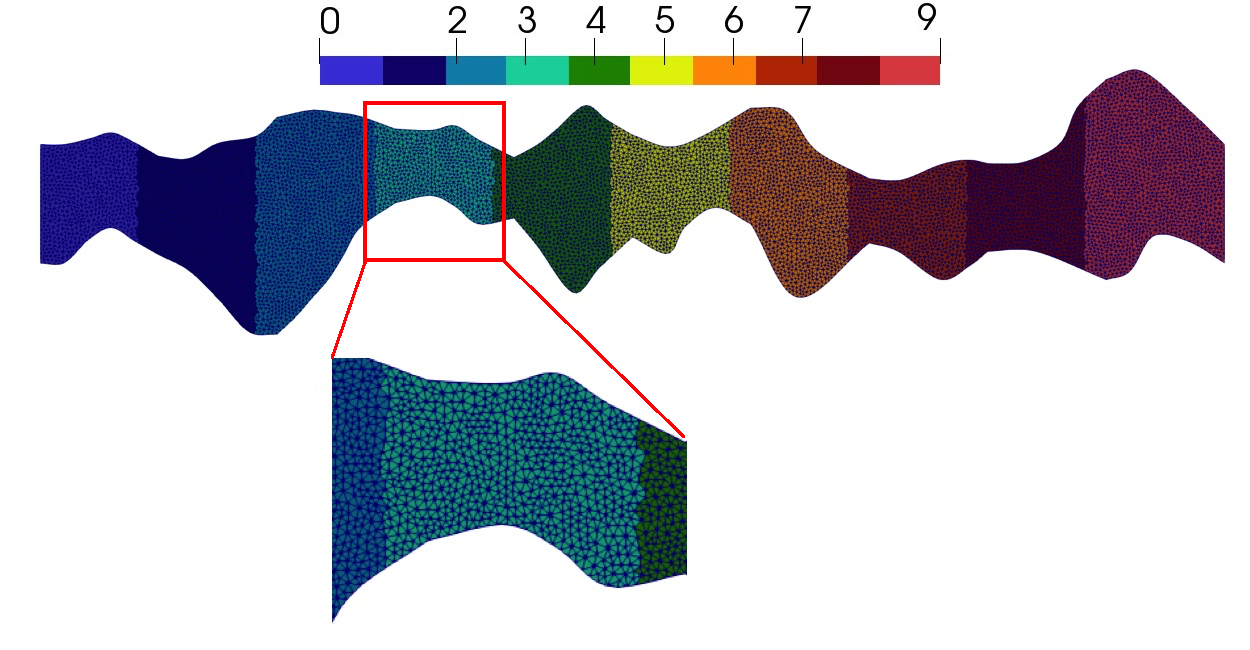}
\caption{Fine grid with subdomain markers for Geometry 1 with 10 local domains}
\label{geom2d-uns}
\end{figure}

\begin{figure}[h!]
\begin{subfigure}{1.0\textwidth}
\centering
\includegraphics[width=1\textwidth]{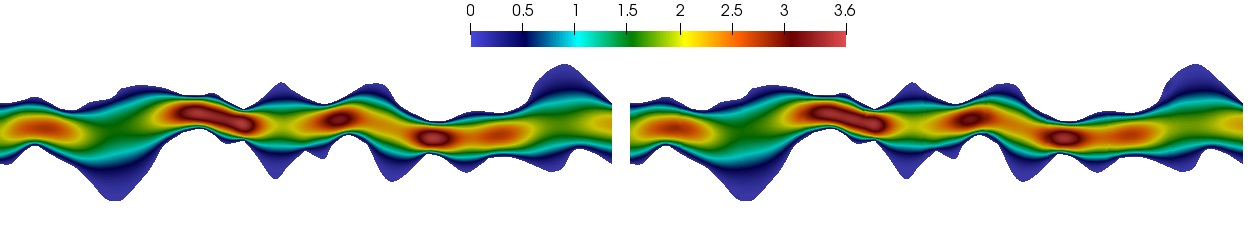}
\caption{Velocity at final time. Left: reference solution. Right: multiscale solution}
\end{subfigure}
\begin{subfigure}{1.0\textwidth}
\centering
\includegraphics[width=0.33\textwidth]{bar01}\\
\includegraphics[width=1\textwidth]{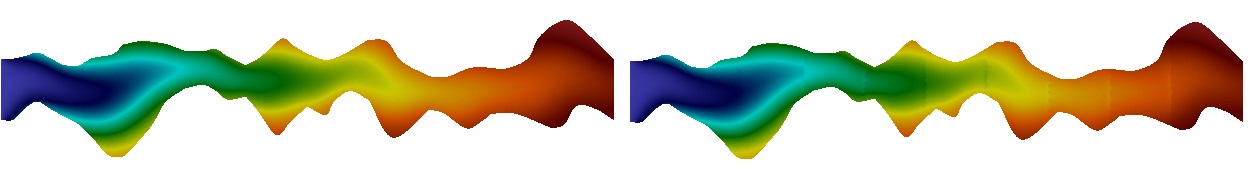}
\caption{Concentration at final time. Left: reference solution. Right: multiscale solution}
\end{subfigure}
\caption{Reference and multiscale solutions of velocity (first row)  and concentration (second row) at the final time.  
Unstructured coarse grid solution (Test 3). 
Left: reference solution with  
$DOF^u_h = 121 226$ and $DOF^c_h = 51954$. 
Right: multiscale solution with multiscale basis functions of Type 2,  
$DOF^u_H = 810$ and $DOF^c_H = 610$. }
\label{sol-uns}
\end{figure}

\begin{table}[h!]
\center
\begin{tabular}{|c|c|}
\hline
\multicolumn{2}{|c|}{Structured}  \\ 
\hline
$DOF_H^u (M^u)$  & $e(u)$ (\%) \\
\hline
110 (5)		& 19.856 \\
210 (10)	& 10.529 \\
410 (20)	& 3.230 \\
610 (30)	& 1.756 \\
810 (40)	& 1.346 \\
\hline
\end{tabular}
\,\,\,\,
\begin{tabular}{|c|c|}
\hline 
\multicolumn{2}{|c|}{Unstructured}  \\ 
\hline
$DOF_H^u (M^u)$  & $e(u)$ (\%) \\
\hline
110 (5)		& 24.824 \\
210 (10)	& 16.777 \\
410 (20)	& 7.981 \\
610 (30)	& 5.986 \\
810 (40)	& 3.846 \\
\hline
\end{tabular}
\caption{Relative $L_2$ error for the velocity at the final time. 
Unstructured coarse grid solution (Test 3). 
Type 2 multiscale basis functions for velocity. 
$DOF^u_h = 121 226$. 
}
\label{err-uns-u}
\end{table}

\begin{table}[h!]
\center
\begin{tabular}{ | c | c | c | c |   c  | }
\hline
\multicolumn{5}{|c|}{Structured}  \\ 
\hline
$DOF_c^c (M^c)$  
& $e(c_{10}$) & $e(c_{20}$)  
& $e(c_{30}$) & $e(c_{40}$) \\\hline
\multicolumn{5}{|c|}{u with 5 basis} \\
\hline
30 (1)  		& 49.86 	& 53.64	& 52.60 	& 50.31 \\
70 (3)  		& 26.77 	& 24.59	& 21.24 	& 19.56 \\
110 (5)  	& 15.25	& 13.78	& 10.12 	& 7.744  \\
210 (10) 	& 8.775 	& 7.873	& 6.543 	& 5.582 \\
410 (20) 	& 5.185 	& 5.501	& 5.294 	& 4.947 \\
610 (30) 	& 5.185	& 5.495	& 5.290 	& 4.947 \\
810 (40) 	& 5.184	& 5.493	& 5.289 	& 4.947 \\
\hline
\multicolumn{5}{|c|}{u with 10 basis} \\
\hline
30 (1)  		& 49.94 	& 53.73	& 52.69 	&  50.41 \\
70 (3)  		& 26.42	& 24.39	& 21.14 	&  19.49 \\
110 (5)  	& 14.87 	& 13.39	& 9.707 	&  7.362 \\
210 (10) 	& 7.543 	& 6.480	& 4.788 	&  3.503 \\
410 (20) 	& 3.061 	& 3.319	& 3.203 	&  2.937 \\
610 (30) 	& 3.056 	& 3.309	& 3.193 	&  2.931 \\
810 (40) 	& 3.054 	& 3.306	& 3.190 	&  2.929 \\
\hline
\multicolumn{5}{|c|}{u with 20 basis} \\
\hline
30 (1)  		& 49.91	& 53.68 	& 52.64 	& 50.36 \\
70 (3)  		& 26.22	& 24.27	& 21.11 	& 19.47 \\
110 (5)  	& 14.63	& 13.24	& 9.693 	& 7.351 \\
210 (10) 	& 6.909	& 5.760	& 3.943 	& 2.372 \\
410 (20) 	& 1.867	& 2.003	& 1.870 	& 1.701 \\
610 (30) 	& 1.850	& 1.978	& 1.842 	& 1.674 \\
810 (40) 	& 1.846	& 1.972	& 1.836 	& 1.669 \\
\hline
\multicolumn{5}{|c|}{u with 40 basis} \\
\hline
30 (1)  		& 49.89 	& 53.64	& 52.60 	& 50.32 \\
70 (3)  		& 26.18	& 24.23 	& 21.09 	& 19.44 \\
110 (5)  	& 14.60	& 13.21	& 9.708 	& 7.396 \\
210 (10) 	& 6.812	& 5.646	& 3.806 	& 2.237 \\
410 (20) 	& 1.689	& 1.820 	& 1.701 	& 1.600 \\
610 (30) 	& 1.669	& 1.790 	& 1.667 	& 1.569 \\
810 (40) 	& 1.664	& 1.783	& 1.660 	& 1.562 \\
\hline
\multicolumn{5}{|c|}{u fine} \\
\hline
30 (1)  		& 49.89	& 53.64	& 52.60 	& 50.32 \\
70 (3)  		& 26.17	& 24.22	& 21.07 	& 19.42 \\
110 (5)  	& 14.58	& 13.17	& 9.620 	& 7.312 \\
210 (10) 	& 6.814	& 5.621	& 3.714 	& 2.068 \\
410 (20) 	& 1.538	& 1.647	& 1.437 	& 1.307 \\
610 (30) 	& 1.517	& 1.617	& 1.402 	& 1.275 \\
810 (40) 	& 1.512	& 1.611	& 1.395 	& 1.269 \\
\hline
\end{tabular}
\,\,
\begin{tabular}{ | c | c | c | c |   c  | }
\hline
\multicolumn{5}{|c|}{Unstructured}  \\ 
\hline
$DOF_c^c (M^c)$  
& $e(c_{10}$) & $e(c_{20}$)  
& $e(c_{30}$) & $e(c_{40}$) \\\hline
\multicolumn{5}{|c|}{u with 5 basis} \\
\hline
30 (1)  		& 50.856	& 56.271	& 55.245  & 52.700  \\
70 (3)  		& 22.335	& 20.899	& 17.647  & 15.467 \\
110 (5)  	& 16.503	& 14.952	& 11.308  & 8.805\\
210 (10) 	& 13.910	& 13.089	& 10.954  & 9.543 \\
410 (20) 	& 11.157	& 10.444	& 9.704   & 9.628 \\
610 (30) 	& 11.627 	& 10.835	& 10.103  & 10.070 \\
810 (40) 	& 11.947	& 11.174	& 10.448  & 10.412 \\
\hline
\multicolumn{5}{|c|}{u with 10 basis} \\
\hline
30 (1)  		& 50.920 & 56.299	& 55.222 & 52.667 \\
70 (3)  		& 21.705 & 19.859	& 16.463 & 14.420 \\
110 (5)  	& 15.006 & 13.277	& 9.427 & 6.236 \\
210 (10) 	& 10.891 & 9.877	& 7.186 & 4.499 \\
410 (20) 	& 3.275  & 3.438	& 3.026 & 2.431 \\
610 (30) 	& 3.132  & 3.334	& 3.004 & 2.515 \\
810 (40) 	& 3.103 & 3.346	& 3.070 & 2.616 \\
\hline
\multicolumn{5}{|c|}{u with 20 basis} \\
\hline
30 (1)  		& 50.997	& 56.365	& 55.282 	& 52.732  \\
70 (3)  		& 21.744	& 20.012	& 16.703 	& 14.762 \\
110 (5)  	& 14.917	& 13.181	& 9.2779 	& 6.281  \\
210 (10) 	& 10.609	& 9.420	& 6.659  	& 4.097\\
410 (20) 	& 2.186   	& 2.284 	& 1.889 	& 1.438 \\
610 (30) 	& 2.004	& 2.159	& 1.833 	& 1.485 \\
810 (40) 	& 1.895	& 2.075	& 1.816 	& 1.517 \\
\hline
\multicolumn{5}{|c|}{u with 40 basis} \\
\hline
30 (1)  		& 50.999	& 56.357	& 55.267  	& 52.717  \\
70 (3)  		&  21.735 	& 21.734	& 16.768	& 14.822  \\
110 (5)  	& 14.897 	& 13.165	& 9.307 	& 6.338 \\
210 (10) 	& 10.552 & 9.308	& 6.583 	& 4.088 \\
410 (20) 	& 2.054	& 2.133	& 1.775 	& 1.424 \\
610 (30) 	& 1.850 	&1.998	&  1.7076 	& 1.452 \\
810 (40) 	& 1.738	& 1.919 	& 1.688 	& 1.482 \\
\hline
\multicolumn{5}{|c|}{u fine} \\
\hline
30 (1)  		& 50.968 	& 50.968	& 55.240 &  52.692 \\
70 (3)  		& 21.805	& 20.080	& 16.763 &  14.858 \\
110 (5)	  	& 15.055	& 13.196 	& 9.225  &  6.391 \\
210 (10) 	& 10.6257	& 9.221	& 6.366 &  3.951 \\
410 (20) 	& 1.844 	& 1.737	& 1.342 &  1.058 \\
610 (30) 	& 1.575	& 1.549	& 1.247  &  1.074 \\
810 (40) 	& 1.420	& 1.433	& 1.215 &  1.100 \\
\hline
\end{tabular}
\caption{Relative $L_2$ error for concentration at $t_m$ ($m=10, 20, 30$ and $40$). 
Unstructured coarse grid solution (Test 3). 
Type 2 multiscale basis functions for velocity. 
$DOF^u_h = 51 954$. 
}
\label{err-uns-c}
\end{table}

We consider a test problem with non-homogeneous Robin boundary condition for concentration with  $c_w = 1$, $\alpha = 0.01$ and $D = 0.01$. We set $c_0 = 1$ and $u_0 = 0$ as initial conditions.  
For the inflow (left) boundary, we set $c_{in} = 0$ for $\Gamma_{in}$. 
We perform simulations for $t_{max} = 0.7$ with 40 time iterations. 
We consider a structured and unstructured coarse grid with 10 local domains. The structured coarse grid is similar to Test 1.  The unstructured grid is depicted in Figure \ref{geom2d-uns}, where we show local domain markers. 

In Figure \ref{sol-uns}, we present reference and multiscale solutions for an unstructured grid. We depicted the magnitude of the velocity field and concentration at the final time.  In multiscale solver, we used $M^u = 40$ multiscale basis functions for velocity and $M^u = 20$ multiscale basis functions for concentration. 
For reference solution, we have $DOF^u_h = 121 226$ and $DOF^c_h = 51954$. For multiscale solution, we have $DOF^u_H = 810$ and $DOF^c_H = 610$.
We observe a good accuracy of the presented method on an unstructured grid. 

In Table \ref{err-uns-u}, we present relative errors in \% for velocity between reference solution and multiscale solution with different numbers of the multiscale basis functions at the final times. Results are presented for structured and unstructured grids. We observe good results for the unstructured grid, however, velocity errors are smaller on a structured grid. 
For example for 40 multiscale basis functions, we have 
$1.3$ \% of error for a structured grid and 
$3.8$ \% of error for an unstructured grid. 

In Table \ref{err-uns-c}, we investigate the influence of the velocity accuracy on the concentration solution for structured and unstructured coarse grids.  
We present numerical results  for different number of multiscale basis functions for concentration of Type 2. Relative errors for concentrations are presented for four time layers $t_m$ with $m = 10, 20, 30$ and $40$. 
We can obtain good solution for concentration when we take at least 10 multiscale basis functions for velocity. 
For example, we obtain near 10 \% of concentration error ($M^c=20$), when we use multiscale velocity solution with 5 multiscale basis functions and near 1 \%  using velocity with $M^u = 20$. 
For concentration, we obtain similar results for structured and unstructured coarse grids using a multiscale solution of velocity with a sufficient number of basis functions.

\section{Conclusions}

We developed a multiscale model order reduction technique for the solution of the flow and transport problem in thin domains. Our motivation stems from reducing the problem dimension in thin layer applications. 
For the fine grid approximation, we apply the discontinuous Galerkin method and use the solution as a reference solution. 
Our multiscale approach for solving problems in complex thin geometries gives an accurate approximation of the velocity space and transport processes. 
We presented two types of the local multiscale basis functions for velocity and concentration, where the first is based on the approach combining all possible flows and transports directions in the local domain and the second approach is based on the separation of the macroscale parameters by the flow direction for velocity and by boundary type for transport. 
Moreover, our multiscale spaces can accurately capture complex processes on the rough wall boundaries with non-homogeneous boundary conditions.  
We use numerical simulations for three test geometries for two and three-dimensional problems to demonstrate the performance of our method.
Numerical investigation of the presented method was performed (1) for different geometries of the computational domain, (2) for different boundary conditions and diffusion coefficients, and (3) for the unstructured coarse grid.  
The proposed multiscale method provides good results with small errors and gives a huge reduction of the system size.

%

\bibliographystyle{unsrt}
\bibliography{lit}

\end{document}